\documentclass[10pt]{article}
\usepackage[utf8]{inputenc}
\usepackage{lmodern}
\usepackage[T1]{fontenc}
\usepackage{amsmath}
\usepackage{amsfonts}
\usepackage{amssymb}
\usepackage{amsthm}
\usepackage{mathrsfs}
\usepackage{stmaryrd}
\usepackage{centernot} 
\usepackage{mathtools}
\usepackage{graphicx}
\usepackage[left=2.9cm,right=2.9cm,top=3cm,bottom=3cm]{geometry}
\usepackage{fancyhdr}
\usepackage{enumitem}
\usepackage{hyperref}
\usepackage{tikz}
\usepackage{dsfont}
\usepackage{multicol}
\usepackage{xcolor}

\usepackage[maxbibnames=99]{biblatex}
\renewcommand{\Im}{\operatorname{Im}}
\DeclareMathOperator{\Ker}{Ker}
\renewcommand{\Re}{\operatorname{Re}}
\renewcommand{\epsilon}{\varepsilon}
\DeclareMathOperator{\supp}{supp}
\DeclareMathOperator{\dist}{d}

\DeclareMathOperator{\Div}{div}

\DeclareMathOperator{\NL}{NL}
\DeclareMathOperator{\Lin}{L} 

\newcommand\blfootnote[1]{
  \begingroup
  \renewcommand\thefootnote{}\footnote{#1}
  \addtocounter{footnote}{-1}
  \endgroup
}

\theoremstyle{plain}
\newtheorem{thm}{Theorem}
\newtheorem{lem}[thm]{Lemma}
\newtheorem{prop}[thm]{Proposition}
\newtheorem{cor}[thm]{Corollary}

\theoremstyle{definition}
\newtheorem{defn}[thm]{Definition}

\theoremstyle{remark}
\newtheorem{rem}[thm]{Remark}

\begin{document}

\begin{center}
\large{\textbf{\uppercase{Local controllability around a regular solution and null-controllability of scattering solutions for semilinear wave equations}}}

\vspace{\baselineskip}

\large{\textsc{Thomas Perrin
\footnote{\textit{Université de Rennes, ENS Rennes, INRIA, CNRS, IRMAR - UMR 6625, F-35000 Rennes, France.} The majority of this work was carried out at \textit{Laboratoire Analyse Géométrie et Application, Institut Galilée - UMR 7539, CNRS/Université Sorbonne Paris Nord, 99 avenue J.B. Clément, 93430 Villetaneuse, France. This is the author’s accepted manuscript. The final version of this article will appear in \emph{Journal of Differential Equations} (Elsevier).}}}
\blfootnote{\textit{Keywords:} local controllability, exact controllability, nonlinear wave equation, nonlinear Klein-Gordon equation, blow-up, scattering solutions.} \blfootnote{\textit{MSC2020:} 35B40, 35B44, 35L71, 93D20. }}

\vspace{\baselineskip}
\end{center}

\noindent
\textbf{Abstract.} On a Riemannian manifold with or without boundary, and whether bounded or unbounded, we consider a semilinear wave (or Klein-Gordon) equation with a subcritical nonlinearity (either defocusing or focusing). We establish local controllability around a partially analytic solution, under the Geometric Control Condition. Specifically, some blow-up solutions can be controlled. In the case of a Klein-Gordon equation on a non-trapping exterior domain of small dimension, we prove the null-controllability of scattering solutions. The proof is based on local energy decay and global-in-time Strichartz estimates. Several consequences are presented, including the null-controllability of a solution initiated near the ground state in some focusing cases, and exact controllability in some defocusing cases.

\section*{Introduction}

\paragraph{Setting and main results.} 
Let $d \geq 2$. Let $\Omega$ be the interior of a smooth $d$-dimensional Riemannian manifold of metric $\mathsf{g}$, with or without boundary, which is either a compact Riemannian manifold, or a compact perturbation of $\mathbb{R}^d$, that is, the complement in $\mathbb{R}^d$ of a smooth bounded (possibly empty) open set, with a metric equal to the Euclidean one outside a ball. In short, we write $\partial \Omega = \emptyset$ if $\Omega$ is either $\mathbb{R}^d$ or a compact Riemannian manifold without boundary, and $\partial \Omega \neq \emptyset$ if $\Omega$ is either a compact perturbation of $\mathbb{R}^d$ (with $\Omega \neq \mathbb{R}^d$) or a compact Riemannian manifold with nonempty boundary. If $\partial \Omega \neq \emptyset$, then we denote by $\partial \Omega$ the boundary of $\Omega$, and we write $\overline{\Omega} = \Omega \cup \partial \Omega$, and if $\partial \Omega = \emptyset$, then we write $\overline{\Omega} = \Omega$. In addition, we say that $\Omega$ is unbounded if $\Omega$ is a compact perturbation of $\mathbb{R}^d$ (or if $\Omega = \mathbb{R}^d$).

Write $H_0^1(\Omega)$ for the closure of $\mathscr{C}_\mathrm{c}^\infty(\Omega)$ in $H^1(\Omega)$, with respect to the norm $\int_\Omega \left(\left\vert \nabla u \right\vert^2 + \vert u \vert^2 \right) \mathrm{d}x$. Let $\beta \in \mathbb{R}$ be such that the Poincaré inequality
\begin{equation*}
    \int_\Omega \left(\left\vert \nabla u \right\vert^2 + \beta \vert u \vert^2 \right) \mathrm{d}x \gtrsim \int_\Omega \vert u \vert^2 \mathrm{d}x, \quad u \in H_0^1(\Omega),
\end{equation*}
is satisfied. This specifically requires $\beta > 0$ if $\partial \Omega = \emptyset$ or if $\Omega$ is unbounded. For $u \in H_0^1(\Omega)$, we write $\left\Vert u \right\Vert_{H_0^1(\Omega)}^2 = \int_\Omega \left(\left\vert \nabla u \right\vert^2 + \beta \vert u \vert^2 \right) \mathrm{d}x$. If $\Omega$ is bounded, $H_0^1(\Omega)$ is also the the closure of $\mathscr{C}_\mathrm{c}^\infty(\Omega)$ in $H^1(\Omega)$ with respect to the norm $\int_\Omega \left\vert \nabla u \right\vert^2 \mathrm{d}x$. If $\Omega$ is an exterior domain, there is no standard definition of $H_0^1(\Omega)$: some authors (for example, \cite{Melrose79}) define it using the norm $\int_\Omega \left\vert \nabla u \right\vert^2 \mathrm{d}x$, whereas others (for example, \cite{BURQGerardTze}) define it using the norm $\int_\Omega \left(\left\vert \nabla u \right\vert^2 + \vert u \vert^2 \right) \mathrm{d}x$.

This articles contains a local controllability and a null-controllability result. We consider a power-like nonlinearity $f \in \mathscr{C}^2(\mathbb{R}, \mathbb{R})$ satisfying $f(0) = f^\prime(0) = 0$, and the following assumptions. For the local controllability result, we assume that there exist $C_0 > 0$ and $\alpha$ such that
\begin{equation}\label{eq_def_nonlinearity_local}
\begin{array}{c}
\left\vert f^{\prime \prime}(s) \right\vert \leq C_0 \left( 1 + \vert s \vert \right)^{\alpha - 2} \text{ for all } s \in \mathbb{R}, \\
1 < \alpha < \frac{d + 2}{d - 2}, \quad \text{ and } \quad \left( \alpha \leq \frac{d - 1}{d - 3} \text{ if } d \geq 5 \text{ and } \partial \Omega \neq \emptyset \right). 
\end{array}
\end{equation}
For the null-controllability result, we assume that $\Omega$ is unbounded, with $3 \leq d \leq 5$, and that there exist $C_0 > 0$ and $\alpha_0 \leq \alpha_1$ such that 
\begin{equation}\label{eq_def_nonlinearity_scatt}
\begin{array}{c}
\left\vert f^{\prime \prime}(s) \right\vert \leq C_0 \left( \vert s \vert^{\alpha_0 - 2} + \vert s \vert ^{\alpha_1 - 2} \right) \text{ for all } s \in \mathbb{R}, \\
2 < \alpha_0 \leq \alpha_1 < \frac{d + 2}{d - 2}, \quad \text{ and } \quad \left( d \neq 5 \text{ if } \partial \Omega \neq \emptyset\right).
\end{array}
\end{equation}
Note that $\alpha$, $\alpha_0$ and $\alpha_1$ can be arbitrarily large if $d = 2$, and that (\ref{eq_def_nonlinearity_scatt}) implies $f^{\prime \prime}(0) = 0$. Note also that the existence of $(\alpha_0, \alpha_1)$ satisfying (\ref{eq_def_nonlinearity_scatt}) implies $2 < \frac{d + 2}{d - 2}$, and thus $d \leq 5$. 

A typical example of such a nonlinearity $f$ is given by
\[f(s) = \sum_{j = 0}^n \lambda_j s \vert s \vert^{\alpha_j - 1}, \quad s \in \mathbb{R},\]
with $n \in \mathbb{N}$, $\lambda_0, \cdots, \lambda_n \in \mathbb{R}$, and $2 \leq \alpha_0 \leq \cdots \leq \alpha_n$ such that (\ref{eq_def_nonlinearity_local}) or (\ref{eq_def_nonlinearity_scatt}) holds. 

\begin{rem}
If $f$ satisfies (\ref{eq_def_nonlinearity_scatt}) for some $\alpha_0 \leq \alpha_1$, then it also satisfies (\ref{eq_def_nonlinearity_local}) for $\alpha = \alpha_1$. Hence, any result stated with $f$ satisfying (\ref{eq_def_nonlinearity_scatt}) can be applied with $f$ satisfying (\ref{eq_def_nonlinearity_local}).
\end{rem}

Write $\Delta_{\mathsf{g}}$ for the Laplace-Beltrami operator on $\Omega$, and $\square = \partial_t^2 - \Delta_{\mathsf{g}}$. For a nonlinearity $f$ satisfying (\ref{eq_def_nonlinearity_local}), we consider the semilinear wave (or Klein-Gordon) equation
\[\left \{
\begin{array}{rcccl}\tag{$\ast$}\label{eq_NLKG}
\square u + \beta u & = & f(u) & \quad & \text{in } \mathbb{R} \times \Omega, \\
(u(0), \partial_t u(0)) & = & \left( u^0, u^1 \right) & \quad & \text{in } \Omega, \\
u & = & 0 & \quad & \text{on } \mathbb{R} \times \partial \Omega,
\end{array}
\right.\]
with real-valued initial data $\left( u^0, u^1 \right) \in H_0^1(\Omega) \times L^2(\Omega)$. If $\partial \Omega = \emptyset$, then the Dirichlet boundary condition can be removed. The local Cauchy theory for (\ref{eq_NLKG}) is well-known, and is given in Theorem \ref{thm_existence_nL_waves}. We say that $u \in H_0^1(\Omega)$ is a stationary solution of (\ref{eq_NLKG}) if $u$ is the solution of (\ref{eq_NLKG}) with initial data $(u, 0)$ and if $u$ is time-independent. If $s f(s) \leq 0$ for $s \in \mathbb{R}$, then (\ref{eq_NLKG}) is said to be defocusing. In this case, solutions of (\ref{eq_NLKG}) are globally defined, and the only stationary solution is $0$. If $s f(s) \geq 0$ for $s \in \mathbb{R}$, then (\ref{eq_NLKG}) is said to be focusing, and blow-up solutions and non-zero stationary solutions may exist (see, for example, \cite{Payne-Sattinger}). In this article, we make no assumption on the sign of $f$.

Consider $T > 0$, $a \in \mathscr{C}^\infty(\Omega, \mathbb{R})$ and $(\mathbf{u}^0, \mathbf{u}^1) \in H_0^1(\Omega) \times L^2(\Omega)$ such that the solution $\mathbf{u}$ of (\ref{eq_NLKG}) with initial data $(\mathbf{u}^0, \mathbf{u}^1)$ exists on the interval $[0, T]$.

\begin{defn}[Local controllability around $\mathbf{u}$ at time $T$]
We say that local controllability around $\mathbf{u}$ at time $T$ holds if there exists $\delta > 0$ such that for all $\left( u^0, u^1 \right) \in H_0^1(\Omega) \times L^2(\Omega)$ satisfying
\[\left\Vert \left( u^0, u^1 \right) - (\mathbf{u}^0, \mathbf{u}^1) \right\Vert_{H_0^1(\Omega) \times L^2(\Omega)} \leq \delta,\]
there exists $g \in L^1((0, T), L^2(\Omega))$ such that the solution $u$ of
\[\left \{
\begin{array}{rcccl}
\square u + \beta u & = & f(u) + ag & \quad & \text{in } (0, T) \times \Omega, \\
(u(0), \partial_t u(0)) & = & \left( u^0, u^1 \right) & \quad & \text{in } \Omega, \\
u & = & 0 & \quad & \text{on } (0, T) \times \partial \Omega,
\end{array}
\right.\]
satisfies $(u(T), \partial_t u(T)) = (\mathbf{u}(T), \partial_t \mathbf{u}(T))$.
\end{defn}

We will use the notion of generalized bicharacteristic, for which we refer to \cite{HormanderIII} and \cite{MelroseSjostrand}. In this article, we always assume that no generalized bicharacteristic has a contact of infinite order with $\mathbb{R} \times \partial \Omega$ (see \cite{BLR} for some details about this assumption).

\begin{defn}
For $\omega \subset \Omega$ open, we say that $(\omega, T)$ satisfies the \emph{Geometric Control Condition} (in short, GCC) if for every generalized bicharacteristic $s \mapsto \left(t(s), x(s), \tau(s), \xi(s) \right)$, there exists $s \in \mathbb{R}$ such that $t(s) \in (0, T)$ and $x(s) \in \omega$.
\end{defn}

Recall that a smooth function $F: (0, T) \times \Omega \rightarrow \mathbb{R}$ is said to be \emph{analytic with respect to $t$} if for all $\left( t_0, x_0 \right) \in (0, T) \times \Omega$, there exists a neighborhood $\mathcal{O} \subset (0, T) \times \Omega$ of $\left( t_0, x_0 \right)$ such that
\[F(t, x) = \sum_{k = 0}^\infty \partial_t^k F(t_0, x) \frac{\left( t - t_0 \right)^k}{k!}, \quad (t, x) \in \mathcal{O}.\]
The first result of this article is the following.

\begin{thm}[Local controllability around a trajectory]\label{thm_main_local_control_trajectory}
Assume that $f$ satisfies \textnormal{(\ref{eq_def_nonlinearity_local})}, and consider $(\mathbf{u}^0, \mathbf{u}^1) \in H_0^1(\Omega) \times L^2(\Omega)$ such that the solution $\mathbf{u}$ of \textnormal{(\ref{eq_NLKG})} with initial data $(\mathbf{u}^0, \mathbf{u}^1)$ exists on the interval $[0, T]$. We make the following assumptions.
\begin{enumerate}[label=(\roman*)]
\item Assume that there exist $\omega \subset \Omega$ open and $c > 0$ such that $a \geq c$ on $\omega$ and such that $(\omega, T)$ satisfies the GCC. In addition, if $\Omega$ is unbounded, assume that there exists $R_0 > 0$ such that $\mathbb{R}^d \backslash B(0, R_0) \subset \omega$. 
\item Assume that $f^\prime(\mathbf{u}) \in L^\infty((0, T) \times \Omega)$, and that $f^\prime(\mathbf{u})$ is smooth, and analytic with respect to $t$. In addition, if $\Omega$ is unbounded, assume that for all $t \in [0, T]$,
\[\left\vert \nabla_{\mathsf{g}} f^\prime(\mathbf{u}(t, x)) \right\vert + \left\vert f^\prime(\mathbf{u}(t, x)) \right\vert \xrightarrow{|x| \rightarrow \infty} 0,\]
where $\nabla_{\mathsf{g}}$ is the gradient with respect to the space variable $x$.
\end{enumerate}
Then, local controllability around $\mathbf{u}$ at time $T$ holds.
\end{thm}

In particular, if there exists a sequence $\left( u_n^0, u_n^1 \right) \in H_0^1(\Omega) \times L^2(\Omega)$ such that
\[\left\Vert \left( u_n^0, u_n^1 \right) - (\mathbf{u}^0, \mathbf{u}^1) \right\Vert_{H_0^1(\Omega) \times L^2(\Omega)} \xrightarrow{n \rightarrow \infty} 0,\]
and such that for all $n \in \mathbb{N}$, the solution of (\ref{eq_NLKG}) with initial data $\left( u_n^0, u_n^1 \right) \in H_0^1(\Omega) \times L^2(\Omega)$ blows up in finite time, then Theorem \ref{thm_main_local_control_trajectory} contains a controllability result for some blow-up solutions. An example of a solution $\mathbf{u}$ satisfying this condition is given below.

The second result of this article concerns the null-controllability in a long time of scattering solutions, in the case of an unbounded domain satisfying the non-trapping condition.

\begin{defn}[Null-controllability in a long time]
We say that null-controllability in a long time for $\left( u^0, u^1 \right) \in H_0^1(\Omega) \times L^2(\Omega)$ holds if there exist $T > 0$ and $g \in L^1((0, T), L^2(\Omega))$ such that the solution $u$ of
\[\left \{
\begin{array}{rcccl}
\square u + \beta u & = & f(u) + ag & \quad & \text{in } (0, T) \times \Omega, \\
(u(0), \partial_t u(0)) & = & \left( u^0, u^1 \right) & \quad & \text{in } \Omega, \\
u & = & 0 & \quad & \text{on } (0, T) \times \partial \Omega,
\end{array}
\right.\]
satisfies $(u(T), \partial_t u(T)) = 0$.
\end{defn}

A domain is said to be non-trapping if all generalized geodesics leave any compact set in finite time (see for example \cite{Melrose79} and \cite{MRS77} for a precise definition). When this condition is satisfied, resolvent estimates can be proven (see \cite{Burq03}, Remark 2.6, and references therein). For simplicity, we adopt these resolvent estimates as our definition of the non-trapping condition.

\begin{defn}\label{def_resolvent_estimate_non_trapping}
Assume that $\Omega$ is unbounded. We say that $\Omega$ is non-trapping if for all $\chi \in \mathscr{C}^\infty_\mathrm{c}(\Omega)$, there exists $C > 0$ such that
\[\sqrt{1 + \vert \lambda \vert} \left\Vert \chi \left( - \Delta + \lambda \right)^{-1} \chi u \right\Vert_{L^2(\Omega)} \leq C \Vert u \Vert_{L^2(\Omega)}, \quad u \in L^2(\Omega), \quad \Im \lambda \neq 0.\]
\end{defn}

We recall the definition of a scattering solution.

\begin{defn}
Consider $\left( u^0, u^1 \right) \in H_0^1(\Omega) \times L^2(\Omega)$. We say the solution $u_{\NL}$ of
\begin{equation}\label{eq_def_scatt}
\left \{
\begin{array}{rcccl}
\square u_{\NL} + \beta u_{\NL} & = & f(u_{\NL}) & \quad & \text{in } \mathbb{R}_+ \times \Omega, \\
(u_{\NL}(0), \partial_t u_{\NL}(0)) & = & \left( u^0, u^1 \right) & \quad & \text{in } \Omega, \\
u_{\NL} & = & 0 & \quad & \text{on } \mathbb{R}_+ \times \partial \Omega,
\end{array}
\right.
\end{equation}
is scattering (or that $\left( u^0, u^1 \right)$ is a scattering initial data) if $u_{\NL}$ exists on the whole interval $\mathbb{R}_+$ and satisfies 
\[\left\Vert \left( u_{\NL}(t), \partial_t u_{\NL}(t) \right) - \left( u_{\Lin}(t), \partial_t u_{\Lin}(t) \right) \right\Vert_{H_0^1(\Omega) \times L^2(\Omega)} \xrightarrow{t \rightarrow + \infty} 0\]
for some solution $u_{\Lin}$ of the linear equation
\[\left \{
\begin{array}{rcccl}
\square u_{\Lin} + \beta u_{\Lin} & = & 0 & \quad & \text{in } \mathbb{R}_+ \times \Omega, \\
u_{\Lin} & = & 0 & \quad & \text{on } \mathbb{R}_+ \times \partial \Omega.
\end{array}
\right.\]
\end{defn}

The second result of this article is the following.

\begin{thm}[Null-controllability of scattering solutions]\label{thm_main_null_control}
Assume that $\Omega$ is a non-trapping unbounded domain, and that $f$ satisfies \textnormal{(\ref{eq_def_nonlinearity_scatt})}, implying in particular that $3 \leq d \leq 5$. Consider $\omega \subset \Omega$ open, $R_0 > 0$ and $T > 0$ such that $(\omega, T)$ satisfies the GCC, and $\mathbb{R}^d \backslash B(0, R_0) \subset \omega$. Assume that there exists $c > 0$ such that $a \geq c$ on $\omega$. For $\left( u^0, u^1 \right) \in H_0^1(\Omega) \times L^2(\Omega)$, if the solution $u_{\NL}$ of \textnormal{(\ref{eq_def_scatt})} is scattering, then null-controllability in a long time for $\left( u^0, u^1 \right)$ holds.
\end{thm}

The proof is based on a local energy decay result (Theorem \ref{thm_local_decay_energy}), and global-in-time Strichartz estimates (Theorem \ref{thm_global_strichartz}), both of which can be of their own interest. The assumption that the metric $\mathsf{g}$ coincides with the Euclidean one outside a ball allows us to derive global-in-time Strichartz estimates on $\Omega$ from the local energy decay and the global-in-time Strichartz estimates on $\mathbb{R}^d$.

\begin{rem}
    Note that assuming that $(\omega, T)$ satisfies the GCC is not mandatory, as $\Omega$ is assumed to be non-trapping. However, since we have not explained the connection between the resolvent estimates of Definition \ref{def_resolvent_estimate_non_trapping} and the generalized bicharacteristics, we have chosen to make both assumptions explicit in the statement of Theorem \ref{thm_main_null_control}.
\end{rem}

\begin{rem}
We explain the restriction on the dimension $d$ in Theorem \ref{thm_main_null_control}. The proof relies on global-in-time Strichartz estimates of the form $\Vert f(u) \Vert_{L^1((0, +\infty), L^2(\Omega))} < +\infty$, where $u$ solves \textnormal{(\ref{eq_def_scatt})}. If $f(s) = s \vert s \vert^{\alpha - 1}$ for some $1 < \alpha < \frac{d + 2}{d - 2}$, then $\Vert f(u) \Vert_{L^1((0, +\infty), L^2(\Omega))} = \Vert u \Vert_{L^\alpha((0, +\infty), L^{2\alpha}(\Omega))}$. If $\alpha < 2$, which occurs whenever $d \geq 6$, we cannot rely on such an estimate (see Theorem 1.1 of \cite{HormanderTranslationInvariantOp} for a related result). Hence, a proof of Theorem \ref{thm_main_null_control} for $d \geq 6$ would typically require global-in-time Strichartz estimates involving Besov spaces.
\end{rem}

\begin{rem}
    If $\Omega$ is trapping, local energy decay and global-in-time Strichartz estimates generally fail to hold. Only weaker results are available; see, for example, \cite{burq_decroissance_1998}, \cite{Ikawa88}, \cite{Ikawa82}. Consequently, the method of proof of Theorem \ref{thm_main_null_control} does not apply, and addressing the stabilization of solutions would require different arguments, except perhaps in certain specific geometries where global Strichartz estimates do hold, as in \cite{lafontaine22}.
\end{rem}

\paragraph{Consequences.}
We give some examples of applications of Theorems \ref{thm_main_local_control_trajectory} and \ref{thm_main_null_control}. For shortness, we refer to the solution $u$ of
\[\left \{
\begin{array}{rcccl}\tag{$\ast \ast$}\label{eq_NLKG_controlled}
\square u + \beta u & = & f(u) + ag & \quad & \text{in } (0, + \infty) \times \Omega, \\
(u(0), \partial_t u(0)) & = & \left( u^0, u^1 \right) & \quad & \text{in } \Omega, \\
u & = & 0 & \quad & \text{on } (0, + \infty) \times \partial \Omega,
\end{array}
\right.\]
as the solution of (\ref{eq_NLKG_controlled}) with initial data $\left( u^0, u^1 \right)$ and with control $g$.

\paragraph{Consequences for some focusing nonlinearities.}
If $f$ is a focusing subcritical nonlinearity, then one can prove that there exists a special stationary solution $Q$ of (\ref{eq_NLKG}), called the ground state (see for example \cite{IbrahimMasmoudiNakanishi}, \cite{Payne-Sattinger}). The ground state is smooth, and decays at infinity if $\Omega$ is unbounded.

Firstly, we consider the case $\Omega = \mathbb{R}^d$, $\beta = 1$, with a nonlinearity $f$ satisfying the $H^1$-subcritical case of \cite{IbrahimMasmoudiNakanishi}, and satisfying (\ref{eq_def_nonlinearity_scatt}). An explicit example is $f(s) = s^3$, with $d = 3$. It is shown in \cite{IbrahimMasmoudiNakanishi} that the set of initial data with energy strictly below the energy of the ground state can be partitioned into two disjoint non-empty sets, $\mathcal{K}^+$ and $\mathcal{K}^-$, such that a solution initiated in $\mathcal{K}^+$ is globally defined and is scattering, while a solution initiated in $\mathcal{K}^-$ blows up in finite time. One can check that $\left( (1 \pm \epsilon) Q, 0 \right) \in \mathcal{K}^{\mp}$ if $\epsilon > 0$ is sufficiently small, implying the existence of both scattering solutions and blow-up solutions initiated arbitrarily close to the ground state. Another way to see the existence of scattering solution close to the ground state is to use the existence of a heteroclinic solution in the spirit of \cite{DuyckaertsMerle}, that is, a solution $W$ which is scattering (for positive time), and satisfies
\[\left\Vert \left( W(t), \partial_t W(t) \right) - (Q, 0) \right\Vert_{H^1(\mathbb{R}^d) \times L^2(\mathbb{R}^d)} \xrightarrow{t \rightarrow -\infty} 0.\]
The existence of such a solution $W$ is proved in \cite{NakanishiSchlagArticleRAD} and \cite{NakanishiSchlagArticleNONRAD}, in the case $f(s) = s^3$, $d = 3$. 

We prove that these results imply a controllability result from a neighborhood of $(Q, 0)$ to zero in a long time. If $(u^0, u^1)$ is close to $(Q, 0)$, then by the local controllability around $(Q, 0)$ (Theorem \ref{thm_main_local_control_trajectory}), applied two times, there exists $T > 0$ and a control $g$ such that $\left( u(T), \partial_t u(T) \right)$ is a scattering initial data, where $u$ is the solution of (\ref{eq_NLKG_controlled}) with initial data $\left( u^0, u^1 \right)$ and with control $g$. Then, the null-controllability of scattering solutions (Theorem \ref{thm_main_null_control}) implies the existence of $T^\prime > 0$, $g^\prime$, such that the solution $v$ of (\ref{eq_NLKG_controlled}) with initial data $\left( u(T), \partial_t u(T) \right)$ and with control $g^\prime$ satisfies $\left( v(T^\prime), \partial_t v(T^\prime) \right) = 0$. In particular, $g \mathds{1}_{[0, T]} + g^\prime \mathds{1}_{[T, T+T^\prime]}$ is a control such that the associated solution $u$ of (\ref{eq_NLKG_controlled}) with initial data $\left( u^0, u^1 \right)$ satisfies $\left( u(T+T^\prime), \partial_t u(T+T^\prime) \right) = 0$. Note that this implies the controllability of some blow-up solutions.

Secondly, we consider the case of a bounded domain $\Omega$, with $d = 3$ and $f(s) = s^3$. In this case, the sets $\mathcal{K}^+$ and $\mathcal{K}^-$ are defined by
\[\left \{
\begin{array}{c}
\mathscr{K}^+ = \left\{ \left( u^0, u^1 \right) \in H_0^1(\Omega) \times L^2(\Omega), E\left(u^0, u^1 \right) < E\left(Q, 0 \right), \left\Vert u^0 \right\Vert_{H_0^1(\Omega)}^2 \geq \left\Vert u^0 \right\Vert_{L^4(\Omega)}^4 \right\}, \\
\mathscr{K}^- = \left\{ \left( u^0, u^1 \right) \in H_0^1(\Omega) \times L^2(\Omega), E\left(u^0, u^1 \right) < E\left(Q, 0 \right), \left\Vert u^0 \right\Vert_{H_0^1(\Omega)}^2 < \left\Vert u^0 \right\Vert_{L^4(\Omega)}^4 \right\},
\end{array}
\right.\]
where the energy is given by $E\left( u^0, u^1 \right) = \frac{1}{2} \left\Vert u^0 \right\Vert_{H_0^1}^2 - \frac{1}{4} \left\Vert u^0 \right\Vert_{L^4}^4 + \frac{1}{2} \left\Vert u^1 \right\Vert_{L^2}^2$. A solution initiated in $\mathcal{K}^+$ is globally defined, and a solution initiated in $\mathcal{K}^-$ blows up in finite time (see \cite{Payne-Sattinger}). In \cite{perrin}, a stabilisation property under the GCC is shown for solutions initiated in $\mathcal{K}^+$. As above, using the fact that $\left( (1 \pm \epsilon) Q, 0 \right) \in \mathcal{K}^{\mp}$ if $\epsilon$ is sufficiently small, one obtains the controllability from a neighborhood of $(Q, 0)$ to zero in a long time. Note also that the stabilization of solutions initiated in $\mathscr{K}^+$, together with local controllability around zero, implies that exact controllability in $\mathcal{K}^+$ in a long time holds.

\paragraph{Consequences for some defocusing nonlinearities.}
Firstly, consider the case of an unbounded domain with a defocusing nonlinearity, such that all (finite-energy) solutions are scattering (see \cite{Brenner}, \cite{GinibreVelo}, \cite{NakanishiRemarks} for example). Assume, in addition, that $\Omega$ is nontrapping, and that $f$ satisfies (\ref{eq_def_nonlinearity_scatt}). Then, using the local controllability around zero (Theorem \ref{thm_main_local_control_trajectory}) and the null-controllability of scattering solutions (Theorem \ref{thm_main_null_control}), one concludes that exact controllability in a long time holds: for any initial data $(u^0, u^1)$ and final data $\left( v^0, v^1 \right)$, there exist a time $T > 0$ and a control $g$ such that $\left( u(T), \partial_t u(T) \right) = \left( v^0, v^1 \right)$, where $u$ is the solution of (\ref{eq_NLKG_controlled}) with initial data $\left( u^0, u^1 \right)$ and with control $g$. 

Secondly, consider the case of a domain $\Omega$ and a defocusing nonlinearity $f$ such that a stabilisation property holds, as in \cite{Alaoui}, \cite{DehmanLebeauZua}, \cite{Joly-Laurent} for example. Then, as above, it implies that exact controllability in a long time holds.

\paragraph{Connection with existing literature.}
The controllability of the linear wave equation has been extensively studied; see for example the seminal work of Bardos, Lebeau, and Rauch \cite{BLR}. For nonlinear equations with boundary control, some controllability results are proved in \cite{Lasiecka1991}, \cite{Lasiecka05}, \cite{TonBui}, and \cite{Zhou}. In the one-dimensional case, \cite{Zuazua1993} proves that exact controllability holds true for a quasilinear nonlinearity, and \cite{Fattorini} shows a local controllability result near equilibrium points. Still in dimension one, the authors of \cite{CORON2006} prove an exact controllability property from one stationary solution to another.

An exact controllability result for defocusing subcritical nonlinearities, between two high-frequency states, is proved in \cite{DehmanLebeau}, in the case of a bounded domain of dimension 3, under GCC. In \cite{JolyLaurentANoteOn}, a semi-global controllability result is established for asymptotically defocusing nonlinearities, where non-zero stationary solutions exist, but no blow-up phenomenon can occur. Local control around zero is studied in \cite{Zuazua1990} (and in \cite{chewning} in the case of boundary controllability). For nearly linear nonlinearities, such as those with global Lipschitzian or super-linear growth like $s \ln(s)^\beta$, we refer to \cite{Fu2007} and \cite{LiZhang}. For a local controllability result for a nonlinear Schrödinger equation, see \cite{LaurentSchrodingerControl} and \cite{ROSIER}. 

\paragraph{Outline of the article.} 
In Section 1, we present Strichartz estimates, gather some inequalities which follow from (\ref{eq_def_nonlinearity_local}) and (\ref{eq_def_nonlinearity_scatt}), and construct the solutions of (\ref{eq_NLKG}) and of certain time-dependent linear equations. 
In Section 2, we prove the local controllability result (Theorem \ref{thm_main_local_control_trajectory}): first, we show that the solution of the controlled system exists if the control is sufficiently small; second, we establish the exact controllability of the linearized problem; and finally, we complete the proof using a fixed-point argument. 
In Section 3, we prove the local decay of the energy and the global-in-time Strichartz estimates, and we show that they imply the null-controllability of scattering solutions (Theorem \ref{thm_main_null_control}). 
In the appendix, we provide a statement of the Christ-Kiselev Lemma and of a general lemma, prove an extension of Rellich's Theorem, establish the inequalities of Section 1 (Lemmas \ref{lem_basic_estimate_f_local} and \ref{lem_basic_estimate_f_scatt}), and derive an observability inequality for a linear wave equation on an exterior domain.


\paragraph{Notation.} We use the notation $A \lesssim B$ to denote that there exists a constant $C > 0$, depending only on certain fixed parameters specific to the context, such that $A \leq C B$.

\paragraph{Acknowledgments.} I warmly thank Thomas Duyckaerts and Jérôme Le Rousseau for their constant support and guidance, and Nicolas Burq for the idea used to prove the local decay of energy.

\section{Preliminaries}

\subsection{Strichartz estimates}

\begin{defn}[Local admissible exponents]
Consider $p, q \in \mathbb{R}$.
\begin{itemize}
\item Assume that $d \geq 3$ and $\partial \Omega = \emptyset$. Then, we say that $(p, q)$ is a pair of local admissible exponents for $\Omega$, and we write $(p, q) \in \Lambda_\Omega$, if $1 \leq p \leq \infty$, $2 \leq q \leq \frac{2d}{d - 3}$, $(q, d) \neq (+\infty, 3)$, and $\frac{1}{p} + \frac{d}{q} \geq \frac{d}{2} - 1$.
\item Assume that $d \geq 3$ and $\partial \Omega \neq \emptyset$. Then, we say that $(p, q)$ is a pair of local admissible exponents for $\Omega$, and we write $(p, q) \in \Lambda_\Omega$, if $1 \leq p \leq \infty$,
\[\frac{1}{p} + \frac{d}{q} \geq \frac{d}{2} - 1 \quad \text{and} \quad 
\left\{ 
\begin{array}{cc}
2 \leq q \leq 14 & \text{ if } d = 3 \\
2 \leq q \leq \frac{2 (d - 1)}{d - 3} & \text{ if } d \geq 4
\end{array}
\right. .\]
\item If $d = 2$, then we say that $(p, q)$ is a pair of local admissible exponents for $\Omega$, and we write $(p, q) \in \Lambda_\Omega$, if $1 \leq p \leq \infty$ and $1 \leq q < \infty$.
\end{itemize}
\end{defn}

A picture of the local admissible exponents can be found in Figure \ref{fig_strichartz_exponents}.

\begin{figure}[ht]
\centering
\begin{tabular}{cc}
\resizebox{0.41\textwidth}{!}
{
\begin{tikzpicture}[xscale=4,yscale=6]
    \draw[fill=gray!30] (1/2,0) -- (1,0) -- (1,1/2) -- (0,1/2) -- (0,1/6) -- cycle;
    \draw[->,>=stealth,gray!80] (0.75,0) .. controls +(-0.2,0.2) and +(-0.2,-0.0) .. (0.55,0.55);
    \node[scale=0.7] at (0.97,0.55) {This segment is not included.};
    \draw[->,>=stealth,gray!80] (0.3,1/6 - 1/3*0.3) .. controls +(0.05,0.1) and +(-0.1,0) .. (1.075,0.32);
    \node[scale=0.7] at (1.28,0.32) {$\frac{1}{p} + \frac{d}{q} = \frac{d}{2} - 1$};
    \draw[->,>=stealth] (-0.1,0) -- (1.2,0) node[below] {$\frac{1}{p}$};
    \draw[->,>=stealth] (0,-0.1) -- (0,0.6) node[left] {$\frac{1}{q}$};
    \draw (1/2,0) -- (0,1/6);
    \node[below] at (1/2,0) {$\frac{1}{2} = \frac{d - 2}{2}$};
    \draw (1/2,0-0.01) -- (1/2,0+0.01);
    \node[left] at (0,1/6) {$\frac{1}{6} = \frac{d - 2}{2d}$};
    \node at (0,1/6) {$\bullet$};
    \node[left] at (0,1/2) {$\frac{1}{2}$};
    \node at (0,1/2) {$\bullet$};
    \node[below] at (1,0) {$1$};
    \draw (1,0-0.01) -- (1,0+0.01);
    \draw[dashed,gray!30] (1/2,0) -- (1,0);
    \node at (0.75,-0.2) {$\underline{d = 3, \ \partial \Omega = \emptyset}$};
\end{tikzpicture}
}
&
\resizebox{0.47\textwidth}{!}
{
\begin{tikzpicture}[xscale=4,yscale=6]
    \draw[fill=gray!30] (1/2,1/5) -- (1,1/5) -- (1,1/2) -- (0,1/2) -- (0,3/10) -- cycle;
    \node[scale=0.7] at (1.4,-0.13) {$1 = \frac{d-2}{2} \text{ if } d = 4$};
    \draw[->,>=stealth,gray!80] (1,0) .. controls +(0.05,-0.02) and +(-0.1,0) .. (1.15,-0.13);
    \draw[->,>=stealth,gray!80] (3/2,0) .. controls +(-0.1,-0.05) and +(-0.15,0.1) .. (1.16,-0.126);
    \draw[->,>=stealth,gray!80] (0.8,3/10 - 1/5*0.8) .. controls +(0.05,0.1) and +(-0.1,0) .. (1.155,0.32);
    \node[scale=0.7] at (1.36,0.32) {$\frac{1}{p} + \frac{d}{q} = \frac{d}{2} - 1$};
    \draw[->,>=stealth] (-0.1,0) -- (1.7,0) node[below] {$\frac{1}{p}$};
    \draw[->,>=stealth] (0,-0.1) -- (0,0.6) node[left] {$\frac{1}{q}$};
    \draw (3/2,0) -- (0,3/10);
    \node[below] at (3/2,0) {$\frac{d - 2}{2}$};
    \node at (3/2,0) {$\bullet$};
    \node[left] at (0,3/10) {$\frac{d - 2}{2d}$};
    \node at (0,3/10) {$\bullet$};
    \node at (1/2,1/5) {$\bullet$};
    \draw[dotted] (1/2,1/5) -- (1/2,0);
    \draw[dotted] (1/2,1/5) -- (0,1/5);
    \node[below] at (1/2,0) {$\frac{1}{2}$};
    \node at (1/2,0) {$\bullet$};
    \node[left] at (0,1/5) {$\frac{d-3}{2d}$};
    \node at (0,1/5) {$\bullet$};
    \node[left] at (0,1/2) {$\frac{1}{2}$};
    \node at (0,1/2) {$\bullet$};
    \node[below] at (1,0) {$1$};
    \node at (1,0) {$\bullet$};
    \draw[dotted] (1,0) -- (1,1/5);
    \node at (0.75,-0.2) {$\underline{d \geq 4, \ \partial \Omega = \emptyset}$};
\end{tikzpicture}
}
\\
\resizebox{0.47\textwidth}{!}
{
\begin{tikzpicture}[xscale=8,yscale=6]
    \draw[fill=gray!30] (2/7,1/14) -- (0.7,1/14) -- (0.7,1/2) -- (0,1/2) -- (0,1/6) -- cycle;
    \draw[->,>=stealth,gray!80] (0.08,1/2 - 3/2*0.08) .. controls +(0.05,0.1) and +(-0.1,0) .. (0.18,0.55);
    \node[scale=0.7] at (0.25,0.55) {$\frac{3}{p} + \frac{2}{q} = 1$};
    \draw[->,>=stealth,gray!80] (0.1,1/6 - 1/3*0.1) .. controls +(0.05,0.1) and +(-0.1,0) .. (0.445,0.55);
    \node[scale=0.7] at (0.55,0.55) {$\frac{1}{p} + \frac{d}{q} = \frac{d}{2} - 1$};
    \draw[->,>=stealth] (-0.1,0) -- (0.8,0) node[below] {$\frac{1}{p}$};
    \draw[->,>=stealth] (0,-0.1) -- (0,0.6) node[left] {$\frac{1}{q}$};
    \draw (1/3,0) -- (0,1/2);
    \node[below] at (1/3,0) {$\frac{1}{3}$};
    \node at (1/3,0) {$\bullet$};
    \node[left] at (0,1/2) {$\frac{1}{2}$};
    \node at (0,1/2) {$\bullet$};
    \draw (1/2,0) -- (0,1/6);
    \node[below] at (1/2,0) {$\frac{1}{2} = \frac{d - 2}{2}$};
    \node at (1/2,0) {$\bullet$};
    \node[left] at (0,1/6) {$\frac{1}{6} = \frac{d - 2}{2d}$};
    \node at (0,1/6) {$\bullet$};
    \node at (2/7,1/14) {$\bullet$};
    \draw[dotted] (2/7,1/14) -- (2/7,0);
    \draw[dotted] (2/7,1/14) -- (0,1/14);
    \node[below] at (2/7,0) {$\frac{2}{7}$};
    \node at (2/7,0) {$\bullet$};
    \node[left] at (0,1/14) {$\frac{1}{14}$};
    \node at (0,1/14) {$\bullet$};
    \node[below] at (0.7,0) {$1$};
    \draw[dotted] (0.7,0) -- (0.7,1/14);
    \node at (0.4,-0.2) {$\underline{d = 3, \ \partial \Omega \neq \emptyset}$};
\end{tikzpicture}
}
&
\resizebox{0.47\textwidth}{!}
{
\begin{tikzpicture}[xscale=4,yscale=6]
    \draw[fill=gray!30] (1/4,1/4) -- (1,1/4) -- (1,1/2) -- (0,1/2) -- (0,3/10) -- cycle;
    \node[scale=0.7] at (1.4,-0.13) {$1 = \frac{d-2}{2} \text{ if } d = 4$};
    \draw[->,>=stealth,gray!80] (1,0) .. controls +(0.05,-0.02) and +(-0.1,0) .. (1.15,-0.13);
    \draw[->,>=stealth,gray!80] (3/2,0) .. controls +(-0.1,-0.05) and +(-0.15,0.1) .. (1.16,-0.126);
    \draw[->,>=stealth,gray!80] (0.08,1/2 - 0.08) .. controls +(0.05,0.1) and +(-0.1,0) .. (0.30,0.55);
    \node[scale=0.7] at (0.45,0.55) {$\frac{1}{p} + \frac{1}{q} = \frac{1}{2}$};
    \draw[->,>=stealth,gray!80] (0.8,3/10 - 1/5*0.8) .. controls +(0.05,0.1) and +(-0.1,0) .. (1.155,0.22);
    \node[scale=0.7] at (1.36,0.22) {$\frac{1}{p} + \frac{d}{q} = \frac{d}{2} - 1$};
    \draw[->,>=stealth] (-0.1,0) -- (1.7,0) node[below] {$\frac{1}{p}$};
    \draw[->,>=stealth] (0,-0.1) -- (0,0.6) node[left] {$\frac{1}{q}$};
    \draw (1/2,0) -- (0,1/2);
    \node[below] at (1/2,0) {$\frac{1}{2}$};
    \node at (1/2,0) {$\bullet$};
    \node[left] at (0,1/2) {$\frac{1}{2}$};
    \node at (0,1/2) {$\bullet$};
    \draw (3/2,0) -- (0,3/10);
    \node[below] at (3/2,0) {$\frac{d-2}{2}$};
    \node at (3/2,0) {$\bullet$};
    \node[left] at (0,3/10+0.02) {$\frac{d-2}{2d}$};
    \node at (0,3/10) {$\bullet$};
    \node at (1/4,1/4) {$\bullet$};
    \draw[dotted] (1/4,1/4) -- (1/4,0);
    \draw[dotted] (1/4,1/4) -- (0,1/4);
    \node[below] at (1/4,0) {$\frac{1}{d-1}$};
    \node at (1/4,0) {$\bullet$};
    \node[left] at (0,1/4-0.02) {$\frac{d-3}{2(d-1)}$};
    \node at (0,1/4) {$\bullet$};
    \node[below] at (1,0) {$1$};
    \node at (1,0) {$\bullet$};
    \draw[dotted] (1,0) -- (1,1/4);
    \node at (0.75,-0.2) {$\underline{d \geq 4, \ \partial \Omega \neq \emptyset}$};
\end{tikzpicture}
}
\end{tabular}
\caption{The local admissible exponents $\Lambda_\Omega$, in gray.}
\label{fig_strichartz_exponents}
\end{figure}

\begin{thm}[Local-in-time Strichartz estimates]\label{thm_local_strichartz_estimates}
Consider $T > 0$. There exists $C > 0$ such that for all $(p, q) \in \Lambda_\Omega$, $\left( u^0, u^1 \right) \in H_0^1(\Omega) \times L^2(\Omega)$, and $g \in L^1([0, T], L^2(\Omega))$, the unique solution $u$ of 
\[\left \{
\begin{array}{rcccl}
\square u + \beta u & = & g & \quad & \text{in } \mathbb{R} \times \Omega, \\
(u(0), \partial_t u(0)) & = & \left( u^0, u^1 \right) & \quad & \text{in } \Omega, \\
u & = & 0 & \quad & \text{on } \mathbb{R} \times \partial \Omega,
\end{array}
\right.\]
satisfies
\begin{equation}\label{eq_thm_strichartz_est}
\Vert u \Vert_{L^p([0, T], L^q)} \leq C \left( \left\Vert \left( u^0, u^1 \right) \right\Vert_{H_0^1(\Omega) \times L^2(\Omega)} + \Vert g \Vert_{L^1([0, T], L^2)} \right).
\end{equation}
\end{thm}

\begin{proof}
If $d = 2$, then the Sobolev embedding $H^1(\Omega) \hookrightarrow L^q(\Omega)$ holds true for $1 \leq q < + \infty$ (see for example \cite{Adams-Fournier}, 4.12 Part I Case B, with $n = p = 2$ and $m = 1$). Hence, in that case, one has
\[\Vert u \Vert_{L^p([0, T], L^q)} \lesssim \Vert u \Vert_{L^\infty([0, T], H_0^1)} \lesssim \left\Vert \left( u^0, u^1 \right) \right\Vert_{H_0^1(\Omega) \times L^2(\Omega)} + \Vert g \Vert_{L^1([0, T], L^2)},\]
by classical semi-group theory, for all $(p, q) \in \Lambda_\Omega$.

Now, we assume that $d \geq 3$. Note that an estimate for the wave equation can be used for the Klein-Gordon equation, as one can absorb the low-order term for $T$ sufficiently small, and iterate to get the result for large $T$.

\paragraph{Case 1:} $\Omega = \mathbb{R}^d$. Here, we rely on \cite{KeelTao} (which is a generalisation of \cite{Kapitanski}, \cite{LindbladSogge} and \cite{MockenhauptSeegerSogge}). See \cite{TataruStrich} for the case of non-smooth coefficients. The original result in $\mathbb{R}^3$ was proved by Strichartz \cite{Strichartz_original}. By Corollary 1.3 of \cite{KeelTao} with $(\tilde{q}, \tilde{r}, \gamma) = (+ \infty, 2, 1)$, (\ref{eq_thm_strichartz_est}) holds true if $2 \leq p \leq \infty$, $2 \leq q < \infty$, $(p, q, d) \neq (2, \infty, 3)$, 
\[\frac{1}{p} + \frac{d}{q} = \frac{d}{2} - 1 \quad \text{and} \quad \frac{2}{p} + \frac{d - 1}{q} \leq \frac{d - 1}{2}.\]
Note that if $\frac{1}{p} + \frac{d}{q} = \frac{d}{2} - 1$, then $\frac{2}{p} + \frac{d - 1}{q} \leq \frac{d - 1}{2}$ is equivalent with $\frac{1}{p} \leq \frac{d - 1}{d + 1}$, a condition which is weaker than $p \geq 2$. Note also that $p \geq 2$ implies $\frac{1}{d} \left(  \frac{d}{2} - 1 - \frac{1}{p} \right) \geq \frac{d - 3}{2d}$. Hence, (\ref{eq_thm_strichartz_est}) holds true if $2 \leq p \leq \infty$, $2 \leq q < \frac{2d}{d - 3}$, $(p, q, d) \neq (2, \infty, 3)$, and $\frac{1}{p} + \frac{d}{q} = \frac{d}{2} - 1$. Classical semi-group theory gives
\[\Vert u \Vert_{L^\infty([0, T], H^1)} \leq C \left( \left\Vert \left( u^0, u^1 \right) \right\Vert_{H_0^1(\Omega) \times L^2(\Omega)} + \Vert g \Vert_{L^1([0, T], L^2)} \right).\]
Using also the Sobolev embedding $H^1(\Omega) \hookrightarrow L^{\frac{2d}{d - 2}}(\Omega)$, one finds (\ref{eq_thm_strichartz_est}) with $(p, q) = (+\infty, 2)$ and $(p, q) = \left( +\infty, \frac{2d}{d - 2} \right)$. By interpolation, (\ref{eq_thm_strichartz_est}) holds true with $(p, q) = (+\infty, q)$, for all $q \in \left[2, \frac{2d}{d - 2} \right]$. Finally, as the estimate is local in time, if (\ref{eq_thm_strichartz_est}) holds true for some $(p, q)$, then it also holds true for $(p_0, q)$, for all $p_0 \in [1, p]$. 

\paragraph{Case 2:} $\Omega$ is a compact manifold with boundary. In that case, we refer to the results of \cite{BlairSmithSogge} (which extend those of \cite{BurqLebeauPlanchon}). Note that Ivanovici's counterexamples in \cite{Ivanovici} show that Strichartz estimates are not true for the full range of exponents in the case of a manifold with boundary. By Corollary 1.2 of \cite{BlairSmithSogge}, applied with $(r, s, \gamma) = (1, 2, 1)$, (\ref{eq_thm_strichartz_est}) holds true if $2 < p \leq \infty$, $2 \leq q < \infty$, 
\[\frac{1}{p} + \frac{d}{q} = \frac{d}{2} - 1 \quad \text{and} \quad 
\left\{ 
\begin{array}{cc}
\frac{3}{p} + \frac{2}{q} \leq 1 & \text{ if } d = 3 \\
\frac{1}{p} + \frac{1}{q} \leq \frac{1}{2} & \text{ if } d \geq 4
\end{array}
\right. .\]
If $\frac{1}{p} + \frac{d}{q} = \frac{d}{2} - 1$ and $d = 3$, then $\frac{3}{p} + \frac{2}{q} \leq 1$ is equivalent with $q \leq 14$. Similarly, if $\frac{1}{p} + \frac{d}{q} = \frac{d}{2} - 1$ and $d \geq 4$, then $\frac{1}{p} + \frac{1}{q} \leq \frac{1}{2}$ is equivalent with $q \leq \frac{2(d - 1)}{d - 3}$. As above, (\ref{eq_thm_strichartz_est}) holds true with $(p, q) = (+ \infty, q)$, for all $q \in \left[2, \frac{2d}{d - 2} \right]$, and if (\ref{eq_thm_strichartz_est}) holds true for some $(p, q)$, then it also holds true for $(p_0, q)$, for all $p_0 \in [1, p]$. This proves that (\ref{eq_thm_strichartz_est}) holds true for all $(p, q) \in \Lambda_\Omega$.

\paragraph{Case 3:} Next, assume that $\Omega$ is a compact manifold without boundary. Let $\left( O_j \right)_{j \in J}$ be a finite family of open subsets of $\Omega$ such that each $O_j$ is included in a coordinate chart of $\Omega$, and such that 
\[\Omega = \bigcup_{j \in J} O_j.\]
Let $\left( \psi_j \right)_{j \in J}$ be such that $\psi_j \in \mathscr{C}^\infty_\mathrm{c}(O_j, [0, 1])$ for $j \in J$, and $\sum_{j \in J} \psi_j = 1$ on $\Omega$. For $j \in J$, let $u_j$ be the solution of 
\[\left \{
\begin{array}{rcccl}
\square u_j + u_j & = & \psi_j g & \quad & \text{in } \mathbb{R} \times \Omega, \\
(u_j(0), \partial_t u_j(0)) & = & (\psi_j u^0, \psi_j u^1) & \quad & \text{in } \Omega.
\end{array}
\right.\]
One has $u = \sum_{j \in J} u_j$. For $(p, q) \in \Lambda_\Omega$ and $T > 0$, write 
\[\Vert u \Vert_{L^p([0, T], L^q)} \leq \sum_{j \in J} \left\Vert u_j \right\Vert_{L^p([0, T], L^q)}.\]
If $T$ is sufficiently small, then by finite speed of propagation, $u_j$ is supported in $O^j$ for all $j \in J$. As $O^j$ is supported in a coordinate chart of $\Omega$, we can apply Strichartz estimate in the case of $\mathbb{R}^d$ (with variable coefficients), to find
\begin{align*}
\Vert u \Vert_{L^p([0, T], L^q)} & \lesssim \sum_{j \in J} \left( \left\Vert \left( \psi_j u^0, \psi_j u^1 \right) \right\Vert_{H_0^1(\Omega) \times L^2(\Omega)} + \Vert \psi_j g \Vert_{L^1([0, T], L^2)} \right) \\
& \lesssim \left\Vert \left( u^0, u^1 \right) \right\Vert_{H_0^1(\Omega) \times L^2(\Omega)} + \Vert g \Vert_{L^1([0, T], L^2)}.
\end{align*}

\paragraph{Case 4:} Finally, assume that $\Omega$ is a compact perturbation of $\mathbb{R}^d$, and write $\Omega = \mathbb{R}^d \backslash U$. Fix $R > 0$ such that $U \subset B(0, R)$ and such that the metric of $\Omega \cap B(0, R)^\complement$ is equal to the Euclidean one. Let $\left( O_j \right)_{j \in J}$ be a finite family of open subsets of $\Omega \cap B(0, R + 2)$ such that 
\[\Omega \cap B(0, R + 1) \subset \bigcup_{j \in J} O_j.\]
There exist $\psi_0$ and $\left( \psi_j \right)_{j \in J}$, satisfying the following properties: $\psi_j \in \mathscr{C}^\infty_\mathrm{c}(O_j, [0, 1])$ for $j \in J$, $\psi_0 \in \mathscr{C}^\infty(\mathbb{R}^d, [0, 1])$, with $\psi_0 = 0$ on $B(0, R)$ and $\psi_0 = 1$ on $\mathbb{R}^d \backslash B(0, R + 1)$, and $\psi_0 + \sum_{j \in J} \psi_j = 1$ on $\Omega$. Write $u_0$ and $(u_j)_{j \in J}$ as in the case of a manifold without boundary. Note that $\Lambda_{\Omega} \subset \Lambda_{\mathbb{R}^d}$. Hence, applying Strichartz estimates in the case of a compact manifold for the functions $u_j$, $j \in J$, and Strichartz estimates in the case of $\mathbb{R}^d$ for $u_0$, one completes the proof as in the case of a manifold without boundary.
\end{proof}

\begin{thm}[Global-in-time Strichartz estimates]\label{thm_global_strichartz}
Assume that $\Omega$ is either a non-trapping exterior domain in $\mathbb{R}^d$, with $d = 3, 4$, or $\Omega = \mathbb{R}^d$, with $3 \leq d \leq 5$. Consider $2 < \alpha < \frac{d + 2}{d - 2}$. There exists $C > 0$ such that for all $\left( u^0, u^1 \right) \in H_0^1(\Omega) \times L^2(\Omega)$ and $F \in L^1(\mathbb{R}, L^2(\Omega))$, the solution $u$ of 
\begin{equation}\label{eq_thm_global_strichartz_1}
\left \{
\begin{array}{rcccl}
\square u + \beta u & = & F & \quad & \text{in } \mathbb{R} \times \Omega, \\
(u(0), \partial_t u(0)) & = & \left( u^0, u^1 \right) & \quad & \text{in } \Omega, \\
u & = & 0 & \quad & \text{on } \mathbb{R} \times \partial \Omega,
\end{array}
\right.
\end{equation}
satisfies
\begin{equation}\label{eq_thm_global_strichartz_2}
\left\Vert u \right\Vert_{L^\alpha(\mathbb{R}, L^{2 \alpha})} \leq C \left( \left\Vert \left( u^0, u^1 \right) \right\Vert_{H_0^1(\Omega) \times L^2(\Omega)} + \Vert F \Vert_{L^1(\mathbb{R}, L^2)} \right).
\end{equation}
\end{thm}

A proof of Theorem \ref{thm_global_strichartz} is given in Section 3. Note that a global-in-time Strichartz estimates implies a local-in-time Strichartz estimate with a constant independent of the time. 

\subsection{Basic nonlinear estimates}

Here, we gather some estimates involving nonlinearities $f$ satisfying (\ref{eq_def_nonlinearity_local}) or (\ref{eq_def_nonlinearity_scatt}), which essentially result from Hölder's inequality. The proofs are given in Appendix \ref{appendix_proof_lem_Holder}. For $\mathbf{u}, h: (0, T) \times \Omega \rightarrow \mathbb{R}$, set 
\begin{equation}\label{eq_def_NL}
\NL_\mathbf{u}(h) = f(\mathbf{u} + h) - f(\mathbf{u}) - f^\prime(\mathbf{u}) h.
\end{equation}

\begin{lem}[Basic nonlinear estimates - 1]\label{lem_basic_estimate_f_local}
Consider $f$ satisfying \textnormal{(\ref{eq_def_nonlinearity_local})} for some $\alpha$. For $T > 0$, set
\begin{equation}\label{eq_def_X_T}
X_T = \mathscr{C}^0([0, T], H_0^1(\Omega)) \ \cap \ \mathscr{C}^1([0, T], L^2(\Omega)) \ \cap \ L^\alpha((0, T), L^{2 \alpha}(\Omega))
\end{equation}
with 
\[\left\Vert u \right\Vert_{X_T} = \max \left( \left\Vert u \right\Vert_{L^\infty([0, T], H_0^1)}, \left\Vert \partial_t u \right\Vert_{L^\infty([0, T], L^2)}, \left\Vert u \right\Vert_{L^\alpha((0, T), L^{2 \alpha})} \right).\]
\begin{enumerate}[label=(\roman*)]
\item There exists $C > 0$ such that for all $T > 0$, and all $u, v \in X_T$, one has
\[\left\Vert f(u) - f(v) \right\Vert_{L^1((0, T), L^2)} \leq C \left\Vert u - v \right\Vert_{X_T} \left(T + \left\Vert u \right\Vert_{L^\alpha((0, T), L^{2 \alpha})}^{\alpha - 1} + \left\Vert v \right\Vert_{L^\alpha((0, T), L^{2 \alpha})}^{\alpha - 1} \right).\]
\item Consider $T > 0$ and $\mathbf{u} \in L^\alpha((0, T), L^{2 \alpha}(\Omega))$. There exists $C = C(T, \mathbf{u}) > 0$ such that for $u, v \in X_T$, one has
\[\left\Vert \NL_\mathbf{u}(u) - \NL_\mathbf{u}(v) \right\Vert_{L^1((0, T), L^2)} \leq C \left\Vert u - v \right\Vert_{X_T} \left(\left\Vert u \right\Vert_{X_T} + \left\Vert v \right\Vert_{X_T} + \left\Vert u \right\Vert_{X_T}^{\alpha - 1} + \left\Vert v \right\Vert_{X_T}^{\alpha - 1} \right).\]
\item There exists $C > 0$ such that for $T > 0$ and $\mathbf{u} \in L^{\alpha}((0, T), L^{2 \alpha}(\Omega))$, one has
\[\left\Vert f^\prime(\mathbf{u}) u \right\Vert_{L^1((0, T), L^2)} \leq C \left\Vert u \right\Vert_{X_T} \left( T + \left\Vert \mathbf{u} \right\Vert_{L^{\alpha}((0, T), L^{2 \alpha})}^{\alpha - 1} \right),\]
for all $u \in X_T$.
\end{enumerate}
\end{lem}

\begin{rem}
If $d = 2$, then one has $X_T = \mathscr{C}^0([0, T], H^1(\Omega)) \cap \mathscr{C}^1([0, T], L^2(\Omega))$, and 
\[\left\Vert u \right\Vert_{X_T} \lesssim \max \left( \left\Vert u \right\Vert_{L^\infty([0, T], H^1)}, \left\Vert \partial_t u \right\Vert_{L^\infty([0, T], L^2)} \right), \quad u \in X_T,\]
by the Sobolev embedding $H^1(\Omega) \hookrightarrow L^p(\Omega)$ for $1 \leq p < + \infty$ (see the beginning of the proof of Theorem \ref{thm_local_strichartz_estimates}).
\end{rem}

\begin{lem}[Basic nonlinear estimates - 2]\label{lem_basic_estimate_f_scatt}
Consider $f$ satisfying \textnormal{(\ref{eq_def_nonlinearity_scatt})} for some $\alpha_0 \leq \alpha_1$. 
\begin{enumerate}[label=(\roman*)]
\item There exists $C > 0$ such that for all $T > 0$, one has
\begin{align*}
\left\Vert f(u) - f(v) \right\Vert_{L^1((0, T), L^2)} \leq \ & C \left\Vert u - v \right\Vert_{L^{\alpha_0}((0, T), L^{2 \alpha_0})} \left(\left\Vert u \right\Vert_{L^{\alpha_0}((0, T), L^{2 \alpha_0})}^{\alpha_0 - 1} + \left\Vert v \right\Vert_{L^{\alpha_0}((0, T), L^{2 \alpha_0})}^{\alpha_0 - 1} \right) \\
+ \, & C \left\Vert u - v \right\Vert_{L^{\alpha_1}((0, T), L^{2 \alpha_1})} \left(\left\Vert u \right\Vert_{L^{\alpha_1}((0, T), L^{2 \alpha_1})}^{\alpha_1 - 1} + \left\Vert v \right\Vert_{L^{\alpha_1}((0, T), L^{2 \alpha_1})}^{\alpha_1 - 1} \right),
\end{align*}
for all $u, v \in L^{\alpha_0}((0, T), L^{2 \alpha_0}(\Omega)) \cap L^{\alpha_1}((0, T), L^{2 \alpha_1}(\Omega))$.
\item Consider $T > 0$ and $\mathbf{u} \in L^{\alpha_0}((0, T), L^{2 \alpha_0}(\Omega)) \cap L^{\alpha_1}((0, T), L^{2 \alpha_1}(\Omega))$. Recall that $\NL_\mathbf{u}$ is defined by \textnormal{(\ref{eq_def_NL})}. There exists $C = C(T, \mathbf{u}) > 0$ such that 
\begin{align*}
& \left\Vert \NL_\mathbf{u}(u) - \NL_\mathbf{u}(v) \right\Vert_{L^1((0, T), L^2)} \leq C \sum_{i = 0, 1} \left\Vert u - v \right\Vert_{L^{\alpha_i}((0, T), L^{2 \alpha_i})} \left(\left\Vert u \right\Vert_{L^{\alpha_i}((0, T), L^{2 \alpha_i})} \right. \\
& \hspace{3cm} \left. + \left\Vert v \right\Vert_{L^{\alpha_i}((0, T), L^{2 \alpha_i})} + \left\Vert u \right\Vert_{L^{\alpha_i}((0, T), L^{2 \alpha_i})}^{\alpha_i - 1} + \left\Vert v \right\Vert_{L^{\alpha_i}((0, T), L^{2 \alpha_i})}^{\alpha_i - 1} \right),
\end{align*}
for all $u, v \in L^{\alpha_0}((0, T), L^{2 \alpha_0}(\Omega)) \cap L^{\alpha_1}((0, T), L^{2 \alpha_1}(\Omega))$. In addition, one can assume that $C(T, \mathbf{u}) \leq  C(1, \mathbf{u})$ if $T \in (0, 1]$. 
\item There exists $C > 0$ such that for all $T > 0$ and $\mathbf{u} \in L^{\alpha_0}((0, T), L^{2 \alpha_0}(\Omega)) \cap L^{\alpha_1}((0, T), L^{2 \alpha_1}(\Omega))$, one has
\[\left\Vert f^\prime(\mathbf{u}) u \right\Vert_{L^1((0, T), L^2)} \leq C \sum_{i = 0, 1} \left\Vert \mathbf{u} \right\Vert_{L^{\alpha_i}((0, T), L^{2 \alpha_i})}^{\alpha_i - 1} \left\Vert u \right\Vert_{L^{\alpha_i}((0, T), L^{2 \alpha_i})},\]
for all $u \in L^{\alpha_0}((0, T), L^{2 \alpha_0}(\Omega)) \cap L^{\alpha_1}((0, T), L^{2 \alpha_1}(\Omega))$.
\end{enumerate}
\end{lem}

\subsection{Solutions of the wave equations}

In this article, we only consider real-valued solutions of wave equations. We start with the case of linear wave equations with time-dependent potential. 

\begin{prop}[Solution of a linear wave equation with time-dependent potential]\label{prop_LKG_time_dependent_potential}
\
\begin{enumerate}[label=(\roman*)]
\item Consider $f$ satisfying \textnormal{(\ref{eq_def_nonlinearity_local})} for some $\alpha$, $T > 0$. Write $X_T$ for the set defined by \textnormal{(\ref{eq_def_X_T})}, and consider $\mathbf{u} \in X_T$. Then, for all $\left( u^0, u^1 \right) \in H_0^1(\Omega) \times L^2(\Omega)$ and $g \in L^1((0, T), L^2(\Omega))$, there exists a unique solution $u \in X_T$ of
\[\left \{
\begin{array}{rcccl}
\square u + \beta u & = & f^\prime(\mathbf{u}) u + g & \quad & \text{in } (0, T) \times \Omega, \\
(u(0), \partial_t u(0)) & = & \left( u^0, u^1 \right) & \quad & \text{in } \Omega, \\
u & = & 0 & \quad & \text{on } (0, T) \times \partial \Omega.
\end{array}
\right.\]
In addition, there exists $C > 0$, independent of $\left( u^0, u^1 \right)$ and $g$, such that
\[\Vert u \Vert_{X_T} \leq C \left( \left\Vert \left( u^0, u^1 \right) \right\Vert_{H_0^1(\Omega) \times L^2(\Omega)} + \Vert g \Vert_{L^1((0, T), L^2)} \right).\]

\item Consider $f$ satisfying \textnormal{(\ref{eq_def_nonlinearity_scatt})} for some $\alpha_0 \leq \alpha_1$, and $T_1, T_2 \in \mathbb{R}$, with $T_1 \leq T_2$. Set
\begin{align}
Y_{\left[ T_1, T_2 \right]} = \ & \mathscr{C}^0(\left[ T_1, T_2 \right], H_0^1(\Omega)) \, \cap \, \mathscr{C}^1(\left[ T_1, T_2 \right], L^2(\Omega)) \nonumber \\
& \cap \, L^{\alpha_0}(\left( T_1, T_2 \right), L^{2 \alpha_0}(\Omega)) \, \cap \, L^{\alpha_1}(\left( T_1, T_2 \right), L^{2 \alpha_1}(\Omega)), \label{eq_def_Y_T}
\end{align}
with
\begin{align*}
\left\Vert u \right\Vert_{Y_{\left[ T_1, T_2 \right]}} = \max & \left( \left\Vert u \right\Vert_{L^\infty(\left[ T_1, T_2 \right], H_0^1)}, \left\Vert \partial_t u \right\Vert_{L^\infty(\left[ T_1, T_2 \right], L^2)}, \right. \\
& \ \, \left. \left\Vert u \right\Vert_{L^{\alpha_0}(\left( T_1, T_2 \right), L^{2 \alpha_0})}, \left\Vert u \right\Vert_{L^{\alpha_1}(\left( T_1, T_2 \right), L^{2 \alpha_1})} \right).
\end{align*}
Fix $\mathbf{u} \in Y_{\left[ T_1, T_2 \right]}$. For all $\left( u^0, u^1 \right) \in H_0^1(\Omega) \times L^2(\Omega)$ and $g \in L^1(\left( T_1, T_2 \right), L^2(\Omega))$, there exists a unique solution $u \in Y_{\left[ T_1, T_2 \right]}$ of 
\[\left \{
\begin{array}{rcccl}
\square u + \beta u & = & f^\prime(\mathbf{u}) u + g & \quad & \text{in } \left( T_1, T_2 \right) \times \Omega, \\
\left( u(T_1), \partial_t u(T_1) \right) & = & \left( u^0, u^1 \right) & \quad & \text{in } \Omega, \\
u & = & 0 & \quad & \text{on } \left( T_1, T_2 \right) \times \partial \Omega.
\end{array}
\right.\]
In addition, there exists $C > 0$, independent of $\left( u^0, u^1 \right)$ and $g$, such that
\[\Vert u \Vert_{Y_{\left[ T_1, T_2 \right]}} \leq C \left( \left\Vert \left( u^0, u^1 \right) \right\Vert_{H_0^1(\Omega) \times L^2(\Omega)} + \Vert g \Vert_{L^1(\left( T_1, T_2 \right), L^2)} \right).\]
\item Consider $V \in L^\infty((0, T), \mathscr{C}^1(\Omega))$ satisfying
\begin{equation}\label{eq_prop_LKG_time_dependent_potential_1}
\left\Vert V \right\Vert_{L^\infty((0, T) \times \Omega)} + \sum_{j = 1}^d \left\Vert \partial_{x^j} V \right\Vert_{L^\infty((0, T) \times \Omega)} < + \infty.
\end{equation}
Then for all $\left( u^0, u^1 \right) \in L^2(\Omega) \times H^{-1}(\Omega)$ and $g \in L^1((0, T), H^{-1}(\Omega))$, there exists a unique solution $u \in \mathscr{C}^0([0, T], L^2(\Omega)) \cap \mathscr{C}^1([0, T], H^{-1}(\Omega))$ of 
\[\left \{
\begin{array}{rcccl}
\square u + \beta u & = & V u + g & \quad & \text{in } (0, T) \times \Omega, \\
(u(0), \partial_t u(0)) & = & \left( u^0, u^1 \right) & \quad & \text{in } \Omega, \\
u & = & 0 & \quad & \text{on } (0, T) \times \partial \Omega.
\end{array}
\right.\]
In addition, there exists $C > 0$, independent of $\left( u^0, u^1 \right)$ and $g$, such that
\[\left\Vert u \right\Vert_{L^\infty([0, T], L^2)} + \left\Vert \partial_t u \right\Vert_{L^\infty([0, T], H^{-1})} \leq C \left( \left\Vert \left( u^0, u^1 \right) \right\Vert_{L^2(\Omega) \times H^{-1}(\Omega)} + \Vert g \Vert_{L^1((0, T), H^{-1})} \right).\]
\end{enumerate} 
\end{prop}

\begin{proof}
We start with the proof of \emph{(i)} in the case $\mathbf{u} = 0$. In that case, classical semi-group theory gives $u \in \mathscr{C}^0([0, T], H_0^1(\Omega)) \cap \mathscr{C}^1([0, T], L^2(\Omega))$. Hence, it suffices to prove that
\[\Vert u \Vert_{L^\alpha((0, T), L^{2 \alpha})} \lesssim \left( \left\Vert \left( u^0, u^1 \right) \right\Vert_{H_0^1(\Omega) \times L^2(\Omega)} + \Vert g \Vert_{L^1((0, T), L^2)} \right).\]
This follows from Theorem \ref{thm_local_strichartz_estimates}, if we show that (\ref{eq_def_nonlinearity_local}) implies $(\alpha, 2 \alpha) \in \Lambda_\Omega$. This holds true if $d = 2$. Assume that $d \geq 3$. One has $\alpha < \frac{d + 2}{d - 2}$, implying
\begin{equation}\label{proof_prop_time_dep_potential_eq_for_later}
\frac{1}{\alpha} + \frac{d}{2 \alpha} = \frac{d + 2}{2 \alpha} > \frac{d - 2}{2}.
\end{equation}
If $\partial \Omega = \emptyset$ and $d \geq 3$, then $\alpha \leq \frac{d + 2}{d - 2}$ implies $2 \alpha \leq \frac{2 d}{d - 3}$, yielding $(\alpha, 2 \alpha) \in \Lambda_\Omega$. Now, assume that $\partial \Omega \neq \emptyset$. If $d = 3$, then $\alpha \leq \frac{d + 2}{d - 2} = 5$ implies $2 \alpha \leq 14$. If $d = 4$, then $\alpha \leq \frac{d + 2}{d - 2} = 3$ implies $2 \alpha \leq \frac{2 ( d - 1 )}{d - 3} = 6$. If $d \geq 5$, then one has $2 \alpha \leq \frac{2 ( d - 1 )}{d - 3}$ by assumption. Hence, in any case, $(\alpha, 2 \alpha) \in \Lambda_\Omega$. 

Now, we prove \emph{(i)} in the case $\mathbf{u} \neq 0$, using the case $\mathbf{u} = 0$ and Picard's fixed point theorem in $X_\epsilon$, for $\epsilon > 0$ sufficiently small. Consider $0 < \epsilon < 1$, $\epsilon \leq T$. Set 
\[\delta = \left\Vert \left( u^0, u^1 \right) \right\Vert_{H_0^1(\Omega) \times L^2(\Omega)} + \left\Vert g \right\Vert_{L^1((0, T), L^2)}.\]
For $U \in X_\epsilon$, write $u = L(U)$ for the solution of 
\[\left \{
\begin{array}{rcccl}
\square u + \beta u & = & f^\prime(\mathbf{u}) U + g & \quad & \text{in } (0, \epsilon) \times \Omega, \\
(u(0), \partial_t u(0)) & = & \left( u^0, u^1 \right) & \quad & \text{in } \Omega, \\
u & = & 0 & \quad & \text{on } (0, \epsilon) \times \partial \Omega.
\end{array}
\right.\]
Using the case $\mathbf{u} = 0$, one finds
\[\left\Vert L(U) \right\Vert_{X_\epsilon} \lesssim \left\Vert \left( u^0, u^1 \right) \right\Vert_{H_0^1(\Omega) \times L^2(\Omega)} + \left\Vert f^\prime(\mathbf{u}) U + g \right\Vert_{L^1((0, \epsilon), L^2)}.\]
Lemma \ref{lem_basic_estimate_f_local}-\emph{(iii)} gives
\[\left\Vert f^\prime(\mathbf{u}) U \right\Vert_{L^1((0, \epsilon), L^2)} \leq C \left\Vert U \right\Vert_{X_\epsilon} \left( \epsilon + \left\Vert \mathbf{u} \right\Vert_{L^{\alpha}((0, \epsilon), L^{2 \alpha})}^{\alpha - 1} \right),\]
yielding
\begin{equation}\label{proof_prop_time_dependent_potential_eq1}
\left\Vert L(U) \right\Vert_{X_\epsilon} \leq C \left( \delta + \left( \epsilon + \left\Vert \mathbf{u} \right\Vert_{L^{\alpha}((0, \epsilon), L^{2 \alpha})}^{\alpha - 1} \right) \left\Vert U \right\Vert_{X_\epsilon}\right)
\end{equation}
for some $C > 0$.

Consider $U, \tilde{U} \in X_\epsilon$ and write $u = L(U)$ and $\tilde{u} = L\left(\tilde{U}\right)$. One has
\[\left \{
\begin{array}{rcccl}
\square (u - \tilde{u}) + (u - \tilde{u}) & = & f^\prime(\mathbf{u}) (U - \tilde{U}) & \quad & \text{in } (0, \epsilon) \times \Omega, \\
((u - \tilde{u})(0), \partial_t (u - \tilde{u})(0)) & = & 0 & \quad & \text{in } \Omega, \\
u - \tilde{u} & = & 0 & \quad & \text{on } (0, \epsilon) \times \partial \Omega.
\end{array}
\right.\]
Hence, the case $\mathbf{u} = 0$ gives
\[\left\Vert L(U) - L\left( \tilde{U} \right) \right\Vert_{X_\epsilon} \lesssim \left\Vert f^\prime(\mathbf{u}) (U - \tilde{U}) \right\Vert_{L^1((0, \epsilon), L^2)},\]
and as above, it implies
\begin{equation}\label{proof_prop_time_dependent_potential_eq2}
\left\Vert L(U) - L\left( \tilde{U} \right) \right\Vert_{X_\epsilon} \leq C^\prime \left( \epsilon + \left\Vert \mathbf{u} \right\Vert_{L^{\alpha}((0, \epsilon), L^{2 \alpha})}^{\alpha - 1} \right) \left\Vert U - \tilde{U} \right\Vert_{X_\epsilon}.
\end{equation}

The constants in (\ref{proof_prop_time_dependent_potential_eq1}) and (\ref{proof_prop_time_dependent_potential_eq2}) do not depend on $\epsilon$, and up to increasing $C$ or $C^\prime$, we can assume that $C = C^\prime$. To apply Picard's fixed point theorem in a ball of radius $R > 0$ in $X_\epsilon$, one needs
\[\left \{
\begin{array}{c}
C \left( \delta + \left( \epsilon + \left\Vert \mathbf{u} \right\Vert_{L^{\alpha}((0, \epsilon), L^{2 \alpha})}^{\alpha - 1} \right) R \right) \leq R .\\
C \left( \epsilon + \left\Vert \mathbf{u} \right\Vert_{L^{\alpha}((0, \epsilon), L^{2 \alpha})}^{\alpha - 1} \right) < 1
\end{array}
\right.\]
We choose $\epsilon$ such that $C \left( \epsilon + \left\Vert \mathbf{u} \right\Vert_{L^{\alpha}((0, \epsilon), L^{2 \alpha})}^{\alpha - 1} \right) \leq \frac{1}{2}$. Then, we can simply chose $R = 2 C \delta$. By Picard's fixed point theorem, the solution $u$ is constructed on $[0, \epsilon]$, and one has
\[\Vert u \Vert_{X_\epsilon} \leq 2 C \left( \left\Vert \left( u^0, u^1 \right) \right\Vert_{H_0^1(\Omega) \times L^2(\Omega)} + \left\Vert g \right\Vert_{L^1((0, \epsilon), L^2)} \right).\] 

In particular, this implies
\[\left\Vert (u(\epsilon), \partial_t u(\epsilon)) \right\Vert_{H_0^1(\Omega) \times L^2(\Omega)} \leq 2 C \left( \left\Vert \left( u^0, u^1 \right) \right\Vert_{H_0^1(\Omega) \times L^2(\Omega)} + \left\Vert g \right\Vert_{L^1((0, T), L^2)} \right).\] 
Hence, this process can be iterated to construct the solution $u$ on $[\epsilon, 2 \epsilon]$. After a finite number of iterations, the solution $u$ is constructed in the space $X_T$, and satisfies
\[\Vert u \Vert_{X_T} \leq (2C)^{m} \left( \left\Vert \left( u^0, u^1 \right) \right\Vert_{H_0^1(\Omega) \times L^2(\Omega)} + \Vert g \Vert_{L^1((0, T), L^2)}\right)\]
for some $m \in \mathbb{N}$. This completes the proof of \emph{(i)}.

Now, we prove \emph{(ii)}. By a basic time-translation, we can assume that $\left[ T_1, T_2 \right] = [0, T]$. Note that both $\alpha_0$ and $\alpha_1$ satisfy the conditions of $\alpha$ in \emph{(i)}, so that \emph{(ii)} in the case $\mathbf{u} = 0$ is a direct consequence of \emph{(i)} in the case $\mathbf{u} = 0$. To prove that \emph{(ii)} is a consequence of \emph{(ii)} in the case $\mathbf{u} = 0$, one argue as above, by constructing the solution in $Y_{[0, \epsilon]}$ if $\epsilon > 0$ is sufficiently small. For $U \in Y_{[0, \epsilon]}$, by Lemma \ref{lem_basic_estimate_f_scatt}-\emph{(iii)}, one has
\[\left\Vert f^\prime(\mathbf{u}) U \right\Vert_{L^1((0, \epsilon), L^2)} \lesssim \left( \left\Vert \mathbf{u} \right\Vert_{L^{\alpha_0}((0, \epsilon), L^{2 \alpha_0})}^{\alpha_0 - 1} + \left\Vert \mathbf{u} \right\Vert_{L^{\alpha_1}((0, \epsilon), L^{2 \alpha_1})}^{\alpha_1 - 1} \right) \left\Vert U \right\Vert_{Y_{[0, \epsilon]}}.\]
Set $\eta(\epsilon) = \left\Vert \mathbf{u} \right\Vert_{L^{\alpha_0}((0, \epsilon), L^{2 \alpha_0})}^{\alpha_0 - 1} + \left\Vert \mathbf{u} \right\Vert_{L^{\alpha_1}((0, \epsilon), L^{2 \alpha_1})}^{\alpha_1 - 1}$. To apply Picard's fixed point theorem in a ball of radius $R > 0$ in $Y_{[0, \epsilon]}$, one needs
\[\left \{
\begin{array}{c}
C \left( \delta + \eta(\epsilon) R \right) \leq R ,\\
C \eta(\epsilon) < 1
\end{array}
\right.\]
where $\delta > 0$ is defined as above. If $\epsilon$ is sufficiently small, then one has $C \eta(\epsilon) \leq \frac{1}{2}$. If $R = 2 C \delta$, then the previous conditions are satisfied. By Picard's fixed point theorem, the solution $u$ is constructed on $[0, \epsilon]$. One can iterate this process as above.

Finally, we prove \emph{(iii)}. The case $V = 0$ is well-known, and as above, using Picard's fixed point theorem, we prove that it implies the case $V \neq 0$. Write $Z_\epsilon = \mathscr{C}^0([0, \epsilon], L^2(\Omega)) \cap \mathscr{C}^1([0, \epsilon], H^{-1}(\Omega))$. By (\ref{eq_prop_LKG_time_dependent_potential_1}), one has 
\[\left\Vert V U \right\Vert_{L^1((0, \epsilon), H^{-1})} \lesssim \epsilon \left\Vert U \right\Vert_{L^\infty((0, \epsilon), H^{-1})} \lesssim \epsilon \Vert U \Vert_{Z_\epsilon}, \quad U \in Z_\epsilon.\]
Using this estimate, the rest of the proof of \emph{(iii)} is similar to the proof of \emph{(i)}.
\end{proof}

\begin{rem}\label{rem_const_temps_petits}
Note that if $C(T_1, T_2)$ denotes the constant of \emph{(ii)}, then one can assume that $C(T_1, T_2) \leq C(T_1, T_1 + 1)$, for $T_1 \leq T_2 \leq T_1 + 1$. 
\end{rem}

We recall the local Cauchy theory for (\ref{eq_NLKG}).

\begin{thm}\label{thm_existence_nL_waves}
Consider $f$ satisfying \textnormal{(\ref{eq_def_nonlinearity_local})}. For any (real-valued) initial data $\left( u^0, u^1 \right) \in H_0^1(\Omega) \times L^2(\Omega)$, there exist a maximal time of existence $T \in (0, + \infty]$ and a unique solution $u$ of \textnormal{(\ref{eq_NLKG})} in $\mathscr{C}^0([0, T), H_0^1(\Omega)) \cap \mathscr{C}^1([0, T), L^2(\Omega))$. If $T < + \infty$, then 
\[\left\Vert \left( u(t), \partial_t u(t) \right) \right\Vert_{H_0^1(\Omega) \times L^2(\Omega)} \xrightarrow{t \rightarrow T^-} + \infty.\]
For $T^\prime < T$, if $f$ satisfies \textnormal{(\ref{eq_def_nonlinearity_local})} for some $\alpha$, then $u \in L^{\alpha}((0, T^\prime), L^{2 \alpha}(\Omega))$, and if $f$ satisfies \textnormal{(\ref{eq_def_nonlinearity_scatt})} for some $\alpha_0 \leq \alpha_1$, then 
\begin{equation}\label{eq_thm_existence_nL_waves_1}
u \in L^{\alpha_0}((0, T^\prime), L^{2 \alpha_0}(\Omega)) \, \cap \, L^{\alpha_1}((0, T^\prime), L^{2 \alpha_1}(\Omega)).
\end{equation}
\end{thm}

\begin{proof}
We only recall that the solution exists on $[0, \epsilon]$ if $\epsilon$ is sufficiently small, in the case $d \geq 3$. Set $q = 2 \alpha$. By (\ref{proof_prop_time_dep_potential_eq_for_later}), there exists $p > \alpha$ such that $(p, q) \in \Lambda_\Omega$. One has
\begin{equation}\label{proof_thm_existence_nL_waves_1}
\left\Vert u \right\Vert_{L^\alpha((0, \epsilon), L^{2 \alpha})} \leq \epsilon^\theta \left\Vert u \right\Vert_{L^p((0, \epsilon), L^q)}, \quad u \in L^p((0, \epsilon), L^q(\Omega)),
\end{equation}
for some $\theta > 0$. Write
\[\tilde{X}_\epsilon = \mathscr{C}^0([0, \epsilon], H_0^1(\Omega)) \ \cap \ \mathscr{C}^1([0, \epsilon], L^2(\Omega)) \ \cap \ L^p((0, \epsilon), L^q(\Omega)).\]
Note that (\ref{proof_thm_existence_nL_waves_1}) implies that $\left\Vert u \right\Vert_{X_\epsilon} \lesssim \left\Vert u \right\Vert_{\tilde{X}_\epsilon}$ for $u \in \tilde{X}_\epsilon$, with a constant independent of $\epsilon \in (0, 1)$. Together with Lemma \ref{lem_basic_estimate_f_local}-\emph{(i)} and (\ref{proof_thm_existence_nL_waves_1}), this gives
\[\left\Vert f(u) - f(v) \right\Vert_{L^1((0, \epsilon), L^2)} \lesssim \left\Vert u - v \right\Vert_{\tilde{X}_\epsilon} \left( \epsilon + \epsilon^{\theta (\alpha - 1)} \left\Vert u \right\Vert_{\tilde{X}_\epsilon}^{\alpha - 1} + \epsilon^{\theta (\alpha - 1)} \left\Vert v \right\Vert_{\tilde{X}_\epsilon}^{\alpha - 1} \right), \quad u, v \in \tilde{X}_\epsilon.\]
Using this estimate, Theorem \ref{thm_local_strichartz_estimates}, and Picard's fixed point theorem, one can construct the solution in $\tilde{X}_\epsilon$ if $\epsilon \in (0, 1)$ is sufficiently small. Note that if $f$ satisfies (\ref{eq_def_nonlinearity_scatt}) for some $\alpha_0 \leq \alpha_1$, then (\ref{proof_thm_existence_nL_waves_1}) holds for $\alpha_0$ and $\alpha_1$, implying (\ref{eq_thm_existence_nL_waves_1}).
\end{proof}

\section{Local controllability around a trajectory}

The proof of Theorem \ref{thm_main_local_control_trajectory} is organized as follows. First, we show that local controllability around $\mathbf{u}$ can be reduced to local controllability around $0$ for a modified nonlinear equation, and we prove that the solution of the controlled equation exists if the control is sufficiently small. Next, we establish an exact controllability result for the linearized equation. Finally, we complete the proof of local controllability using a fixed-point argument.

Consider $f$, $a$, $T$ and $\mathbf{u}$ satisfying the assumptions of Theorem \ref{thm_main_local_control_trajectory}. 

\subsection{The linearized equation}

Local controllability around $\mathbf{u}$ can be reformulated as follows. Consider $g \in L^1((0, T), L^2(\Omega))$, and denote by $u$ the solution of
\[\left \{
\begin{array}{rcccl}
\square u + \beta u & = & f(u) + ag & \quad & \text{in } (0, T) \times \Omega, \\
(u(T), \partial_t u(T)) & = & (\mathbf{u}(T), \partial_t \mathbf{u}(T)) & \quad & \text{in } \Omega, \\
u & = & 0 & \quad & \text{on } (0, T) \times \partial \Omega.
\end{array}
\right.\]
We prove below that $u$ exists on $[0, T]$ if $g$ is sufficiently small. Set $h = u - \mathbf{u}$. Then $(u(0), \partial_t u(0)) = \left( u^0, u^1 \right) \in H_0^1(\Omega) \times L^2(\Omega)$ if and only if $(h(0), \partial_t h (0)) = \left( u^0, u^1 \right) - (\mathbf{u}^0, \mathbf{u}^1)$, and $h$ solves
\begin{equation}\label{eq_linearized_h}
\left \{
\begin{array}{rcccl}
\square h + \beta h & = & f^\prime(\mathbf{u}) h + \NL_\mathbf{u}(h) + ag & \quad & \text{in } (0, T) \times \Omega, \\
(h(T), \partial_t h(T)) & = & (0, 0) & \quad & \text{in } \Omega, \\
h & = & 0 & \quad & \text{on } (0, T) \times \partial \Omega,
\end{array}
\right.
\end{equation}
where $\NL_\mathbf{u}(h)$ is defined by (\ref{eq_def_NL}). Hence, local controllability around $\mathbf{u}$ is equivalent to local controllability around zero for this modified equation. We prove that $h$ (and so $u$) exists on $[0, T]$ if $g$ is sufficiently small.

\begin{lem}
If $g \in L^1((0, T), L^2(\Omega))$ is sufficiently small, then the solution $h$ of \textnormal{(\ref{eq_linearized_h})} is well-defined on $[0, T]$.
\end{lem}

\begin{proof}
We use Picard's fixed point theorem in $X_T$ (defined by (\ref{eq_def_X_T})). For $H \in X_T$, write $h = L(H)$ for the solution of 
\[\left \{
\begin{array}{rcccl}
\square h + \beta h & = & f^\prime(\mathbf{u}) h + \NL_\mathbf{u}(H) + ag & \quad & \text{in } (0, T) \times \Omega, \\
(h(T), \partial_t h(T)) & = & 0 & \quad & \text{in } \Omega, \\
h & = & 0 & \quad & \text{on } (0, T) \times \partial \Omega.
\end{array}
\right.\]

By Proposition \ref{prop_LKG_time_dependent_potential}-\emph{(i)}, one has
\[\left\Vert L(H) \right\Vert_{X_T} \lesssim \left\Vert \NL_\mathbf{u}(H) + ag \right\Vert_{L^1((0, T), L^2)} \lesssim \left\Vert \NL_\mathbf{u}(H) \right\Vert_{L^1((0, T), L^2)} + \left\Vert g \right\Vert_{L^1((0, T), L^2)},\]
and by Lemma \ref{lem_basic_estimate_f_local}-\emph{(ii)}, this gives
\[\left\Vert L(H) \right\Vert_{X_T} \lesssim \left\Vert H \right\Vert_{X_T}^2 + \left\Vert H \right\Vert_{X_T}^{\alpha} + \left\Vert g \right\Vert_{L^1((0, T), L^2)}.\]
Similarly, for $H, \tilde{H} \in X_T$, one finds
\[\left\Vert L(H) - L(\tilde{H}) \right\Vert_{X_T} \lesssim \left\Vert H - \tilde{H} \right\Vert_{X_T} \left( \left\Vert H \right\Vert_{X_T} + \left\Vert \tilde{H} \right\Vert_{X_T} + \left\Vert H \right\Vert_{X_T}^{\alpha - 1} + \left\Vert \tilde{H} \right\Vert_{X_T}^{\alpha - 1} \right).\]

To apply Picard's fixed point theorem in a ball of radius $R \in (0, 1)$ in $X_T$, as $\alpha > 2$, one needs
\[\left \{
\begin{array}{c}
C \left( \left\Vert g \right\Vert_{L^1((0, T), L^2)} + R^2 \right) \leq R \\
C R < 1
\end{array}
\right.\]
for some $C > 0$. We choose $R = 2 C \left\Vert g \right\Vert_{L^1((0, T), L^2)}$. If $\left\Vert g \right\Vert_{L^1((0, T), L^2)}$ is sufficiently small, then the previous conditions are satisfied. This completes the proof.
\end{proof}

\subsection{Exact controllability for the linearized equation}

Here, we prove the following exact controllability result.

\begin{prop}\label{prop_control_time_dep_potential}
There exists a continuous linear operator 
\[\begin{array}{cccc} 
\mathsf{G}: & H_0^1(\Omega) \times L^2(\Omega) & \longrightarrow & L^1((0, T), L^2(\Omega)) .\\
& \left( u^0, u^1 \right) & \longmapsto & g\left( u^0, u^1 \right)
\end{array}\]
such that for $\left( u^0, u^1 \right) \in H_0^1(\Omega) \times L^2(\Omega)$, the solution of 
\begin{equation}\label{eq_prop_control_time_dep_potential}
\left \{
\begin{array}{rcccl}
\square u + \beta u & = & f^\prime(\mathbf{u}) u + a \mathsf{G}\left( u^0, u^1 \right) & \quad & \text{in } (0, T) \times \Omega, \\
(u(T), \partial_t u(T)) & = & (0, 0) & \quad & \text{in } \Omega, \\
u & = & 0 & \quad & \text{on } (0, T) \times \partial \Omega,
\end{array}
\right.
\end{equation}
satisfies $(u(0), \partial_t u(0)) = \left( u^0, u^1 \right)$.
\end{prop}

\begin{rem}
The two main difficulties of Proposition \ref{prop_control_time_dep_potential} are the following.
\begin{enumerate}[label=(\roman*)]
\item As the potential $f^\prime(\mathbf{u})$ is time-dependent, we need unique continuation for wave equations with partially analytic coefficients, to prove that there is no nonzero solution of the dual equation which is equal to zero on the support of $a$. 
\item As the domain $\Omega$ may be unbounded, the embedding $H^1(\Omega) \hookrightarrow L^2(\Omega)$ may fail to be compact. In that case, our proof relies on the assumptions that $\left\vert \nabla f^\prime(\mathbf{u}(t, x)) \right\vert + \left\vert f^\prime(\mathbf{u}(t, x))\right\vert \longrightarrow 0$ as $|x| \rightarrow \infty$, and that $a \geq c > 0$ on the complement of a bounded region.
\end{enumerate}
\end{rem}

\begin{proof} 
We show that the operator
\[\begin{array}{cccc} 
L: & L^2((0, T) \times \Omega) & \longrightarrow & H_0^1(\Omega) \times L^2(\Omega) .\\
& g & \longmapsto & (u(0), \partial_t u(0))
\end{array}\]
is onto, where $u$ is the solution of 
\begin{equation*}
\left \{
\begin{array}{rcccl}
\square u + \beta u & = & f^\prime(\mathbf{u}) u + a g & \quad & \text{in } (0, T) \times \Omega, \\
(u(T), \partial_t u(T)) & = & (0, 0) & \quad & \text{in } \Omega, \\
u & = & 0 & \quad & \text{on } (0, T) \times \partial \Omega.
\end{array}
\right.
\end{equation*}

\paragraph{Step 1: the dual problem.}
For the duality between $H_0^1(\Omega) \times L^2(\Omega)$ and $L^2(\Omega) \times H^{-1}(\Omega)$, we choose 
\[\left\langle \left( v^0, v^1 \right), \left( u^0, u^1 \right) \right\rangle_{L^2(\Omega) \times H^{-1}(\Omega), H_0^1(\Omega) \times L^2(\Omega)} = \left\langle v^1, u^0 \right\rangle_{H^{-1}(\Omega), H_0^1(\Omega)} - \left\langle v^0, u^1 \right\rangle_{L^2(\Omega)}.\]
Fix $\left( v^0, v^1 \right) \in L^2(\Omega) \times H^{-1}(\Omega)$. By definition, $L^\ast \left( v^0, v^1 \right)$ satisfies
\[\left\langle L^\ast\left( v^0, v^1 \right), g \right\rangle_{L^2((0, T) \times \Omega)} = \left\langle \left( v^0, v^1 \right), L(g) \right\rangle_{L^2(\Omega) \times H^{-1}(\Omega), H_0^1(\Omega) \times L^2(\Omega)}, \quad g \in L^2((0, T) \times \Omega).\]
In particular, if $v$ is a function such that $(v(0), \partial_t v(0)) = \left( v^0, v^1 \right)$, then one has 
\[\left\langle L^\ast\left( v^0, v^1 \right), g \right\rangle_{L^2((0, T) \times \Omega)} = \int_0^T \partial_t \left( - \left\langle \partial_t v(t), u(t) \right\rangle_{H^{-1}(\Omega), H_0^1(\Omega)} + \left\langle v(t), \partial_t u(t) \right\rangle_{L^2(\Omega)} \right) \mathrm{d}t\]
for all $g \in L^2((0, T) \times \Omega)$. If, in addition, $v$ is smooth, then 
\begin{align*}
& \left\langle L^\ast\left( v^0, v^1 \right), g \right\rangle_{L^2((0, T) \times \Omega)} \\
= & \int_0^T \left( \left\langle - \partial_t^2 v(t) + \Delta v(t) - \beta v(t) + f^\prime(\mathbf{u}(t)) v(t), u(t) \right\rangle_{L^2(\Omega)} + \left\langle a v(t), g(t) \right\rangle_{L^2(\Omega)} \right) \mathrm{d}t
\end{align*}
for all $g \in L^2((0, T) \times \Omega)$. This shows that for $\left( v^0, v^1 \right) \in \mathscr{C}^\infty_\mathrm{c}(\Omega)^2$, one has
\begin{equation}\label{proof_prop_exact_control_eq1}
L^\ast\left( v^0, v^1 \right) = a v
\end{equation}
where $v$ is the solution of
\begin{equation}\label{proof_prop_exact_control_eq1bis}
\left \{
\begin{array}{rcccl}
\square v + \beta v & = & f^\prime(\mathbf{u}) v & \quad & \text{in } (0, T) \times \Omega, \\
(v(0), \partial_t v(0)) & = & \left( v^0, v^1 \right) & \quad & \text{in } \Omega, \\
v & = & 0 & \quad & \text{on } (0, T) \times \partial \Omega.
\end{array}
\right.
\end{equation}
By definition, $L^\ast$ is a continuous operator from $L^2(\Omega) \times H^{-1}(\Omega)$ to $L^2((0, T) \times \Omega)$. Thus, for $\left( v^0, v^1 \right) \in L^2(\Omega) \times H^{-1}(\Omega)$, one has
\[L^\ast \left( v^0, v^1 \right) = \lim_{n \rightarrow \infty} L^\ast \left( v_n^0, v_n^1 \right)\]
where $\left( \left( v_n^0, v_n^1 \right) \right)_{n \in \mathbb{N}}$ is a sequence of elements of $\mathscr{C}^\infty_\mathrm{c}(\Omega)^2$ converging to $\left( v^0, v^1 \right)$ in $L^2(\Omega) \times H^{-1}(\Omega)$. This proves that (\ref{proof_prop_exact_control_eq1}) holds for all $\left( v^0, v^1 \right) \in L^2(\Omega) \times H^{-1}(\Omega)$, where $v$ is the solution of (\ref{proof_prop_exact_control_eq1bis}) given by Proposition \ref{prop_LKG_time_dependent_potential}-\emph{(iii)}.

\paragraph{Step 2: a compactness property.}
Write $L^\ast = A + K$, where $A$ is the operator from $L^2(\Omega) \times H^{-1}(\Omega)$ to $L^2((0, T) \times \Omega)$ defined by $A\left( v^0, v^1 \right) = a \phi$, where $\phi$ is the solution of
\[\left \{
\begin{array}{rcccl}
\square \phi + \beta \phi & = & 0 & \quad & \text{in } (0, T) \times \Omega, \\
(\phi(0), \partial_t \phi(0)) & = & \left( v^0, v^1 \right) & \quad & \text{in } \Omega, \\
\phi & = & 0 & \quad & \text{on } (0, T) \times \partial \Omega.
\end{array}
\right.\]
By definition, $K\left( v^0, v^1 \right) = aw$, where $w$ is the solution of 
\[\left \{
\begin{array}{rcccl}
\square w + \beta w & = & f^\prime(\mathbf{u}) v & \quad & \text{in } (0, T) \times \Omega, \\
(w(0), \partial_t w(0)) & = & 0 & \quad & \text{in } \Omega, \\
w & = & 0 & \quad & \text{on } (0, T) \times \partial \Omega,
\end{array}
\right.\]
and $v$ is the solution of (\ref{proof_prop_exact_control_eq1bis}).

We show that $K$ is compact. Let $\left( \left( v_n^0, v_n^1 \right) \right)_{n \in \mathbb{N}}$ be a bounded sequence of elements of $L^2(\Omega) \times H^{-1}(\Omega)$. We want to show that there exists a subsequence of $\left( K\left( v_n^0, v_n^1 \right) \right)_{n \in \mathbb{N}}$ which converges in $L^2((0, T) \times \Omega)$. If $\Omega$ is compact, then the proof is a consequence of Rellich's theorem. Indeed, in that case, we can assume that $\left( \left( v_n^0, v_n^1 \right) \right)_{n \in \mathbb{N}}$ converges in $H^{-1}(\Omega) \times H^{-2}(\Omega)$ up to a subsequence. One has 
\[\Vert a w_n - a w_m \Vert_{L^2((0, T) \times \Omega)} \lesssim \left\Vert f^\prime(\mathbf{u}) (v_n - v_m) \right\Vert_{L^1((0, T), H^{-1})} \lesssim \left\Vert \left( v_n^0, v_n^1 \right) - (v_m^0, v_m^1) \right\Vert_{L^2(\Omega) \times H^{-1}(\Omega)},\]
implying that the sequence $\left( a w_n \right)_{n \in \mathbb{N}}$ converges.

Now, assume that $\Omega$ is not compact. We use the following extension of Rellich's theorem.

\begin{lem}\label{lem_rellich_extension}
Consider $U$ a (possibly empty) smooth bounded open subset of $\mathbb{R}^d$ and $s \in \mathbb{R}$. Let $V \in \mathscr{C}^\infty(\mathbb{R}^d \backslash U)$ be such that 
\[\sum_{\vert \beta \vert \leq \vert s - 1 \vert} \left\vert \partial_x^\beta V(x) \right\vert \xrightarrow{|x| \rightarrow \infty} 0.\]
Then the operator
\[\begin{array}{ccc} 
H^s(\mathbb{R}^d \backslash U) & \longrightarrow & H^{s-1}(\mathbb{R}^d \backslash U)\\
u & \longmapsto & Vu
\end{array}\]
is compact.
\end{lem}

A proof of Lemma \ref{lem_rellich_extension} can be found in Appendix \ref{appendix_proof_rellich}. We apply Ascoli's theorem to the sequence $\left( f^\prime(\mathbf{u}) v_n \right)_{n \in \mathbb{N}}$. For all $n \in \mathbb{N}$, one has $f^\prime(\mathbf{u}) v_n \in \mathscr{C}^0([0, T], H^{-1}(\Omega))$, and 
\[\left\Vert \partial_t \left( f^\prime(\mathbf{u}) v_n \right) \right\Vert_{L^\infty([0, T], H^{-1})} \lesssim \left\Vert \left( v_n^0, v_n^1 \right) \right\Vert_{L^2(\Omega) \times H^{-1}(\Omega)}.\]
Hence, the sequence $\left( f^\prime(\mathbf{u}) v_n \right)_{n \in \mathbb{N}}$ is equicontinuous. Applying Lemma \ref{lem_rellich_extension} with $s = 0$, one finds that for all $t \in [0, T]$, the set 
\[\left\{ f^\prime(\mathbf{u}(t)) v_n(t), n\in \mathbb{N} \right\}\]
is relatively compact in $H^{-1}(\Omega)$. Hence, by Ascoli's theorem, the sequence $\left( f^\prime(\mathbf{u}) v_n \right)_{n \in \mathbb{N}}$ converges in $L^\infty([0, T], H^{-1}(\Omega))$, up to a subsequence. Then, as 
\[\Vert a w_n - a w_m \Vert_{L^2((0, T) \times \Omega)} \lesssim \left\Vert f^\prime(\mathbf{u}) (v_n - v_m) \right\Vert_{L^1((0, T), H^{-1})} \lesssim \left\Vert f^\prime(\mathbf{u}) (v_n - v_m) \right\Vert_{L^\infty((0, T), H^{-1})},\]
one finds that the sequence $\left( a w_n \right)_{n \in \mathbb{N}}$ converges. Hence, $K$ is compact.

\paragraph{Step 3: observability for a wave equation with constant coefficients.}

One has 
\[\left\Vert \left( v^0, v^1 \right) \right\Vert_{L^2(\Omega) \times H^{-1}(\Omega)} \lesssim \left\Vert A \left( v^0, v^1 \right) \right\Vert_{L^2((0, T) \times \Omega)},\]
by the following theorem.

\begin{thm}\label{thm_obs_constant_coeff}
Assume that there exist $\omega \subset \Omega$ and $c > 0$ such that $a \geq c$ on $\omega$, and such that $(\omega, T)$ satisfies the GCC. In addition, if $\Omega$ is unbounded, assume that there exists $R_0 > 0$ such that $\mathbb{R}^d \backslash B(0, R_0) \subset \omega$. Then, there exists $C > 0$ such that for all $\left( u^0, u^1 \right) \in L^2(\Omega) \times H^{-1}(\Omega)$, the solution $u$ of
\[\left \{
\begin{array}{rcccl}
\square u + \beta u & = & 0 & \quad & \text{in } (0, T) \times \Omega, \\
(u(0), \partial_t u(0)) & = & \left( u^0, u^1 \right) & \quad & \text{in } \Omega, \\
u & = & 0 & \quad & \text{on } (0, T) \times \partial \Omega,
\end{array}
\right.\]
given by Proposition \ref{prop_LKG_time_dependent_potential}-\emph{(iii)} (with $V = 0$), satisfies
\[\left\Vert \left( u^0, u^1 \right) \right\Vert_{L^2(\Omega) \times H^{-1}(\Omega)} \leq C\left\Vert a u \right\Vert_{L^2((0, T) \times \Omega)}.\]
\end{thm}

In the case of a compact domain $\Omega$, it is well-known that the GCC implies Theorem \ref{thm_obs_constant_coeff}, since the work of Bardos, Leabeau and Rauch (see \cite{BLR}, Theorem 3.8). If $\Omega$ is not compact, we provide a proof of Theorem \ref{thm_obs_constant_coeff} in Appendix \ref{appendix_proof_obs_unbounded}. A stabilisation property in a similar context can be found in \cite{Joly-Laurent}. 

\paragraph{Step 4: invisible solutions of the dual of the linearized equation.}
In that step, we prove that the operator $L^\ast$ is one-to-one. Let $\left( v^0, v^1 \right) \in L^2(\Omega) \times H^{-1}(\Omega)$ be such that $L^\ast \left( v^0, v^1 \right) = av = 0$. One has $v(t) = 0$ on $\omega$, for all $t \in [0, T]$, and by assumption, $(\omega, T)$ fulfils the GCC. By the theorem of propagation of singularities of Melrose and Sjöstrand (see \cite{MelroseSjostrand}), $v$ is smooth. In particular, we can use Theorem 6.1 of \cite{LaurentLeautaud}, which we copy here for convenience. We write $\dist$ for the geodesic distance on a Riemannian manifold $M$, and 
\[\dist(x_1, \omega) = \inf_{x_0 \in \omega} \dist(x_0, x_1), \quad x_1 \in M,\]
for the distance to a subset $\omega \subset M$.

\begin{thm}[Theorem 6.1 of \cite{LaurentLeautaud}]
Let $(M, g)$ be a compact Riemannian manifold with (or without) boundary and write $\Delta_g$ for the Laplace-Beltrami operator on $M$. Let $\omega$ be an open subset of $M$, and consider $T > 0$ such that
\begin{equation}\label{eq_thm_LaurentLeautaud_1}
T > \sup_{x_1 \in M} \dist(x_1, \omega).
\end{equation}
Set $P = \partial_t^2 - \Delta_g + V$, where $V \in \mathscr{C}^\infty([-T, T] \times \Omega)$ depends analytically on the variable $t$. There exist $C$, $\kappa$, $\mu_0 > 0$ such that for any $\left( u^0, u^1 \right) \in H^1_0(M) \times L^2(M)$, if $u$ is the solution of
\[\left \{
\begin{array}{rcccl}
P u & = & 0 & \quad & \text{in } (-T, T) \times M, \\
(u(0), \partial_t u(0)) & = & \left( u^0, u^1 \right) & \quad & \text{in } M, \\
u & = & 0 & \quad & \text{on } (-T, T) \times \partial M,
\end{array}
\right.\]
then for any $\mu \geq \mu_0$, one has
\[\left\Vert \left( u^0, u^1 \right) \right\Vert_{L^2(\Omega) \times H^{-1}(\Omega)} \leq C e^{\kappa \mu} \Vert u \Vert_{L^2((-T ,T), H^1(\omega))} + \frac{C}{\mu} \left\Vert \left( u^0, u^1 \right) \right\Vert_{H_0^1(\Omega) \times L^2(\Omega)}.\]
\end{thm}

If $\Omega$ is compact, then this theorem immediately gives $\left( v^0, v^1 \right) = 0$. Note that (\ref{eq_thm_LaurentLeautaud_1}) is a consequence of the fact that $(\omega, T)$ satisfies the GCC: see, for example, Lemma B.4 of \cite{LaurentLeautaudOBS}. If $\Omega$ is not compact, then by assumption there exists $R_0 > 0$ such that $a > 0$ on $\mathbb{R}^d \backslash B(0, R_0)$. Hence, $v$ is the solution of (\ref{proof_prop_exact_control_eq1bis}) on the compact domain $\Omega \cap B(0, R_0)$, and we can also apply the previous theorem. That proves that $L^\ast$ is one-to-one in all cases.

\paragraph{Step 5: conclusion.}

By Step 3, one has
\[\left\Vert \left( v^0, v^1 \right) \right\Vert_{L^2(\Omega) \times H^{-1}(\Omega)} \lesssim \left\Vert L^\ast \left( v^0, v^1 \right) \right\Vert_{L^2((0, T) \times \Omega)} + \left\Vert K \left( v^0, v^1 \right) \right\Vert_{L^2((0, T) \times \Omega)}.\]
By Step 2 and Step 4, $K$ is compact and $L^\ast$ is one-to-one. We apply the following classical result, which proves that $L$ is onto.

\begin{thm}
Let $X$, $Y$ and $Z$ be Hilbert spaces, and $L: X \rightarrow Y$ and $K: Y \rightarrow Z$ be linear continuous operators. Assume that $L^\ast$ is one-to-one and $K$ is compact, and that there exists $C > 0$ such that
\[\Vert y \Vert_{Y} \leq C \left( \left\Vert L^\ast y \right\Vert_{X} + \left\Vert K y \right\Vert_{Z} \right), \quad y \in Y.\]
Then, $L$ is onto.
\end{thm}

\begin{proof}
By Lemma \ref{lem_general_compact_term_2} in Appendix \ref{appendix_lemme_compact_term}, one has 
\[\Vert y \Vert_{Y} \lesssim \Vert L^\ast y \Vert_{X}, \quad y \in Y.\]
A standard argument of functional analysis implies that $L$ is onto (see \cite{LeRousseau}, Corollary 11.20, for example).
\end{proof}

Then, for any linear subspace $E$ of $L^2((0, T) \times \Omega)$ such that
\[E \oplus \Ker L = L^2((0, T) \times \Omega),\]
the operator $\mathsf{G}$ can be constructed as a continuous linear operator from $H_0^1(\Omega) \times L^2(\Omega)$ to $E$. This completes the proof of Proposition \ref{prop_control_time_dep_potential}.
\end{proof}

\subsection{Local controllability for the non-linear equation}

Here, we prove Theorem \ref{thm_main_local_control_trajectory}. Fix $(h^0, h^1) \in H_0^1(\Omega) \times L^2(\Omega)$. Let $X = X_T$ be the space defined by (\ref{eq_def_X_T}). Write $\mathsf{G}$ for the operator of Proposition \ref{prop_control_time_dep_potential}. For $H \in X$, write $\phi_H$ for the solution of 
\[\left \{
\begin{array}{rcccl}
\square \phi_H + \beta \phi_H & = & f^\prime(\mathbf{u}) \phi_H + \NL_\mathbf{u}(H) & \quad & \text{in } (0, T) \times \Omega, \\
(\phi_H(T), \partial_t \phi_H(T)) & = & 0 & \quad & \text{in } \Omega, \\
\phi_H & = & 0 & \quad & \text{on } (0, T) \times \partial \Omega,
\end{array}
\right.\]
given by Proposition \ref{prop_LKG_time_dependent_potential}-\emph{(i)}. We claim that the solution $h = \Gamma(H)$ of
\[\left \{
\begin{array}{rcccl}
\square h + \beta h & = & f^\prime(\mathbf{u}) h + \NL_\mathbf{u}(H) + a \mathsf{G} \left((h^0, h^1) - (\phi_H(0), \partial_t \phi_H(0))\right) & \quad & \text{in } (0, T) \times \Omega, \\
(h(T), \partial_t h(T)) & = & 0 & \quad & \text{in } \Omega, \\
h & = & 0 & \quad & \text{on } (0, T) \times \partial \Omega,
\end{array}
\right.\]
satisfies $(h(0), \partial_t h(0)) = (h^0, h^1)$. Indeed, $w = h - \phi_H$ solves
\[\left \{
\begin{array}{rcccl}
\square w + \beta w & = & f^\prime(\mathbf{u}) w + a \mathsf{G} \left((h^0, h^1) - (\phi_H(0), \partial_t \phi_H(0))\right) & \quad & \text{in } (0, T) \times \Omega, \\
(w(T), \partial_t w(T)) & = & 0 & \quad & \text{in } \Omega, \\
w & = & 0 & \quad & \text{on } (0, T) \times \partial \Omega,
\end{array}
\right.\]
implying $(h(0), \partial_t h(0)) - (\phi_H(0), \partial_t \phi_H(0)) = (w(0), \partial_t w(0)) = (h^0, h^1) - (\phi_H(0), \partial_t \phi_H(0))$, by definition of $\mathsf{G}$.

We show that if
\[\delta = \left\Vert (h^0, h^1) \right\Vert_{H_0^1(\Omega) \times L^2(\Omega)}\]
is sufficiently small, then $\Gamma$ has a unique fixed-point in a small neighbourhood of zero in $X$. By Proposition \ref{prop_LKG_time_dependent_potential}-\emph{(i)} (applied to $h = \Gamma(H)$), one has
\[\left\Vert \Gamma(H) \right\Vert_X \lesssim \left\Vert \NL_\mathbf{u}(H) + a \mathsf{G}\left((h^0, h^1) - (\phi_H(0), \partial_t \phi_H(0))\right) \right\Vert_{L^1((0, T), L^2)}.\]
Using the continuity of $\mathsf{G}$ (see Proposition \ref{prop_control_time_dep_potential}), one finds
\[\left\Vert \Gamma(H) \right\Vert_X \lesssim \left\Vert \NL_\mathbf{u}(H) \right\Vert_{L^1((0, T), L^2)} + \left\Vert (h^0, h^1) \right\Vert_{H_0^1(\Omega) \times L^2(\Omega)} + \left\Vert (\phi_H(0), \partial_t \phi_H(0)) \right\Vert_{H_0^1(\Omega) \times L^2(\Omega)}.\]
By Proposition \ref{prop_LKG_time_dependent_potential}-\emph{(i)} (applied to $\phi_H$), one has
\[\left\Vert (\phi_H(0), \partial_t \phi_H(0)) \right\Vert_{H_0^1(\Omega) \times L^2(\Omega)} \lesssim \left\Vert \NL_\mathbf{u}(H) \right\Vert_{L^1((0, T), L^2)}\]
implying
\[\left\Vert \Gamma(H) \right\Vert_X \lesssim \left\Vert \NL_\mathbf{u}(H) \right\Vert_{L^1((0, T), L^2)} + \left\Vert (h^0, h^1) \right\Vert_{H_0^1(\Omega) \times L^2(\Omega)}.\]
Thus, using Lemma \ref{lem_basic_estimate_f_local}-\emph{(ii)}, one obtains
\[\left\Vert \Gamma(H) \right\Vert_X \lesssim \Vert H \Vert_X^2 + \Vert H \Vert_X^{\alpha} + \delta.\]

Similarly, for $H, \tilde{H} \in X$, one has 
\begin{align*}
& \left\Vert \Gamma(H) - \Gamma(\tilde{H}) \right\Vert_X \\
\lesssim & \left\Vert \NL_\mathbf{u}(H) - \NL_\mathbf{u}(\tilde{H}) + a \mathsf{G}\left((\phi_H(0), \partial_t \phi_H(0)) - (\phi_{\tilde{H}}(0), \partial_t \phi_{\tilde{H}}(0))\right) \right\Vert_{L^1((0, T), L^2)} \\
\lesssim & \left\Vert \NL_\mathbf{u}(H)- \NL_\mathbf{u}(\tilde{H}) \right\Vert_{L^1((0, T), L^2)} + \left\Vert (\phi_H(0), \partial_t \phi_H(0)) - (\phi_{\tilde{H}}(0), \partial_t \phi_{\tilde{H}}(0)) \right\Vert_{H_0^1(\Omega) \times L^2(\Omega)} \\
\lesssim & \left\Vert \NL_\mathbf{u}(H)- \NL_\mathbf{u}(\tilde{H}) \right\Vert_{L^1((0, T), L^2)} \\
\lesssim & \left\Vert H - \tilde{H} \right\Vert_X \left( \left\Vert H \right\Vert_X + \left\Vert H \right\Vert_X^{\alpha - 1} + \left\Vert \tilde{H} \right\Vert_X + \left\Vert \tilde{H} \right\Vert_X^{\alpha - 1} \right).
\end{align*}

To apply Picard's fixed point theorem in a ball of radius $R > 0$ in $X$, one needs
\[\left \{
\begin{array}{c}
C \left( \delta + R^2 + R^{\alpha} \right) \leq R \\
C \left(R + R^{\alpha - 1} \right) < 1
\end{array}
\right.\]
where $C > 0$ is a constant. We choose $R = 2 C \delta$. As $\alpha > 1$, the previous conditions are satisfied if $\delta$ is sufficiently small. This completes the proof.

\section{Null-controllability of a scattering solution in a long time}

In this section, we prove Theorem \ref{thm_main_null_control}. In particular, we only consider $f$ satisfying (\ref{eq_def_nonlinearity_scatt}), implying that $\Omega$ is unbounded, and that $3 \leq d \leq 5$. Note also that this requires $\beta > 0$. The proof of Theorem \ref{thm_main_null_control} is organized as follows. First, we prove a local energy decay result for solutions of the linear equation. Second, we prove that together with local-in-time Strichartz estimates and global-in-time Strichartz estimates on $\mathbb{R}^d$, it implies global-in-time Strichartz estimates on $\Omega$ (Theorem \ref{thm_global_strichartz}). Finally, using local energy decay, global-in-time Strichartz estimates, and local controllability around zero, we prove Theorem \ref{thm_main_null_control}.

\subsection{Local energy decay}

Here, we prove the following result.

\begin{thm}\label{thm_local_decay_energy}
Assume that $\Omega$ is unbounded and non-trapping, with $d \geq 3$, and consider $\chi \in \mathscr{C}^\infty_\mathrm{c}(\Omega)$ and $R_0 > 0$. There exists $C > 0$ such that for all $\left( u^0, u^1 \right) \in H_0^1(\Omega) \times L^2(\Omega)$ and $F \in L^2(\mathbb{R} \times \Omega)$ supported in $\mathbb{R} \times \left( \Omega \cap B(0, R_0) \right)$, the solution $u$ of 
\begin{equation}\label{eq_thm_local_decay_energy}
\left \{
\begin{array}{rcccl}
\square u + \beta u & = & F & \quad & \text{in } \mathbb{R} \times \Omega, \\
(u(0), \partial_t u(0)) & = & \left( u^0, u^1 \right) & \quad & \text{in } \Omega, \\
u & = & 0 & \quad & \text{on } \mathbb{R} \times \partial \Omega,
\end{array}
\right.
\end{equation}
satisfies $\left( \chi u, \chi \partial_t u \right) \in L^2(\mathbb{R}, H_0^1(\Omega) \times L^2(\Omega))$, with
\[\left\Vert \left( \chi u, \chi \partial_t u \right) \right\Vert_{L^2(\mathbb{R}, H_0^1 \times L^2)} \leq C \left( \left\Vert \left( u^0, u^1 \right) \right\Vert_{H_0^1(\Omega) \times L^2(\Omega)} + \Vert F \Vert_{L^2(\mathbb{R} \times \Omega)} \right).\]
\end{thm}

\begin{proof}
The proof is based on \cite{Burq03} and \cite{BurqPersonnal}. There is a small mistake in the $TT^\ast$ argument in \cite{Burq03}: formula (2.6) is incorrect, because the operator 
\[\begin{array}{ccc} 
H^s(\Omega) & \longrightarrow & H^s(\Omega) \\
u & \longmapsto & \chi u
\end{array}\]
is not self-adjoint when $s \neq 0$. Carrying out the argument with the adjoint of this operator requires the use of more complicated resolvent estimates than those employed in \cite{Burq03}. Instead of doing that, we rely on \cite{BurqPersonnal}: we use two $TT^\ast$ arguments, at two different levels of regularity, and we conclude using interpolation.

We split the proof in 4 steps.

\paragraph{Step 1: a first $TT^\ast$ argument.} 
Here, we prove that 
\begin{equation}\label{proof_thm_local_decay_energy_eq_1}
\left\Vert \chi \partial_t u \right\Vert_{L^2(\mathbb{R} \times \Omega)} \leq C \left\Vert \left( u^0, u^1 \right) \right\Vert_{H_0^1(\Omega) \times L^2(\Omega)}, \quad \left( u^0, u^1 \right) \in H_0^1(\Omega) \times L^2(\Omega),
\end{equation}
where $u$ is the solution of (\ref{eq_thm_local_decay_energy}) with $F = 0$. Write $\mathcal{H} = H_0^1(\Omega) \times L^2(\Omega)$, which is a Hilbert space for the scalar product 
\begin{align*}
\left\langle \left( u^0, u^1 \right), \left( v^0, v^1 \right) \right\rangle_{\mathcal{H}} & = \left\langle u^0, v^0 \right\rangle_{H_0^1(\Omega)} + \left\langle u^1, v^1 \right\rangle_{L^2(\Omega)} \\
& = \left\langle \nabla u^0, \nabla v^0 \right\rangle_{L^2(\Omega)} + \beta \left\langle u^0, v^0 \right\rangle_{L^2(\Omega)} + \left\langle u^1, v^1 \right\rangle_{L^2(\Omega)},
\end{align*}
and $S(t): \mathcal{H} \rightarrow \mathcal{H}$ for the linear semi-group associated with (\ref{eq_thm_local_decay_energy}), of infinitesimal generator
\[A = \left(
\begin{array}{cc}
0 & \mathrm{Id} \\ 
\Delta - \beta & 0
\end{array}
\right): D(A) \subset \mathcal{H} \rightarrow \mathcal{H}, \quad D(A) = \left( H^2(\Omega) \cap H_0^1(\Omega) \right) \times H_0^1(\Omega).\]
By conservation of the energy, one has
\begin{equation}\label{proof_thm_local_decay_energy_eq_2}
\left\Vert S(t) \left( u^0, u^1 \right) \right\Vert_{\mathcal{H}} = \left\Vert \left( u^0, u^1 \right) \right\Vert_{\mathcal{H}}, \quad t \in \mathbb{R}, \quad \left( u^0, u^1 \right) \in \mathcal{H}.
\end{equation}
Denote by $\pi_1: \mathbb{R}^2 \rightarrow \mathbb{R}$ the projection on the second coordinate. For $t \in \mathbb{R}$, consider the linear continuous operator
\[\begin{array}{cccc} 
T(t): & \mathcal{H} & \longrightarrow & L^2(\Omega) \\
& \left( u^0, u^1 \right) & \longmapsto & \chi \partial_t u(t) = \chi \pi_1 S(t) \left( u^0, u^1 \right)
\end{array}.\]

We start with the computation of $T(t)^\ast$. One has
\[\left\langle T(t) \left( u^0, u^1 \right), v^1 \right\rangle_{L^2(\Omega)} = \left\langle \left( u^0, u^1 \right), S(t)^\ast \left( 0, \chi v^1 \right) \right\rangle_{\mathcal{H}}, \quad t \in \mathbb{R}, \quad v^1 \in L^2(\Omega).\]
By Corollary 10.6 of \cite{Pazy}, the adjoint semigroup of $S$ is a $C_0$-semigroup, generated by $A^\ast$. As all functions considered are real-valued, an integration by parts gives
\[\left\langle A \left( u^0, u^1 \right), \left( v^0, v^1 \right) \right\rangle_{\mathcal{H}} = \left\langle u^1, v^0 \right\rangle_{H_0^1(\Omega)} - \left\langle u^0, v^1 \right\rangle_{H_0^1(\Omega)} = - \left\langle \left( u^0, u^1 \right), A \left( v^0, v^1 \right) \right\rangle_{\mathcal{H}},\]
for $\left( u^0, u^1 \right), \left( v^0, v^1 \right) \in D(A)$, yielding $A^\ast = - A$. Hence, one finds
\[T(t)^\ast v^1 = S(- t) \left( 0, \chi v^1 \right), \quad t \in \mathbb{R}, \quad v^1 \in L^2(\Omega).\]

Fix $\left( u^0, u^1 \right) \in \mathcal{H}$ and $\phi \in \mathscr{C}_\mathrm{c}^\infty(\mathbb{R} \times \Omega)$. Note that (\ref{proof_thm_local_decay_energy_eq_1}) is equivalent to
\[\left\vert \int_\mathbb{R} \left\langle T(t) \left( u^0, u^1 \right), \phi(t) \right\rangle_{L^2(\Omega)} \mathrm{d}t \right\vert \leq C \left\Vert \phi \right\Vert_{L^2(\mathbb{R} \times \Omega)} \left\Vert \left( u^0, u^1 \right) \right\Vert_{\mathcal{H}},\]
for some $C > 0$ independent of $\left( u^0, u^1 \right)$ and $\phi$. The Cauchy-Schwarz inequality gives
\[\left\vert \int_\mathbb{R} \left\langle T(t) \left( u^0, u^1 \right), \phi(t) \right\rangle_{L^2(\Omega)} \mathrm{d}t \right\vert \leq \left\Vert \left( u^0, u^1 \right) \right\Vert_{\mathcal{H}} \left\Vert \int_\mathbb{R} T(t)^\ast \phi(t) \mathrm{d}t \right\Vert_{\mathcal{H}},\]
and 
\begin{align*}
\left\Vert \int_\mathbb{R} T(t)^\ast \phi(t) \mathrm{d}t \right\Vert_{\mathcal{H}}^2 & = \left\langle \int_\mathbb{R} T(t)^\ast \phi(t) \mathrm{d}t, \int_\mathbb{R} T(s)^\ast \phi(s) \mathrm{d}s \right\rangle_{\mathcal{H}} \\
& = \int_\mathbb{R} \left\langle \phi(t), \int_\mathbb{R} T(t) T(s)^\ast \phi(s) \mathrm{d}s \right\rangle_{L^2(\Omega)} \mathrm{d}t \\
& \leq \left\Vert \phi \right\Vert_{L^2(\mathbb{R} \times \Omega)} \left\Vert T_0(\phi) \right\Vert_{L^2(\mathbb{R} \times \Omega)},
\end{align*}
where $T_0: L^2(\mathbb{R} \times \Omega) \rightarrow L^2(\mathbb{R} \times \Omega)$ is the operator given by
\[T_0(\phi): t \longmapsto \int_\mathbb{R} T(t) T(s)^\ast \phi(s) \mathrm{d}s = \int_\mathbb{R} \chi \pi_1 S(t - s) \left( 0, \chi \phi(s) \right) \mathrm{d}s, \quad \phi \in L^2(\mathbb{R} \times \Omega).\]
Write $T_0^\pm$ for the operator 
\[T_0^\pm(\psi): t \longmapsto \int_\mathbb{R} \mathds{1}_{t - s \in \mathbb{R}_\pm} \pi_1 S(t - s) \left( 0, \psi(s) \right) \mathrm{d}s, \quad \psi \in L^2(\mathbb{R} \times \Omega),\]
so that $T_0 \phi = \chi T_0^+(\chi \phi) + \chi T_0^-(\chi \phi)$. To prove (\ref{proof_thm_local_decay_energy_eq_1}), we show that
\begin{equation}\label{proof_thm_local_decay_energy_eq_2b}
\left\Vert \chi T_0^\pm(\chi \phi) \right\Vert_{L^2(\mathbb{R} \times \Omega)} \lesssim \left\Vert \phi \right\Vert_{L^2(\mathbb{R} \times \Omega)}, \quad \phi \in L^2(\mathbb{R} \times \Omega).
\end{equation}

We start with the contribution of $\chi T_0^+(\chi \phi)$. Set 
\[U(t) = \int_\mathbb{R} \mathds{1}_{t - s > 0} S(t - s) \left( 0, \chi \phi(s) \right) \mathrm{d}s, \quad t \in \mathbb{R},\]
so that $\pi_1 U = T_0^+(\chi \phi)$. Let $R_1 > 0$ be such that $\supp \phi \subset \left(- R_1, R_1 \right) \times \Omega$. One has $U(t) \in D(A)$ for $t \in \mathbb{R}$, and
\begin{equation}\label{proof_thm_local_decay_energy_eq_3}
\left \{
\begin{array}{rcccl}
\partial_t U & = & A U + \left(0, \chi \phi \right) & \quad & \text{in } \mathbb{R} \times \Omega, \\
U & = & 0 & \quad & \text{in } (- \infty, - R_1) \times \Omega.
\end{array}
\right.
\end{equation}
By (\ref{proof_thm_local_decay_energy_eq_2}), one has $\sup_{t \in \mathbb{R}} \left\Vert U(t) \right\Vert_{\mathcal{H}} < \infty$, yielding
\[\int_{\mathbb{R}} \left\Vert U(t) e^{- i \tau t} \right\Vert_{\mathcal{H}} \mathrm{d}t \lesssim \int_{- R_1}^{+ \infty} e^{\Im(\tau) t} \mathrm{d} t < \infty, \quad \Im \tau < 0.\]
This implies that the Fourier transform of $U$ with respect to $t$, defined by
\[\widehat{U}(\tau) = \int_{\mathbb{R}} U(t) e^{- i \tau t} \mathrm{d}t, \quad \Im \tau < 0,\]
is holomorphic in the half-plane $\left\{ \Im \tau < 0 \right\}$. Using (\ref{proof_thm_local_decay_energy_eq_3}), one obtains
\[\left( i \tau - A \right) \widehat{U}(\tau) = \left( 0, \chi \widehat{\phi}(\tau) \right), \quad \Im \tau < 0.\]
If $\tau^2 \in \mathbb{C} \backslash \left[ \beta, + \infty \right)$, then the operator $i \tau - A$ is invertible, with
\[\left( i \tau - A \right)^{-1} =
\left(
\begin{array}{cc}
i \tau \left( - \Delta + \beta - \tau^2 \right)^{-1} & \left( - \Delta + \beta - \tau^2 \right)^{-1} \\ 
\left( \Delta - \beta \right) \left( - \Delta + \beta - \tau^2 \right)^{-1} & i \tau \left( - \Delta + \beta - \tau^2 \right)^{-1}
\end{array}
\right).\]
Note that for $\tau = \tau_0 + i \tau_1$ with $\tau_1 < 0$, one has $\tau^2 = \tau_0^2 - \tau_1^2 + 2 i \tau_0 \tau_1$, implying $\tau^2 \notin \left[ \beta, + \infty \right)$ as $\beta \geq 0$. In particular, one has
\[\widehat{U}(\tau) = \left( i \tau - A \right)^{-1} \left( 0, \chi \widehat{\phi}(\tau) \right), \quad \Im \tau < 0,\]
yielding
\begin{equation}\label{proof_thm_local_decay_energy_eq_4}
\left\Vert \chi \pi_1 \widehat{U}(\tau) \right\Vert_{L^2(\Omega)} = \left\Vert \chi \tau \left( - \Delta + \beta - \tau^2 \right)^{-1} \chi \widehat{\phi}(\tau) \right\Vert_{L^2(\Omega)}, \quad \Im \tau < 0.
\end{equation}
We use the following lemma.

\begin{lem}\label{lem_resolvent_estimate}
Assume that $\beta > 0$. Then there exists $C > 0$ such that 
\[\left(1 + \vert \tau \vert \right) \left\Vert \chi \left( - \Delta + \beta - \tau^2 \right)^{-1} \chi w \right\Vert_{L^2(\Omega)} \leq C \Vert w \Vert_{L^2(\Omega)}, \quad \Im \tau \neq 0, \quad w \in L^2(\Omega).\]
\end{lem}

\begin{proof}
We prove
\begin{equation}\label{proof_lem_resolvent_eq_1}
\left(1 + \left\vert \tau_0 + i \tau_1 \right\vert \right) \left\Vert \chi \left( - \Delta + \beta - \tau_0^2 + \tau_1^2 - 2 i \tau_0 \tau_1 \right)^{-1} \chi w \right\Vert_{L^2(\Omega)} \leq C \Vert w \Vert_{L^2(\Omega)},
\end{equation}
for $w \in L^2(\Omega)$, $\tau_0 \in \mathbb{R}$ and $\tau_1 \neq 0$. We start with the case $\tau_0 = 0$. For $u \in H^2(\Omega) \cap H_0^1(\Omega)$ and $\tau_1 \in \mathbb{R}$, integrating by parts, one finds
\[\left\Vert \left( - \Delta + \beta + \tau_1^2 \right) u \right\Vert_{L^2(\Omega)}^2 = \left\Vert \left( - \Delta + \beta \right) u \right\Vert_{L^2(\Omega)}^2 + \tau_1^4 \left\Vert u \right\Vert_{L^2(\Omega)}^2 + 2 \tau_1^2 \left\Vert u \right\Vert_{H_0^1(\Omega)}^2,\]
implying
\[
\left\Vert \left( - \Delta + \beta + \tau_1^2 \right) u \right\Vert_{L^2(\Omega)}^2 \gtrsim \left\Vert u \right\Vert_{L^2(\Omega)}^2 + \tau_1^2 \left\Vert u \right\Vert_{L^2(\Omega)}^2,\]
by the Poincaré inequality and the ellipticity of $- \Delta + \beta$. This gives 
\[\left(1 + \left\vert \tau_1 \right\vert \right) \left\Vert \left( - \Delta + \beta + \tau_1^2 \right)^{-1} w \right\Vert_{L^2(\Omega)} \lesssim \Vert w \Vert_{L^2(\Omega)}, \quad \tau_1 \in \mathbb{R}, \quad w \in L^2(\Omega),\]
implying (\ref{proof_lem_resolvent_eq_1}) in the case $\tau_0 = 0$. 

For $\tau = \tau_0 + i \tau_1 \in \mathbb{C}$, with $\tau_0 \neq 0$ and $\tau_1 \neq 0$, one has $\Im \left( \beta - \tau^2 \right) \neq 0$, implying
\[\sqrt{1 + \left\vert \beta - \tau^2 \right\vert} \left\Vert \chi \left( - \Delta + \beta - \tau^2 \right)^{-1} \chi w \right\Vert_{L^2(\Omega)} \leq C \Vert w \Vert_{L^2(\Omega)}, \quad w \in L^2(\Omega),\]
as $\Omega$ is non-trapping (see Definition \ref{def_resolvent_estimate_non_trapping}). As $\sqrt{1 + \left\vert \beta - \tau^2 \right\vert} \gtrsim 1 + \vert \tau \vert$ for $\tau \in \mathbb{C}$, this completes the proof.
\end{proof}

Coming back to (\ref{proof_thm_local_decay_energy_eq_4}) and using Lemma \ref{lem_resolvent_estimate}, one finds
\[\left\Vert \chi \pi_1 \widehat{U}(\tau) \right\Vert_{L^2(\Omega)} \lesssim \left\Vert \widehat{\phi}(\tau) \right\Vert_{L^2(\Omega)}, \quad \Im \tau < 0.\]
Writing $\tau = \tau_0 + i \tau_1$ and letting $\tau_1$ tends to zero, one obtains 
\[\left\Vert \chi \pi_1 \widehat{U}(\tau_0) \right\Vert_{L^2(\Omega)} \lesssim \left\Vert \widehat{\phi}(\tau_0) \right\Vert_{L^2(\Omega)}, \quad \tau_0 \in \mathbb{R}.\]
The Plancherel theorem gives
\begin{equation}\label{proof_thm_local_decay_energy_eq_5}
\left\Vert \chi T_0^+(\chi \phi) \right\Vert_{L^2(\mathbb{R} \times \Omega)}^2 = \left\Vert \chi \pi_1 U \right\Vert_{L^2(\mathbb{R} \times \Omega)}^2 \lesssim \left\Vert \phi \right\Vert_{L^2(\mathbb{R} \times \Omega)}^2.
\end{equation}

To estimate the contribution of $\chi T_0^-(\chi \phi)$, one argues similarly. Set 
\[V(t) = \int_\mathbb{R} \mathds{1}_{t - s < 0} S(t - s) \left( 0, \chi \phi(s) \right) \mathrm{d}s, \quad t \in \mathbb{R}.\]
One has $\pi_1 V = T_0^-(\chi \phi)$, and
\[\left \{
\begin{array}{rcccl}
\partial_t V & = & A V - \left(0, \chi \phi \right) & \quad & \text{in } \mathbb{R} \times \Omega, \\
V & = & 0 & \quad & \text{in } (R_1, + \infty) \times \Omega.
\end{array}
\right.\]
Arguing as above, one finds
\[\widehat{V}(\tau) = - \left( i \tau - A \right)^{-1} \left( 0, \chi \widehat{\phi}(\tau) \right), \quad \Im \tau > 0,\]
and with Lemma \ref{lem_resolvent_estimate} and the Plancherel theorem, this implies
\[\left\Vert \chi T_0^-(\chi \phi) \right\Vert_{L^2(\mathbb{R} \times \Omega)}^2 = \left\Vert \chi \pi_1 V \right\Vert_{L^2(\mathbb{R} \times \Omega)}^2 \lesssim \left\Vert \phi \right\Vert_{L^2(\mathbb{R} \times \Omega)}^2.\]
Together with (\ref{proof_thm_local_decay_energy_eq_5}), this gives (\ref{proof_thm_local_decay_energy_eq_2b}), completing the proof of (\ref{proof_thm_local_decay_energy_eq_1}).

\paragraph{Step 2: a second $TT^\ast$ argument.} 
Here, we prove that 
\begin{equation}\label{proof_thm_local_decay_energy_eq_6}
\left\Vert \chi u \right\Vert_{L^2(\mathbb{R} \times \Omega)} \leq C \left\Vert \left( u^0, u^1 \right) \right\Vert_{L^2(\Omega) \times H^{-1}(\Omega)}, \quad \left( u^0, u^1 \right) \in L^2(\Omega) \times H^{-1}(\Omega),
\end{equation}
where $u$ is the solution of (\ref{eq_thm_local_decay_energy}) with $F = 0$, at another level of regularity (see Proposition \ref{prop_LKG_time_dependent_potential}-\emph{(iii)} with $V = 0$). We define the scalar product on $H^{-1}(\Omega)$ as
\[\left\langle u^1, v^1 \right\rangle_{H^{-1}(\Omega)} = \left\langle \left( - \Delta + \beta \right)^{- \frac{1}{2}} u^1, \left( - \Delta + \beta \right)^{- \frac{1}{2}} v^1 \right\rangle_{L^2(\Omega)}, \quad u^1, v^1 \in H^{-1}(\Omega),\]
and we write $\mathcal{K} = L^2(\Omega) \times H^{-1}(\Omega)$, which is a Hilbert space for the scalar product 
\[\left\langle \left( u^0, u^1 \right), \left( v^0, v^1 \right) \right\rangle_{\mathcal{K}} = \left\langle u^0, v^0 \right\rangle_{L^2(\Omega)} + \left\langle u^1, v^1 \right\rangle_{H^{-1}(\Omega)}.\]
The semi-group $\mathsf{S}(t): \mathcal{K} \rightarrow \mathcal{K}$ associated with (\ref{eq_thm_local_decay_energy}) at the level of regularity $L^2(\Omega) \times H^{-1}(\Omega)$ is generated by 
\[\mathsf{A} = \left(
\begin{array}{cc}
0 & \mathrm{Id} \\ 
\Delta - \beta & 0
\end{array}
\right): D(\mathsf{A}) \subset \mathcal{K} \rightarrow \mathcal{K}, \quad D(\mathsf{A}) = H_0^1(\Omega) \times L^2(\Omega).\]
One has 
\[\left\langle \mathsf{A} \left( u^0, u^1 \right), \left( v^0, v^1 \right) \right\rangle_{\mathcal{K}} = \left\langle u^1, v^0 \right\rangle_{L^2(\Omega)} - \left\langle u^0, v^1 \right\rangle_{L^2(\Omega)} = - \left\langle \left( u^0, u^1 \right), \mathsf{A} \left( v^0, v^1 \right) \right\rangle_{\mathcal{K}},\]
for $\left( u^0, u^1 \right), \left( v^0, v^1 \right) \in D(\mathsf{A})$, that is, $\mathsf{A}^\ast = - \mathsf{A}$. Note that $\mathsf{S}(t) = \left( - \Delta + \beta \right)^{\frac{1}{2}} S(t) \left( - \Delta + \beta \right)^{- \frac{1}{2}}$. Hence, (\ref{proof_thm_local_decay_energy_eq_2}) implies
\[\left\Vert \mathsf{S}(t) \left( u^0, u^1 \right) \right\Vert_{\mathcal{K}} \lesssim \left\Vert \left( u^0, u^1 \right) \right\Vert_{\mathcal{K}}, \quad t \in \mathbb{R}, \quad \left( u^0, u^1 \right) \in \mathcal{K}.\]

Denote by $\pi_0: \mathbb{R}^2 \rightarrow \mathbb{R}$ the projection on the first coordinate, and for $t \in \mathbb{R}$, consider
\[\begin{array}{cccc} 
\mathsf{T}(t): & \mathcal{K} & \longrightarrow & L^2(\Omega) \\
& \left( u^0, u^1 \right) & \longmapsto & \chi u(t) = \chi \pi_0 \mathsf{S}(t) \left( u^0, u^1 \right)
\end{array}.\]
One has $\mathsf{T}(t)^\ast v^0 = \mathsf{S}(- t) \left( \chi v^0, 0 \right)$, for $t \in \mathbb{R}$ and $v^0 \in L^2(\Omega)$. As above, (\ref{proof_thm_local_decay_energy_eq_6}) will follow from
\[\left\Vert \mathsf{T}_0(\phi) \right\Vert_{L^2(\mathbb{R} \times \Omega)} \lesssim \left\Vert \phi \right\Vert_{L^2(\mathbb{R} \times \Omega)}, \quad \phi \in L^2(\mathbb{R} \times \Omega),\]
where $\mathsf{T}_0$ is given by
\[\mathsf{T}_0(\phi): t \longmapsto \int_\mathbb{R} \chi \pi_0 \mathsf{S}(t - s) \left( \chi \phi(s), 0 \right) \mathrm{d}s, \quad \phi \in L^2(\mathbb{R} \times \Omega).\]
One has $\mathsf{T}_0 \phi = \chi \mathsf{T}_0^+(\chi \phi) + \chi \mathsf{T}_0^-(\chi \phi)$, with
\[\mathsf{T}_0^\pm(\psi): t \longmapsto \int_\mathbb{R} \mathds{1}_{t - s \in \mathbb{R}_\pm} \pi_0 \mathsf{S}(t - s) \left( \psi(s), 0 \right) \mathrm{d}s, \quad \psi \in L^2(\mathbb{R} \times \Omega).\]

We only estimate the contribution of $\chi \mathsf{T}_0^+(\chi \phi)$, the corresponding estimate for $\chi \mathsf{T}_0^-(\chi \phi)$ being similar. As above, set 
\[\mathsf{U}(t) = \int_\mathbb{R} \mathds{1}_{t - s > 0} \mathsf{S}(t - s) \left( \chi \phi(s), 0 \right) \mathrm{d}s, \quad t \in \mathbb{R},\]
so that $\pi_0 \mathsf{U} = \mathsf{T}_0^+(\chi \phi)$, and let $R_1 > 0$ be such that $\supp \phi \subset \left(- R_1, R_1 \right) \times \Omega$. One has
\[\left \{
\begin{array}{rcccl}
\partial_t \mathsf{U} & = & \mathsf{A} \mathsf{U} + \left( \chi \phi, 0 \right) & \quad & \text{in } \mathbb{R} \times \Omega, \\
\mathsf{U} & = & 0 & \quad & \text{in } (- \infty, - R_1) \times \Omega.
\end{array}
\right.\]
One has $\widehat{\mathsf{U}}(\tau) = \left( i \tau - \mathsf{A} \right)^{-1} \left( \chi \widehat{\phi}(\tau), 0 \right)$ for $\Im \tau < 0$, implying
\[\left\Vert \chi \pi_0 \widehat{\mathsf{U}}(\tau) \right\Vert_{L^2(\Omega)} = \left\Vert \chi \tau \left( - \Delta + \beta - \tau^2 \right)^{-1} \chi \widehat{\phi}(\tau) \right\Vert_{L^2(\Omega)}, \quad \Im \tau < 0.\]
As above, Lemma \ref{lem_resolvent_estimate} gives
\[\left\Vert \chi \pi_0 \widehat{\mathsf{U}}(\tau) \right\Vert_{L^2(\Omega)} \lesssim \left\Vert \widehat{\phi}(\tau) \right\Vert_{L^2(\Omega)}, \quad \Im \tau < 0.\]
Letting $\Im \tau$ tends to zero, and using the Plancherel theorem as above, one obtains
\[\left\Vert \chi \mathsf{T}_0^+(\chi \phi) \right\Vert_{L^2(\mathbb{R} \times \Omega)}^2 = \left\Vert \chi \pi_0 \mathsf{U} \right\Vert_{L^2(\mathbb{R} \times \Omega)}^2 \lesssim \left\Vert \phi \right\Vert_{L^2(\mathbb{R} \times \Omega)}^2.\]
This proves (\ref{proof_thm_local_decay_energy_eq_6}).

\paragraph{Step 3: interpolation.} 
Here, we prove that 
\begin{equation}\label{proof_thm_local_decay_energy_eq_8}
\left\Vert \chi u \right\Vert_{L^2(\mathbb{R}, H_0^1(\Omega))} \leq C \left\Vert \left( u^0, u^1 \right) \right\Vert_{H_0^1(\Omega) \times L^2(\Omega)}, \quad \left( u^0, u^1 \right) \in H_0^1(\Omega) \times L^2(\Omega),
\end{equation}
where $u$ is the solution of (\ref{eq_thm_local_decay_energy}) with $F = 0$. By interpolation, (\ref{proof_thm_local_decay_energy_eq_8}) follows from Step 2 and 
\begin{equation}\label{proof_thm_local_decay_energy_eq_9}
\left\Vert \chi u \right\Vert_{L^2(\mathbb{R}, H^2 \cap H_0^1)} \leq C \left\Vert \left( u^0, u^1 \right) \right\Vert_{H^2 \cap H_0^1(\Omega) \times H_0^1(\Omega)}, \quad \left( u^0, u^1 \right) \in \left(H^2(\Omega) \cap H_0^1(\Omega) \right) \times H_0^1(\Omega),
\end{equation}
where the norm of $H^2(\Omega) \cap H_0^1(\Omega)$ is given by 
\[\left\Vert u^0 \right\Vert_{H^2 \cap H_0^1 (\Omega)} = \left\Vert \left( - \Delta + \beta \right) u^0 \right\Vert_{L^2(\Omega)}, \quad u^0 \in H^2(\Omega) \cap H_0^1(\Omega),\]
and where $u$ is the solution of (\ref{eq_thm_local_decay_energy}) with $F = 0$.

We use the following elementary lemma.

\begin{lem}\label{lem_basic_ipp_estimate}
There exist $\tilde{\chi} \in \mathscr{C}^\infty_\mathrm{c}(\Omega)$ and $C > 0$ such that
\[\left\Vert \chi u \right\Vert_{H^2 \cap H_0^1 (\Omega)} \leq C \left( \left\Vert \tilde{\chi} \left( - \Delta + \beta \right) u \right\Vert_{L^2(\Omega)} + \left\Vert \tilde{\chi} u \right\Vert_{L^2(\Omega)} \right), \quad u \in H^2(\Omega) \cap H_0^1(\Omega).\]
\end{lem}

\begin{proof}
Fix $u \in H^2(\Omega) \cap H_0^1(\Omega)$, and let $\chi_1 \in \mathscr{C}^\infty_\mathrm{c}(\Omega)$ be such that $\chi_1 \chi = \chi$. One has
\begin{align*}
\left\Vert \chi u \right\Vert_{H^2 \cap H_0^1 (\Omega)} & \leq \left\Vert \chi \left( - \Delta + \beta \right) u \right\Vert_{L^2(\Omega)} + \left\Vert \left( \Delta \chi \right) u \right\Vert_{L^2(\Omega)} + 2 \left\Vert \nabla \chi \cdot \nabla u \right\Vert_{L^2(\Omega)} \\
& \lesssim \left\Vert \chi \left( - \Delta + \beta \right) u \right\Vert_{L^2(\Omega)} + \left\Vert \left( \Delta \chi \right) u \right\Vert_{L^2(\Omega)} + \left\Vert \chi_1 u \right\Vert_{H_0^1(\Omega)}.
\end{align*}
We only need to estimate the last term. Integrating by parts, one finds
\begin{align}
\left\Vert \chi_1 u \right\Vert_{H_0^1(\Omega)}^2 & = \left\langle \left( - \Delta + \beta \right) \left( \chi_1 u \right), \chi_1 u \right\rangle_{L^2(\Omega)} \nonumber \\
& = \left\langle \chi_1 \left( - \Delta + \beta \right) u, \chi_1 u \right\rangle_{L^2(\Omega)} - 2 \left\langle \nabla \chi_1 \cdot \nabla u, \chi_1 u \right\rangle_{L^2(\Omega)} - \left\langle \left( \Delta \chi_1 \right) u, \chi_1 u \right\rangle_{L^2(\Omega)}. \label{proof_lem_basic_ipp_estimate_eq_1}
\end{align}
Another integration by parts gives
\begin{equation}\label{proof_lem_basic_ipp_estimate_eq_2}
- 2 \left\langle \nabla \chi_1 \cdot \nabla u, \chi_1 u \right\rangle_{L^2(\Omega)} = \frac{1}{2} \int_\Omega \Delta \left( \chi_1^2 \right) u^2 \mathrm{d}x,
\end{equation}
and together with (\ref{proof_lem_basic_ipp_estimate_eq_1}) and the Cauchy-Schwarz inequality, this yields
\[\left\Vert \chi_1 u \right\Vert_{H_0^1(\Omega)}^2 \lesssim \left\Vert \chi_2 \left( - \Delta + \beta \right) u \right\Vert_{L^2(\Omega)}^2 + \left\Vert \chi_2 u \right\Vert_{L^2(\Omega)},\]
for some $\chi_2 \in \mathscr{C}^\infty_\mathrm{c}(\Omega)$ satisfying $\chi_2 \chi_1 = \chi_1$. This completes the proof.
\end{proof}

Now, we prove (\ref{proof_thm_local_decay_energy_eq_9}). Consider $\left( u^0, u^1 \right) \in \left(H^2(\Omega) \cap H_0^1(\Omega) \right) \times H_0^1(\Omega)$, and write $u$ for the solution of (\ref{eq_thm_local_decay_energy}) of initial data $\left( u^0, u^1 \right)$, with $F = 0$. Lemma \ref{lem_basic_ipp_estimate} gives
\[\left\Vert \chi u \right\Vert_{L^2((-T, T), H^2 \cap H_0^1)} \lesssim \left\Vert \tilde{\chi} \left( - \Delta + \beta \right) u \right\Vert_{L^2((-T, T) \times \Omega)} + \left\Vert \tilde{\chi} u \right\Vert_{L^2((-T, T) \times \Omega)}, \quad T > 0,\]
for some $\tilde{\chi} \in \mathscr{C}^\infty_\mathrm{c}(\Omega)$. Note that $v = \left( - \Delta + \beta \right) u$ is the solution of (\ref{eq_thm_local_decay_energy}) of initial data 
\[\left( v^0, v^1 \right) = \left( \left( - \Delta + \beta \right) u^0, \left( - \Delta + \beta \right) u^1 \right) \in L^2(\Omega) \times H^{-1}(\Omega),\] 
and with $F = 0$. Applying the estimate of Step 2 (with $\tilde{\chi}$ instead of $\chi$) to $u$ and $v$, one finds
\begin{align*}
\left\Vert \chi u \right\Vert_{L^2((-T, T), H^2 \cap H_0^1)} & \lesssim \left\Vert \left( v^0, v^1 \right) \right\Vert_{L^2(\Omega) \times H^{-1}(\Omega)} + \left\Vert \left( u^0, u^1 \right) \right\Vert_{L^2(\Omega) \times H^{-1}(\Omega)} \\
& \lesssim \left\Vert \left( u^0, u^1 \right) \right\Vert_{H^2 \cap H_0^1(\Omega) \times H_0^1(\Omega)},
\end{align*}
for $T > 0$. This gives (\ref{proof_thm_local_decay_energy_eq_9}). As explained above, this proves that (\ref{proof_thm_local_decay_energy_eq_8}) holds true.

\paragraph{Step 4: the inhomogeneous estimate.} 
Here, we prove
\begin{equation}\label{proof_thm_local_decay_energy_eq_10}
\left\Vert \left( \chi u, \chi \partial_t u \right) \right\Vert_{L^2(\mathbb{R}, H_0^1 \times L^2)} \lesssim \Vert F \Vert_{L^2(\mathbb{R} \times \Omega)}, \quad F \in L^2(\mathbb{R} \times B(0, R_0)),
\end{equation}
where $u$ is the solution of (\ref{eq_thm_local_decay_energy}) of initial data $\left( u^0, u^1 \right) = 0$, with source term $F$. By linearity, together with Step 1 and Step 3, this will complete the proof of Theorem \ref{thm_local_decay_energy}. Writing $F = F \mathds{1}_{[0, + \infty)} + F \mathds{1}_{(- \infty, 0)}$ and using the linearity and the time-reversibility of (\ref{eq_thm_local_decay_energy}), it suffices to prove (\ref{proof_thm_local_decay_energy_eq_10}) for $F$ supported in $\mathbb{R}_+$. By density, one can also assume that $F$ is smooth and compactly supported.

Consider $F \in \mathscr{C}^\infty_\mathrm{c}(\mathbb{R}_+ \times B(0, R_0))$, and write $U = \left(u, \partial_t u \right)$, where $u$ is the solution of (\ref{eq_thm_local_decay_energy}) of initial data $\left( u^0, u^1 \right) = 0$ and with source term $F$. By the Duhamel formula, one has
\[U(t) = \int_0^t S(t - s) \left(0, F(s) \right) \mathrm{d}s, \quad t \in \mathbb{R},\]
yielding
\[\left\Vert U(t) \right\Vert_{H_0^1(\Omega) \times L^2(\Omega)} \lesssim \int_0^t \left\Vert F(s) \right\Vert_{L^2(\Omega)} \mathrm{d}s \lesssim 1, \quad t \in \mathbb{R},\]
by (\ref{proof_thm_local_decay_energy_eq_2}). As $U(t) = 0$ for $t \leq 0$, this implies that the Fourier transform of $U$ is holomorphic in the half-plane $\left\{ \Im \tau < 0 \right\}$. As $u$ is a solution of (\ref{eq_thm_local_decay_energy}), one finds
\[\left( i \tau - A \right) \widehat{U}(\tau) = \left( 0, \widehat{F}(\tau) \right), \quad \Im \tau < 0.\]
As in Step 1, one has
\[\widehat{U}(\tau) = \left( i \tau - A \right)^{-1} \left( 0, \widehat{F}(\tau) \right) = \left( \left( - \Delta + \beta - \tau^2 \right)^{-1} \widehat{F}(\tau), i \tau \left( - \Delta + \beta - \tau^2 \right)^{-1} \widehat{F}(\tau) \right), \quad \Im \tau < 0.\]
Let $\chi_0 \in \mathscr{C}^\infty_\mathrm{c}(\Omega)$ be such that $\chi_0 \chi = \chi$ and $\chi_0 F = F$. One has
\begin{align}
\left\Vert \chi \widehat{U}(\tau) \right\Vert_{H_0^1(\Omega) \times L^2(\Omega)}^2 & \lesssim \left\Vert \chi_0 \left( - \Delta + \beta - \tau^2 \right)^{-1} \chi_0 \widehat{F}(\tau) \right\Vert_{H_0^1(\Omega)}^2 \nonumber \\
& \ + \left\Vert \chi_0 \tau \left( - \Delta + \beta - \tau^2 \right)^{-1} \chi_0 \widehat{F}(\tau) \right\Vert_{L^2(\Omega)}^2, \quad \Im \tau < 0. \label{proof_thm_local_decay_energy_eq_11}
\end{align}

Fix $\tau \in \mathbb{C}$, $\Im \tau < 0$. On the one hand, Lemma \ref{lem_resolvent_estimate} gives
\begin{equation}\label{proof_thm_local_decay_energy_eq_12}
\left\Vert \chi_0 \tau \left( - \Delta + \beta - \tau^2 \right)^{-1} \chi_0 \widehat{F}(\tau) \right\Vert_{L^2(\Omega)} \lesssim \left\Vert \widehat{F}(\tau) \right\Vert_{L^2(\Omega)}.
\end{equation}
On the other hand, to estimate the other term of (\ref{proof_thm_local_decay_energy_eq_11}), we use (\ref{proof_lem_basic_ipp_estimate_eq_1}) and (\ref{proof_lem_basic_ipp_estimate_eq_2}). It gives
\[\left\Vert \chi_0 w \right\Vert_{H_0^1(\Omega)}^2 = \left\langle \chi_0 \left( - \Delta + \beta \right) w, \chi_0 w \right\rangle_{L^2(\Omega)} + \frac{1}{2} \int_\Omega \Delta \left( \chi_0^2 \right) w^2 \mathrm{d}x - \left\langle \left( \Delta \chi_0 \right) w, \chi_0 w \right\rangle_{L^2(\Omega)}.\]
with $w = \left( - \Delta + \beta - \tau^2 \right)^{-1} \chi_0 \widehat{F}(\tau)$. One has 
\begin{align*}
\left\langle \chi_0 \left( - \Delta + \beta \right) w, \chi_0 w \right\rangle_{L^2(\Omega)} & = \left\langle \chi_0 \left( - \Delta + \beta - \tau^2 \right) w, \chi_0 w \right\rangle_{L^2(\Omega)} + \tau^2 \left\langle \chi_0 w, \chi_0 w \right\rangle_{L^2(\Omega)} \\
& = \left\langle \chi_0^2 \widehat{F}(\tau), \chi_0 w \right\rangle_{L^2(\Omega)} + \tau^2 \left\Vert \chi_0 w \right\Vert_{L^2(\Omega)}^2
\end{align*}
implying
\[\left\Vert \chi_0 w \right\Vert_{H_0^1(\Omega)}^2 \lesssim \left\Vert \widehat{F}(\tau) \right\Vert_{L^2(\Omega)}^2 + \left( 1 + \vert \tau \vert^2 \right) \left\Vert \chi_1 \left( - \Delta + \beta - \tau^2 \right)^{-1} \chi_1 \chi_0 \widehat{F}(\tau) \right\Vert_{L^2(\Omega)}^2,\]
for some $\chi_1 \in \mathscr{C}^\infty_\mathrm{c}(\Omega)$. By Lemma \ref{lem_resolvent_estimate}, this gives
\begin{equation}\label{proof_thm_local_decay_energy_eq_13}
\left\Vert \chi_0 w \right\Vert_{H_0^1(\Omega)}^2 \lesssim \left\Vert \widehat{F}(\tau) \right\Vert_{L^2(\Omega)}^2.
\end{equation}
Using (\ref{proof_thm_local_decay_energy_eq_11}), (\ref{proof_thm_local_decay_energy_eq_12}), and (\ref{proof_thm_local_decay_energy_eq_13}), one obtains 
\[\left\Vert \chi \widehat{U}(\tau) \right\Vert_{H_0^1(\Omega) \times L^2(\Omega)}^2 \lesssim \left\Vert \widehat{F}(\tau) \right\Vert_{L^2(\Omega)}, \quad \Im \tau < 0.\]
Using the Plancherel theorem as in Step 1, one finds (\ref{proof_thm_local_decay_energy_eq_10}). This completes the proof of Theorem \ref{thm_local_decay_energy}.
\end{proof} 

\subsection{Global Strichartz estimates for a non-trapping exterior domain}

Here, we prove Theorem \ref{thm_global_strichartz}. As in \cite{Burq03}, Theorem \ref{thm_global_strichartz} will be a consequence of the local energy decay (Theorem \ref{thm_local_decay_energy}) and the global-in-time Strichartz estimate in the case $\Omega = \mathbb{R}^d$. The latter is derived from the following result of \cite{GinibreVelo}. A definition of the Besov spaces can be found, for example, in \cite{Adams-Fournier}, paragraph 7.32. 

\begin{thm}[Proposition 2.2 of \cite{GinibreVelo}]
Consider $d \geq 3$, $2 \leq r \leq \infty$, $\rho \in \mathbb{R}$, $1 \leq m \leq \infty$, and write
\[\delta(r) = \frac{d}{2} - \frac{d}{r}, \quad \gamma(r) = \frac{d - 1}{2} - \frac{d - 1}{r}, \quad \sigma = \rho + \delta(r) - 1, \quad \text{and} \quad \frac{1}{q} = \max(\sigma, 0).\]
Assume
\begin{equation}\label{eq_thm_prop_ginibre_1}
\sigma < \frac{1}{2}, \quad 2 \sigma \leq \gamma(r)
\end{equation}
and
\begin{equation}\label{eq_thm_prop_ginibre_2}
\left \{
\begin{array}{lcl}
\frac{1}{m} = \min \left( \frac{1}{2}, \ \frac{\delta(r)}{2}, \ \gamma(r) - \sigma \right) & \text{ if } & \min \left( \frac{\delta(r)}{2}, \ \gamma(r) - \sigma \right) \neq \frac{1}{2}, \\
\frac{1}{m} < \frac{1}{2} & \text{ if } & \min \left( \frac{\delta(r)}{2}, \ \gamma(r) - \sigma \right) = \frac{1}{2}.
\end{array}
\right.
\end{equation}
Then there exists a constant $C > 0$ such that for all $\left( u^0, u^1 \right)$ in $H^1(\mathbb{R}^d) \times L^2(\mathbb{R}^d)$, the solution $u$ of 
\[\left \{
\begin{array}{rcccl}
\square u + \beta u & = & 0 & \quad & \text{in } \mathbb{R} \times \mathbb{R}^d, \\
(u(0), \partial_t u(0)) & = & \left( u^0, u^1 \right) & \quad & \text{in } \mathbb{R}^d, \\
\end{array}
\right.\]
satisfies
\begin{equation}\label{eq_thm_prop_ginibre_3}
\left( \sum_{z \in \mathbb{Z}} \left( \int_{z - \frac{1}{2}}^{z + \frac{1}{2}} \Vert u(t) \Vert_{B_{r, 2}^\rho}^q \mathrm{d} t \right)^{\frac{m}{q}} \right)^{\frac{1}{m}} \leq C \left\Vert \left(u^0, u^1 \right) \right\Vert_{H^1(\mathbb{R}^d) \times L^2(\mathbb{R}^d)}.
\end{equation}
\end{thm}

We prove the following corollary. 

\begin{cor}\label{cor_global_strichartz_whole_space}
Consider $3 \leq d \leq 5$. For $2 < \alpha \leq \frac{d + 2}{d - 2}$, there exists $C > 0$ such that for all $\left( u^0, u^1 \right)$ in $H^1(\mathbb{R}^d) \times L^2(\mathbb{R}^d)$, and all $F \in L^1(\mathbb{R}, L^2(\mathbb{R}^d))$, the solution $u$ of 
\begin{equation}\label{eq_cor_global_strichart_whole_space}
\left \{
\begin{array}{rcccl}
\square u + \beta u & = & F & \quad & \text{in } \mathbb{R} \times \mathbb{R}^d, \\
(u(0), \partial_t u(0)) & = & \left( u^0, u^1 \right) & \quad & \text{in } \mathbb{R}^d, \\
\end{array}
\right.
\end{equation}
satisfies
\[\Vert u \Vert_{L^\alpha(\mathbb{R}, L^{2 \alpha}(\mathbb{R}^d))} \leq C \left( \left\Vert \left(u^0, u^1 \right) \right\Vert_{H^1(\mathbb{R}^d) \times L^2(\mathbb{R}^d)} + \Vert F \Vert_{L^1(\mathbb{R}, L^2(\mathbb{R}^d))} \right).\]
\end{cor}

\begin{proof}
We start with the case $F = 0$. The following choices are motivated by the case $d = \alpha = 3$, which can be found in \cite{NakanishiSchlagBOOK} (see in particular (2.115) and (2.121)). We choose $m = q = \alpha$, so that the left-hand side of (\ref{eq_thm_prop_ginibre_3}) is $\Vert u \Vert_{L^\alpha(\mathbb{R}, B_{r, 2}^\rho)}$. Then, we choose $\gamma(r) - \sigma = \frac{1}{m} = \frac{1}{\alpha}$, and together with $\frac{1}{\alpha} = \frac{1}{q} = \sigma$, this gives
\[\frac{1}{r} = \frac{1}{2} - \frac{2}{\alpha ( d - 1 )} \quad \text{and} \quad \rho = \frac{1}{\alpha} - \frac{d}{2} + 1 + \frac{d}{r}.\]
Using $\alpha > 2$, one can verify that (\ref{eq_thm_prop_ginibre_1}) and (\ref{eq_thm_prop_ginibre_2}) are satisfied. 

Next, we prove that $B_{r, 2}^\rho \hookrightarrow L^{2 \alpha}(\mathbb{R}^d)$. One has $B_{r_0, 2}^{\rho_0} \hookrightarrow L^p(\mathbb{R}^d)$ for $p = \frac{r_0 d}{d - \rho_0 r_0}$ (see for example (2.121) of \cite{NakanishiSchlagBOOK}). In particular, one has $B_{r, 2}^{\rho_0} \hookrightarrow L^{2 \alpha}(\mathbb{R}^d)$ for $\rho_0 = \frac{d}{r} - \frac{d}{2 \alpha}$. As $\alpha \leq \frac{d + 2}{d - 2}$, one has $\rho \geq \rho_0$, implying
\[B_{r_0, 2}^\rho \hookrightarrow B_{r_0, 2}^{\rho_0} \hookrightarrow L^{2 \alpha}(\mathbb{R}^d).\]
This completes the proof of Corollary \ref{cor_global_strichartz_whole_space} in the case $F = 0$.

Lastly, by linearity, it suffices to prove Corollary \ref{cor_global_strichartz_whole_space} in the case $\left( u^0, u^1 \right) = 0$ and $F \neq 0$. Writing $F = F \mathds{1}_{[0, + \infty)} + F \mathds{1}_{(- \infty, 0)}$ and using the linearity and the time-reversibility of (\ref{eq_cor_global_strichart_whole_space}), we can assume that $F$ is supported in $\mathbb{R}_+$. Consider $F \in L^1(\mathbb{R}_+, L^2(\mathbb{R}^d))$, and let $u$ be the solution of (\ref{eq_cor_global_strichart_whole_space}) associated with $F$ and with $\left( u^0, u^1 \right) = 0$. The Duhamel formula gives
\[\left( u(t), \partial_t u(t) \right) = \int_0^t S(t - s) \left(0, F(s) \right) \mathrm{d} s, \quad t \in \mathbb{R},\]
where $S$ is the semi-group associated with (\ref{eq_cor_global_strichart_whole_space}). Set
\[\begin{array}{cccc} 
T: & L^1(\mathbb{R}, L^2(\mathbb{R}^d)) & \longrightarrow & L^\alpha(\mathbb{R}, L^{2 \alpha}(\mathbb{R}^d)) \\
& F & \longmapsto & \left( t \mapsto \int_0^\infty \pi_0 S(t - s) \left(0, F(s) \right) \mathrm{d} s \right)
\end{array},\]
where $\pi_0: \mathbb{R}^2 \rightarrow \mathbb{R}$ is the projection on the first coordinate, and denote by $\tilde{T}: L^1(\mathbb{R}, L^2(\mathbb{R}^d)) \rightarrow L^\alpha(\mathbb{R}, L^{2 \alpha}(\mathbb{R}^d))$ the operator defined in the Christ-Kiselev lemma (Lemma \ref{lem_christ_kiselev}). One has $u = \tilde{T} F$, implying that Corollary \ref{cor_global_strichartz_whole_space} follows from Lemma \ref{lem_christ_kiselev} and the continuity of $T$. 

Consider $F \in L^1(\mathbb{R}, L^2(\mathbb{R}^d))$, and write $u = TF$. Then $u$ is the solution of 
\[\left \{
\begin{array}{rcccl}
\square u + \beta u & = & 0 & \quad & \text{in } \mathbb{R} \times \mathbb{R}^d, \\
(u(0), \partial_t u(0)) & = & \left( u^0, u^1 \right) & \quad & \text{in } \mathbb{R}^d, \\
\end{array}
\right.\]
with
\[\left( u^0, u^1 \right) = \int_0^\infty S(- s) \left( 0, F(s) \right) \mathrm{d} s.\]
Using (\ref{proof_thm_local_decay_energy_eq_2}), one finds
\[\left\Vert \left( u^0, u^1 \right) \right\Vert_{H^1(\mathbb{R}^d) \times L^2(\mathbb{R}^d)} \lesssim \Vert F \Vert_{L^1(\mathbb{R}, L^2(\mathbb{R}^d))}.\]
Hence, Corollary \ref{cor_global_strichartz_whole_space} in the case $ F = 0$ gives
\[\Vert u \Vert_{L^\alpha(\mathbb{R}, L^{2 \alpha}(\mathbb{R}^d))} \lesssim \Vert F \Vert_{L^1(\mathbb{R}, L^2(\mathbb{R}^d))}.\]
This proves that $T: L^1(\mathbb{R}, L^2(\mathbb{R}^d)) \rightarrow L^\alpha(\mathbb{R}, L^{2 \alpha}(\mathbb{R}^d))$ is well-defined and continuous, completing the proof of Corollary \ref{cor_global_strichartz_whole_space}.
\end{proof}

Now, following the strategy of \cite{Burq03} (and \cite{SmithSogge}), we use Theorem \ref{thm_local_decay_energy} and Corollary \ref{cor_global_strichartz_whole_space} to prove Theorem \ref{thm_global_strichartz}.

\begin{proof}[Proof of Theorem \ref{thm_global_strichartz}]
We split the proof in 2 steps.

\paragraph{Step 1: the homogeneous estimate.} 
Consider $\left( u^0, u^1 \right) \in H_0^1(\Omega) \times L^2(\Omega)$, and write $u$ for the solution $u$ of (\ref{eq_thm_global_strichartz_1}) with initial data $\left( u^0, u^1 \right)$ and with $F = 0$. Consider $\chi \in \mathscr{C}^\infty_\mathrm{c}(\Omega)$ such that $\chi = 1$ on $B(0, R)$, with $R > 0$ such that $\mathbb{R}^d \backslash B(0, R) \subset \Omega$, and such that the metric of $\left( \mathbb{R}^d \backslash B(0, R) \right) \cap \Omega$ is the Euclidean metric. To show that (\ref{eq_thm_global_strichartz_2}) holds true, we estimate separately the contribution of $v = \chi u$ and of $w = (1 - \chi) u$.

\underline{Contribution of $w$.} One has $w = 0$ in $B(0, R)$, and $w$ is the solution of
\[\left \{
\begin{array}{rcccl}
\square w + \beta w & = & 2 \nabla \chi \cdot \nabla u + \Delta \chi u & \quad & \text{in } \mathbb{R} \times \mathbb{R}^d, \\
(w(0), \partial_t w(0)) & = & \left( (1 - \chi) u^0, (1 - \chi) u^1 \right) & \quad & \text{in } \mathbb{R}^d.
\end{array}
\right.\]
Write $w = w_0 + w_1$, where $w_0$ is the solution of
\[\left \{
\begin{array}{rcccl}
\square w_0 + \beta w_0 & = & 0 & \quad & \text{in } \mathbb{R} \times \mathbb{R}^d, \\
(w_0(0), \partial_t w_0(0)) & = & \left( (1 - \chi) u^0, (1 - \chi) u^1 \right) & \quad & \text{in } \mathbb{R}^d,
\end{array}
\right.\]
and $w_1$ is the solution of
\begin{equation}\label{proof_thm_global_strichartz_eq_1}
\left \{
\begin{array}{rcccl}
\square w_1 + \beta w_1 & = & 2 \nabla \chi \cdot \nabla u + \Delta \chi u & \quad & \text{in } \mathbb{R} \times \mathbb{R}^d, \\
(w_1(0), \partial_t w_1(0)) & = & 0 & \quad & \text{in } \mathbb{R}^d.
\end{array}
\right.
\end{equation}
The global-in-time Strichartz estimate in $\mathbb{R}^d$ (Corollary \ref{cor_global_strichartz_whole_space}) gives
\begin{align*}
\left\Vert w_0 \right\Vert_{L^\alpha(\mathbb{R}, L^{2 \alpha}(\Omega))} & \leq \left\Vert w_0 \right\Vert_{L^\alpha(\mathbb{R}, L^{2 \alpha}(\mathbb{R}^d))} \\
& \lesssim \left\Vert \left( (1 - \chi) u^0, (1 - \chi) u^1 \right) \right\Vert_{H^1(\mathbb{R}^d) \times L^2(\mathbb{R}^d)} \\
& \lesssim \left\Vert \left( u^0, u^1 \right) \right\Vert_{H_0^1(\Omega) \times L^2(\Omega)}.
\end{align*}

Next, we estimate the contribution of $w_1$, using twice the local energy decay (Theorem \ref{thm_local_decay_energy}). Write $F_1 = 2 \nabla \chi \cdot \nabla u + \Delta \chi u$. Note that one has
\[\left\Vert F_1 \right\Vert_{L^2(\mathbb{R} \times \mathbb{R}^d)} \lesssim \left\Vert \chi_1 u \right\Vert_{L^2(\mathbb{R}, H_0^1(\Omega))}\]
for some $\chi_1 \in \mathscr{C}^\infty_{\mathrm{c}}(\Omega)$, implying
\begin{equation}\label{proof_thm_global_strichartz_eq_2bis}
\left\Vert F_1 \right\Vert_{L^2(\mathbb{R} \times \mathbb{R}^d)} \lesssim \left\Vert \left( u^0, u^1 \right) \right\Vert_{H_0^1(\Omega) \times L^2(\Omega)}
\end{equation}
by Theorem \ref{thm_local_decay_energy}. 

We prove
\begin{equation}\label{proof_thm_global_strichartz_eq_3}
\left\Vert w_1 \right\Vert_{L^\alpha((0, +\infty), L^{2 \alpha}(\Omega))} \lesssim \left\Vert F_1 \right\Vert_{L^2(\mathbb{R} \times \mathbb{R}^d)}.
\end{equation}
By the Duhamel formula, one has
\[\left( w_1(t), \partial_t w_1(t) \right) = \int_0^t S(t - s) \left(0, F_1(s) \right) \mathrm{d} s, \quad t \in \mathbb{R},\]
where $S$ is the semi-group associated with equation (\ref{proof_thm_global_strichartz_eq_1}). Consider $\chi_2 \in \mathscr{C}^\infty_\mathrm{c}(\mathbb{R}^d)$ such that $\chi_2 \chi = \chi$, set
\[\begin{array}{cccc} 
T: & L^2(\mathbb{R}, L^2(\mathbb{R}^d)) & \longrightarrow & L^\alpha(\mathbb{R}, L^{2 \alpha}(\mathbb{R}^d)) \\
& F & \longmapsto & \left( t \mapsto \int_0^\infty \pi_0 S(t - s) \left(0, \chi_2 F(s) \right) \mathrm{d} s \right)
\end{array},\]
where $\pi_0: \mathbb{R}^2 \rightarrow \mathbb{R}$ is the projection on the first coordinate, and denote by $\tilde{T}: L^2(\mathbb{R}, L^2(\mathbb{R}^d)) \rightarrow L^\alpha(\mathbb{R}, L^{2 \alpha}(\mathbb{R}^d))$ the operator defined in the Christ-Kiselev lemma (Lemma \ref{lem_christ_kiselev}). As $\alpha > 2$, we can apply Lemma \ref{lem_christ_kiselev}: to prove that $\tilde{T}$ is well-defined and continuous, it suffices to prove that $T$ is well-defined and continuous. As $\mathds{1}_{(0, + \infty)} w_1 = \tilde{T} F_1$, this will imply (\ref{proof_thm_global_strichartz_eq_3}).

By definition, one has
\[TF(t) = \int_0^\infty \pi_0 S(t - s) \left(0, \chi_2 F(s) \right) \mathrm{d}s, \quad F \in L^2(\mathbb{R}, L^2(\mathbb{R}^d)), \quad t \in \mathbb{R}.\]
Write $T = T_1 \circ T_0$, with 
\[\begin{array}{cccc} 
T_0: & L^2(\mathbb{R}, L^2(\mathbb{R}^d)) & \longrightarrow & H^1(\mathbb{R}^d) \times L^2(\mathbb{R}^d) \\
& F & \longmapsto & \int_0^\infty S(- s) \left(0, \chi_2 F(s) \right) \mathrm{d}s
\end{array},\]
and
\[\begin{array}{cccc} 
T_1: & H^1(\mathbb{R}^d) \times L^2(\mathbb{R}^d) & \longrightarrow & L^\alpha(\mathbb{R}, L^{2 \alpha}(\mathbb{R}^d)) \\
& \left( u^0, u^1 \right) & \longmapsto & \left( t \mapsto \pi_0 S(t) \left( u^0, u^1 \right) \right)
\end{array}.\]
The operator $T_1$ is continuous by the global Strichartz estimate in the case $\Omega = \mathbb{R}^d$ (Theorem \ref{cor_global_strichartz_whole_space}). By Theorem \ref{thm_local_decay_energy} (applied with $\Omega = \mathbb{R}^d$ and $\chi_2$), the operator
\[\begin{array}{cccc} 
T_2: & H^1(\mathbb{R}^d) \times L^2(\mathbb{R}^d) & \longrightarrow & L^2(\mathbb{R}, L^2(\mathbb{R}^d)) \\
& \left( u^0, u^1 \right) & \longmapsto & \left( s \mapsto \chi_2 \pi_1 S(s) \left( u^0, u^1 \right) \right)
\end{array}\]
is well-defined and continuous, where $\pi_1: \mathbb{R}^2 \rightarrow \mathbb{R}$ is the projection on the second coordinate. We prove that $T_0 = T_2^\ast$, implying that $T_0$ is well-defined and continuous. Consider $F \in \mathscr{C}^\infty_{\mathrm{c}}(\mathbb{R} \times \mathbb{R}^d)$, $\left( u^0, u^1 \right) \in H^1(\mathbb{R}^d) \times L^2(\mathbb{R}^d)$, and write
\[\left\langle T_2 \left( u^0, u^1 \right), F \right\rangle_{L^2(\mathbb{R}, L^2(\mathbb{R}^d))} = \int_\mathbb{R} \left\langle S(s) \left( u^0, u^1 \right), \left( 0, \chi_2 F(s) \right) \right\rangle_{H^1(\mathbb{R}^d) \times L^2(\mathbb{R}^d)} \mathrm{d}s.\]
As explained in the proof of Theorem \ref{thm_local_decay_energy}, one has $S(s)^\ast = S(-s)$, yielding
\[\left\langle T_2 \left( u^0, u^1 \right), F \right\rangle_{L^2(\mathbb{R}, L^2(\mathbb{R}^d))} = \left\langle \left( u^0, u^1 \right), \int_\mathbb{R} S(-s) \left( 0, \chi_2 F(s) \right) \mathrm{d}s \right\rangle_{H^1(\mathbb{R}^d) \times L^2(\mathbb{R}^d)}.\]
This proves that $T_0 = T_2^\ast$, and completes the proof of (\ref{proof_thm_global_strichartz_eq_3}).

Using (\ref{proof_thm_global_strichartz_eq_3}), and also (\ref{proof_thm_global_strichartz_eq_3}) applied to $t \mapsto w_1(-t)$, one obtains
\[\left\Vert w_1 \right\Vert_{L^\alpha(\mathbb{R}, L^{2 \alpha}(\Omega))} \lesssim \left\Vert F_1 \right\Vert_{L^2(\mathbb{R} \times \mathbb{R}^d)}.\]
Together with (\ref{proof_thm_global_strichartz_eq_2bis}), this gives
\[\left\Vert w \right\Vert_{L^\alpha(\mathbb{R}, L^{2 \alpha}(\Omega))} \lesssim \left\Vert \left( u^0, u^1 \right) \right\Vert_{H_0^1(\Omega) \times L^2(\Omega)}.\]

\underline{Contribution of $v$.}
By definition, $v$ is the solution of
\[\left \{
\begin{array}{rcccl}
\square v + \beta v & = & - 2 \nabla \chi \cdot \nabla u - \Delta \chi u & \quad & \text{in } \mathbb{R} \times \Omega, \\
(v(0), \partial_t v(0)) & = & \left( \chi u^0, \chi u^1 \right) & \quad & \text{in } \Omega, \\
v & = & 0 & \quad & \text{on } \mathbb{R} \times \partial \Omega.
\end{array}
\right.\]
Consider $\phi \in \mathscr{C}^\infty_{\mathrm{c}}((0, 1))$ such that $\phi = 1$ on $\left[ \frac{1}{4}, \frac{3}{4} \right]$, and set $v_n(t) = \phi\left( t - \frac{n}{2} \right) v(t)$, for $t \in \mathbb{R}$ and $n \in \mathbb{Z}$. One has
\[\left \{
\begin{array}{rcccl}
\square v_n + \beta v_n & = & F_n & \quad & \text{in } \mathbb{R} \times \Omega, \\
(v_n(t), \partial_t v_n(t)) & = & 0 & \quad & \text{in } \left( \mathbb{R} \backslash \left( \frac{n}{2}, \frac{n}{2} + 1 \right) \right) \times \Omega, \\
v_n & = & 0 & \quad & \text{on } \mathbb{R} \times \partial \Omega,
\end{array}
\right.\]
with 
\[F_n = - \phi\left( \cdot - \frac{n}{2} \right) \left( 2 \nabla \chi \cdot \nabla u + \Delta \chi u \right) + 2 \phi^\prime\left( \cdot - \frac{n}{2} \right) \chi \partial_t u + \phi^{\prime \prime}\left( \cdot - \frac{n}{2} \right) \chi u.\]

Consider $N \in \mathbb{N}$. Using
\[\sum_{n \in \mathbb{Z}} \phi\left( t - \frac{n}{2} \right)^\alpha \gtrsim 1, \quad t \in \mathbb{R},\]
one finds
\begin{align*}
\left\Vert v \right\Vert_{L^\alpha\left( \left(- \frac{N}{2}, \frac{N}{2} + 1 \right), L^{2 \alpha}(\Omega)\right)}^\alpha & \lesssim \int_{- \frac{N}{2}}^{\frac{N}{2} + 1} \sum_{n \in \mathbb{Z}} \phi\left( t - \frac{n}{2} \right)^\alpha \left\Vert v(t) \right\Vert_{L^{2 \alpha}(\Omega)}^\alpha \mathrm{d}t \\
& = \sum_{\vert n \vert \leq N + 1} \left\Vert v_n \right\Vert_{L^\alpha\left( \left(\frac{n}{2}, \frac{n}{2} + 1 \right), L^{2 \alpha}(\Omega)\right)}^\alpha.
\end{align*}
Using the local-in-time Strichartz estimate given by Proposition \ref{prop_LKG_time_dependent_potential}-\emph{(ii)} (with $\mathbf{u} = 0$), this gives
\[\left\Vert v \right\Vert_{L^\alpha\left( \left(- \frac{N}{2}, \frac{N}{2} + 1 \right), L^{2 \alpha}(\Omega)\right)}^\alpha \lesssim \sum_{\vert n \vert \leq N + 1} \left\Vert F_n \right\Vert_{L^1\left( \left(\frac{n}{2}, \frac{n}{2} + 1 \right), L^2(\Omega)\right)}^\alpha,\]
implying
\begin{equation}\label{proof_thm_global_strichartz_eq_4}
\left\Vert v \right\Vert_{L^\alpha\left( \left(- \frac{N}{2}, \frac{N}{2} + 1 \right), L^{2 \alpha}(\Omega)\right)}^\alpha \lesssim \left( \sum_{\vert n \vert \leq N + 1} \left\Vert F_n \right\Vert_{L^1\left( \left(\frac{n}{2}, \frac{n}{2} + 1 \right), L^2(\Omega)\right)}^2 \right)^{\frac{\alpha}{2}},
\end{equation}
as $\alpha \geq 2$. One has
\begin{align*}
\sum_{\vert n \vert \leq N + 1} \left\Vert F_n \right\Vert_{L^1\left( \left(\frac{n}{2}, \frac{n}{2} + 1 \right), L^2(\Omega)\right)}^2 & \leq \sum_{\vert n \vert \leq N + 1} \left\Vert F_n \right\Vert_{L^2\left( \left(\frac{n}{2}, \frac{n}{2} + 1 \right), L^2(\Omega)\right)}^2 \\
& \lesssim \sum_{\vert n \vert \leq N + 1} \left\Vert \left( \chi_2 u, \chi_2 \partial_t u \right) \right\Vert_{L^2\left(\left(\frac{n}{2}, \frac{n}{2} + 1 \right), H_0^1 \times L^2 \right)}^2 \\
& \lesssim \left\Vert \left( \chi_2 u, \chi_2 \partial_t u \right) \right\Vert_{L^2(\mathbb{R}, H_0^1 \times L^2)}^2
\end{align*}
for some $\chi_2 \in \mathscr{C}^\infty_{\mathrm{c}}(\Omega)$. Hence, Theorem \ref{thm_local_decay_energy} implies
\[\sum_{\vert n \vert \leq N + 1} \left\Vert F_n \right\Vert_{L^1\left( \left(\frac{n}{2}, \frac{n}{2} + 1 \right), L^2(\Omega)\right)}^2 \lesssim \left\Vert \left( u^0, u^1 \right) \right\Vert_{H_0^1(\Omega) \times L^2(\Omega)}^2.\]
Together with (\ref{proof_thm_global_strichartz_eq_4}), this gives
\[\left\Vert v \right\Vert_{L^\alpha(\mathbb{R}, L^{2 \alpha}(\Omega))} \lesssim \left\Vert \left( u^0, u^1 \right) \right\Vert_{H_0^1(\Omega) \times L^2(\Omega)},\]
completing the proof of (\ref{eq_thm_global_strichartz_2}) in the case $F = 0$.

\paragraph{Step 2: the inhomogeneous estimate.} 
Here, we prove (\ref{eq_thm_global_strichartz_2}) in the case $F \neq 0$ and $\left( u^0, u^1 \right) = 0$. The proof is similar to that of the inhomogeneous estimate of Corollary \ref{cor_global_strichartz_whole_space}, so we only sketch it. Using the Duhamel formula and the Christ-Kiselev lemma (Lemma \ref{lem_christ_kiselev}), it suffices to prove that the operator
\[\begin{array}{cccc} 
T: & L^1(\mathbb{R}, L^2(\Omega)) & \longrightarrow & L^\alpha(\mathbb{R}, L^{2 \alpha}(\Omega)) \\
& F & \longmapsto & \left(t \mapsto \int_0^\infty \pi_0 S(t - s) \left(0, F(s) \right) \mathrm{d}s \right)
\end{array}\]
is well-defined and continuous.

Consider $F \in L^1(\mathbb{R}, L^2(\Omega))$, and write $u = TF$. Then $u$ is the solution of 
\[\left \{
\begin{array}{rcccl}
\square u + \beta u & = & 0 & \quad & \text{in } \mathbb{R} \times \Omega, \\
(u(0), \partial_t u(0)) & = & \left( u^0, u^1 \right) & \quad & \text{in } \Omega, \\
u & = & 0 & \quad & \text{on } \mathbb{R} \times \partial \Omega,
\end{array}
\right.\]
with
\[\left( u^0, u^1 \right) = \int_0^\infty S(- s) \left( 0, F(s) \right) \mathrm{d} s.\]
One has
\[\left\Vert \left( u^0, u^1 \right) \right\Vert_{H_0^1(\Omega) \times L^2(\Omega)} \lesssim \Vert F \Vert_{L^1(\mathbb{R}, L^2(\Omega))},\]
implying
\[\Vert v \Vert_{L^\alpha(\mathbb{R}, L^{2 \alpha}(\Omega))} \lesssim \Vert F \Vert_{L^1(\mathbb{R}, L^2(\Omega))}\]
by Step 1. This yields (\ref{eq_thm_global_strichartz_2}) in the case $F \neq 0$. By linearity, this completes the proof of Theorem \ref{thm_global_strichartz}.
\end{proof}

\subsection{Proof of the null-controllability of a scattering solution}

Here, we prove Theorem \ref{thm_main_null_control}. We start by proving that a scattering solution is bounded in the energy space and has a finite Strichartz norm.

\begin{lem}\label{lem_scatt_implies_finite_strich_norm}
Assume that $\Omega$ is a non-trapping unbounded domain, and consider $f$ satisfying \textnormal{(\ref{eq_def_nonlinearity_scatt})} for some $\alpha_0 \leq \alpha_1$. Let $\left( u^0, u^1 \right) \in H_0^1(\Omega) \times L^2(\Omega)$ be such that the solution $u_{\NL}$ of
\[\left \{
\begin{array}{rcccl}
\square u_{\NL} + \beta u_{\NL} & = & f(u_{\NL}) & \quad & \text{in } \mathbb{R}_+ \times \Omega, \\
(u_{\NL}(0), \partial_t u_{\NL}(0)) & = & \left( u^0, u^1 \right) & \quad & \text{in } \Omega, \\
u_{\NL} & = & 0 & \quad & \text{on } \mathbb{R}_+ \times \partial \Omega,
\end{array}
\right.\]
is scattering. Then $\left( u_{\NL}, \partial_t u_{\NL} \right) \in L^\infty((0, +\infty), H_0^1(\Omega) \times L^2(\Omega))$ and 
\begin{equation}\label{eq_lem_scatt_implies_finite_strich_norm}
u_{\NL} \in L^{\alpha_0}((0, +\infty), L^{2 \alpha_0}(\Omega)) \cap L^{\alpha_1}((0, +\infty), L^{2 \alpha_1}(\Omega)).
\end{equation}
\end{lem}

\begin{proof}
First, we prove (\ref{eq_lem_scatt_implies_finite_strich_norm}). As $u_{\NL}$ is scattering, there exists a solution $u_{\Lin}$ of the linear equation
\begin{equation}\label{eq_proof_lem_scatt_implies_finite_strich_norm_0}
\left \{
\begin{array}{rcccl}
\square u_{\Lin} + \beta u_{\Lin} & = & 0 & \quad & \text{in } \mathbb{R}_+ \times \Omega, \\
u_{\Lin} & = & 0 & \quad & \text{on } \mathbb{R}_+ \times \partial \Omega.
\end{array}
\right.
\end{equation}
such that
\begin{equation}\label{eq_proof_lem_scatt_implies_finite_strich_norm_1}
\left\Vert \left( u_{\NL}(t), \partial_t u_{\NL}(t) \right) - \left( u_{\Lin}(t), \partial_t u_{\Lin}(t) \right) \right\Vert_{H_0^1(\Omega) \times L^2(\Omega)} \xrightarrow{t \rightarrow + \infty} 0.
\end{equation} 
Consider $\epsilon > 0$. By Theorem \ref{thm_global_strichartz}, one has 
\[u_{\Lin} \in L^{\alpha_0}((0, +\infty), L^{2 \alpha_0}(\Omega)) \cap L^{\alpha_1}((0, +\infty), L^{2 \alpha_1}(\Omega)).\]
Hence, using also (\ref{eq_proof_lem_scatt_implies_finite_strich_norm_1}), there exists $T = T(\epsilon)$ such that 
\begin{equation}\label{eq_proof_lem_scatt_implies_finite_strich_norm_2}
\left\Vert \left( u_{\NL}(t), \partial_t u_{\NL}(t) \right) - \left( u_{\Lin}(t), \partial_t u_{\Lin}(t) \right) \right\Vert_{H_0^1(\Omega) \times L^2(\Omega)} \leq \epsilon, \quad t \geq T,
\end{equation} 
and
\begin{equation}\label{eq_proof_lem_scatt_implies_finite_strich_norm_3}
\left\Vert u_{\Lin} \right\Vert_{L^{\alpha_0}((T, + \infty), L^{2 \alpha_0})} + \left\Vert u_{\Lin} \right\Vert_{L^{\alpha_1}((T, + \infty), L^{2 \alpha_1})} \leq \epsilon.
\end{equation}

For $T^\prime \geq T$, set
\[\eta(T^\prime) = \left\Vert u_{\NL} \right\Vert_{L^{\alpha_0}((T, T^\prime), L^{2 \alpha_0})} + \left\Vert u_{\NL} \right\Vert_{L^{\alpha_1}((T, T^\prime), L^{2 \alpha_1})}.\]
Note that $\eta(T^\prime) < + \infty$ for all $T^\prime \geq T$ by Theorem \ref{thm_existence_nL_waves}, and that $\eta$ is a continuous real function satisfying $\eta(T) = 0$. Set $v = u_{\NL} - u_{\Lin}$. Then $v$ is the solution of 
\[\left \{
\begin{array}{rcccl}
\square v + \beta v & = & f(u_{\NL}) & \quad & \text{in } \mathbb{R}_+ \times \Omega, \\
(v(T), \partial_t v(T)) & = & \left( u_{\NL}(T), \partial_t u_{\NL}(T) \right) - \left( u_{\Lin}(T), \partial_t u_{\Lin}(T) \right) & \quad & \text{in } \Omega, \\
v & = & 0 & \quad & \text{on } \mathbb{R}_+ \times \partial \Omega.
\end{array}
\right.\]
By (\ref{eq_proof_lem_scatt_implies_finite_strich_norm_3}), one has
\[\eta(T^\prime) \leq \epsilon + \left\Vert v \right\Vert_{L^{\alpha_0}((T, T^\prime), L^{2 \alpha_0})} + \left\Vert v \right\Vert_{L^{\alpha_1}((T, T^\prime), L^{2 \alpha_1})}, \quad T^\prime \geq T.\]
Using global Strichartz estimates, together with (\ref{eq_proof_lem_scatt_implies_finite_strich_norm_2}), one finds
\begin{align*}
\eta(T^\prime) & \lesssim \epsilon + \left\Vert \left( v(T), \partial_t v(T) \right) \right\Vert_{H_0^1(\Omega) \times L^2(\Omega)} + \left\Vert f(u_{\NL}) \right\Vert_{L^1((T, T^\prime), L^2)} \\
& \lesssim \epsilon + \left\Vert f(u_{\NL}) \right\Vert_{L^1((T, T^\prime), L^2)}, \\
& \lesssim \epsilon + \left\Vert f(u_{\Lin}) \right\Vert_{L^1((T, T^\prime), L^2)} + \left\Vert f(u_{\NL}) - f(u_{\Lin}) \right\Vert_{L^1((T, T^\prime), L^2)}, \quad T^\prime \geq T.
\end{align*}
Applying Lemma \ref{lem_basic_estimate_f_scatt}-\emph{(i)} twice, assuming that $\epsilon < 1$, and using (\ref{eq_proof_lem_scatt_implies_finite_strich_norm_3}) again, one obtains
\begin{align*}
\eta(T^\prime) & \lesssim \epsilon + \left\Vert u_{\Lin} \right\Vert_{L^{\alpha_0}((T, T^\prime), L^{2 \alpha_0})}^{\alpha_0} + \left\Vert u_{\Lin} \right\Vert_{L^{\alpha_1}((T, T^\prime), L^{2 \alpha_1})}^{\alpha_1} \\
& \quad + \left\Vert v \right\Vert_{L^{\alpha_0}((T, T^\prime), L^{2 \alpha_0})} \left( \left\Vert u_{\Lin} \right\Vert_{L^{\alpha_0}((T, T^\prime), L^{2 \alpha_0})}^{\alpha_0 - 1} + \left\Vert u_{\NL} \right\Vert_{L^{\alpha_0}((T, T^\prime), L^{2 \alpha_0})}^{\alpha_0 - 1} \right) \\
& \quad + \left\Vert v \right\Vert_{L^{\alpha_1}((T, T^\prime), L^{2 \alpha_1})} \left( \left\Vert u_{\Lin} \right\Vert_{L^{\alpha_1}((T, T^\prime), L^{2 \alpha_1})}^{\alpha_1 - 1} + \left\Vert u_{\NL} \right\Vert_{L^{\alpha_1}((T, T^\prime), L^{2 \alpha_1})}^{\alpha_1 - 1} \right) \\
& \lesssim \epsilon + \left( \epsilon + \eta(T^\prime) \right) \left( \epsilon + \eta(T^\prime)^{\alpha_0 - 1} + \eta(T^\prime)^{\alpha_1 - 1} \right), \\
& \lesssim \epsilon + \epsilon \eta(T^\prime) + \eta(T^\prime)^{\alpha_0} + \eta(T^\prime)^{\alpha_1}, \hspace{12em} T^\prime \geq T.
\end{align*}
Hence, for $\epsilon$ sufficiently small, one finds 
\[\eta(T^\prime) \lesssim \epsilon + \eta(T^\prime)^{\alpha_0} + \eta(T^\prime)^{\alpha_1}, \quad T^\prime \geq T.\]
By the mean value theorem, this implies that there exists $c = c(\epsilon)$ such that either $\eta(T^\prime) < c$ for all $T^\prime \geq T$, or $\eta(T^\prime) > c$ for all $T^\prime \geq T$. As $\eta(T) = 0$, this proves that $\eta$ is bounded, yielding (\ref{eq_lem_scatt_implies_finite_strich_norm}).

Second, we prove 
\begin{equation}\label{eq_proof_lem_scatt_implies_finite_strich_norm_4}
\left( u_{\NL}, \partial_t u_{\NL} \right) \in L^\infty((0, +\infty), H_0^1(\Omega) \times L^2(\Omega)).
\end{equation}
The Duhamel formula gives
\[\left( u_{\NL}(T), \partial_t u_{\NL}(T) \right) = S(T) \left( u_{\NL}(0), \partial_t u_{\NL}(0) \right) + \int_0^T S(T - t) \left( 0, f \left( u_{\NL}(t) \right) \right) \mathrm{d}t, \quad T \geq 0,\]
where $S$ is the semi-group associated with (\ref{eq_proof_lem_scatt_implies_finite_strich_norm_0}). Using (\ref{proof_thm_local_decay_energy_eq_2}), one finds
\[\left\Vert \left( u_{\NL}(T), \partial_t u_{\NL}(T) \right) \right\Vert_{H_0^1(\Omega) \times L^2(\Omega)} \lesssim 1 + \left\Vert f(u_{\NL}) \right\Vert_{L^1((0, T), L^2)}, \quad T \geq 0.\]
By Lemma \ref{lem_basic_estimate_f_scatt}-\emph{(i)}, this implies
\[\left\Vert \left( u_{\NL}(T), \partial_t u_{\NL}(T) \right) \right\Vert_{H_0^1(\Omega) \times L^2(\Omega)} \lesssim 1 + \left\Vert u_{\NL} \right\Vert_{L^{\alpha_0}((0, T)), L^{2 \alpha_0})}^{\alpha_0} + \left\Vert u_{\NL} \right\Vert_{L^{\alpha_1}((0, T), L^{2 \alpha_1})}^{\alpha_1}, \quad T \geq 0.\]
Hence, (\ref{eq_proof_lem_scatt_implies_finite_strich_norm_4}) is a consequence of (\ref{eq_lem_scatt_implies_finite_strich_norm}). This completes the proof.
\end{proof}

Now, we prove Theorem \ref{thm_main_null_control}.

\begin{proof}[Proof of Theorem \ref{thm_main_null_control}]
Consider $\left( u^0, u^1 \right) \in H_0^1(\Omega) \times L^2(\Omega)$ such that the solution $u_{\NL}$ of
\[\left\{
\begin{array}{rcccl}
\square u_{\NL} + \beta u_{\NL} & = & f(u_{\NL}) & \quad & \text{in } \mathbb{R}_+ \times \Omega, \\
(u_{\NL}(0), \partial_t u_{\NL}(0)) & = & \left( u^0, u^1 \right) & \quad & \text{in } \Omega, \\
u_{\NL} & = & 0 & \quad & \text{on } \mathbb{R}_+ \times \partial \Omega,
\end{array}
\right.\]
is scattering. Using local controllability around $0$ (Theorem \ref{thm_main_local_control_trajectory}), it suffices to show that for all $\epsilon > 0$, there exist $T$ and $g$ such that $\left\Vert \left( u(T), \partial_t u(T) \right) \right\Vert_{H_0^1(\Omega) \times L^2(\Omega)} \leq \epsilon$, where $u$ is the solution of 
\[\left \{
\begin{array}{rcccl}
\square u + \beta u & = & f(u) + g & \quad & \text{in } \mathbb{R}_+ \times \Omega, \\
(u(0), \partial_t u(0)) & = & \left( u^0, u^1 \right) & \quad & \text{in } \Omega, \\
u & = & 0 & \quad & \text{on } \mathbb{R}_+ \times \partial \Omega.
\end{array}
\right.\]

Consider $\epsilon \in (0, 1)$. As $u_{\NL}$ is scattering, there exist $T > 0$ and $u_{\Lin} \in \mathscr{C}^0(\mathbb{R}, H_0^1(\Omega)) \cap \mathscr{C}^1(\mathbb{R}, L^2(\Omega))$, satisfying $\square u_{\Lin} + \beta u_{\Lin} = 0$, such that
\begin{equation}\label{eq_proof_thm_null_control_scattering_1}
\left\Vert \left( u_{\NL}(t), \partial_t u_{\NL}(t) \right) - \left( u_{\Lin}(t), \partial_t u_{\Lin}(t) \right) \right\Vert_{H_0^1(\Omega) \times L^2(\Omega)} \leq \epsilon, \quad t \geq T.
\end{equation}
Recall that $a \geq c > 0$ on $\mathbb{R}^d \backslash B(0, R_0)$. Hence, there exists a bounded function $\chi \in \mathscr{C}^\infty(\Omega)$ such that $\chi = \frac{1}{a}$ on $\mathbb{R}^d \backslash B(0, R_0 + 1)$. Up to increasing $T$, we can assume that
\begin{equation}\label{eq_proof_thm_null_control_scattering_2}
\left\Vert (1 - a \chi) u_{\Lin} \right\Vert_{L^1((T, T + 1), L^2)} + \left\Vert (1 - a \chi) \partial_t u_{\Lin} \right\Vert_{L^1((T, T + 1), L^2)} \leq \epsilon
\end{equation}
by the local energy decay (Theorem \ref{thm_local_decay_energy}). Up to increasing $T$ again, we can assume that
\begin{equation}\label{eq_proof_thm_null_control_scattering_3}
\left\Vert u_{\NL} \right\Vert_{L^{\alpha_0}((T, T + 1), L^{2 \alpha_0})}^{\alpha_0} + \left\Vert u_{\NL} \right\Vert_{L^{\alpha_1}((T, T + 1), L^{2 \alpha_1})}^{\alpha_1} \leq \epsilon,
\end{equation}
by Lemma \ref{lem_scatt_implies_finite_strich_norm}, as $u_{\NL}$ is scattering.

Let $\varphi \in \mathscr{C}^\infty(\mathbb{R}, [0, 1])$ be such that $\varphi(t) = 1$ for $t \leq T$ and $\varphi(t) = 0$ for $t \geq T + 1$. Set $v(t, x) = \varphi(t) u_{\NL}(t, x)$ and $g = \square v + \beta v - f(v)$. By definition, $g$ is supported in $[T, T + 1]$, and one has
\[g = u_{\NL} \partial_t^2 \varphi + 2 \partial_t u_{\NL} \partial_t \varphi + \varphi f(u_{\NL}) - f(\varphi u_{\NL}).\]
To complete the proof, we prove that the solution $u$ of
\[\left \{
\begin{array}{rcccl}
\square u + \beta u & = & f(u) + a \chi g & \quad & \text{in } (0, T + 1) \times \Omega, \\
(u(0), \partial_t u(0)) & = & \left( u^0, u^1 \right) & \quad & \text{in } \Omega, \\
u & = & 0 & \quad & \text{on } (0, T + 1) \times \partial \Omega,
\end{array}
\right.\]
satisfies
\begin{equation}\label{eq_proof_thm_null_control_scattering_4}
\left\Vert \left( u(T + 1), \partial_t u(T + 1) \right) \right\Vert_{H_0^1(\Omega) \times L^2(\Omega)} \lesssim \epsilon.
\end{equation}

As $g = 0$ on $[0, T]$, one has $u = u_{\NL}$ on $[0, T]$. Note that (\ref{eq_proof_thm_null_control_scattering_4}) is equivalent to 
\begin{equation}\label{eq_proof_thm_null_control_scattering_5}
\left\Vert \left( h(T + 1), \partial_t h(T + 1) \right) \right\Vert_{H_0^1(\Omega) \times L^2(\Omega)} \lesssim \epsilon,
\end{equation}
where $h = u - v$ is the solution of
\[\left \{
\begin{array}{rcccl}
\square h + \beta h & = & f(v + h) - f(v) + (a \chi - 1) g & \quad & \text{in } (T, T + 1) \times \Omega, \\
(h(T), \partial_t h(T)) & = & 0 & \quad & \text{in } \Omega, \\
h & = & 0 & \quad & \text{on } (T, T + 1) \times \partial \Omega.
\end{array}
\right.\]

Now, we prove 
\begin{equation}\label{eq_proof_thm_null_control_scattering_6}
\left\Vert (1 - a \chi) g \right\Vert_{L^1((T, T + 1), L^2)} \lesssim \epsilon.
\end{equation}
The triangular inequality gives
\begin{align}
\left\Vert (1 - a \chi) g \right\Vert_{L^1((T, T + 1), L^2)} \lesssim & \ \left\Vert (1 - a \chi ) u_{\Lin} \right\Vert_{L^1((T, T + 1), L^2)} + \left\Vert (1 - a \chi) \partial_t u_{\Lin} \right\Vert_{L^1((T, T + 1), L^2)} \nonumber \\
& + \left\Vert (u_{\NL}, \partial_t u_{\NL}) - ( u_{\Lin}, \partial_t u_{\Lin}) \right\Vert_{L^\infty((T, T + 1), H_0^1 \times L^2)} \nonumber \\
& + \left\Vert f(u_{\NL}) \right\Vert_{L^1((T, T + 1), L^2)} + \left\Vert f(\varphi u_{\NL}) \right\Vert_{L^1((T, T + 1), L^2)}. \label{eq_proof_thm_null_control_scattering_7}
\end{align}
Using Lemma \ref{lem_basic_estimate_f_scatt}-\emph{(i)}, together with (\ref{eq_proof_thm_null_control_scattering_3}), one finds
\begin{align}
& \left\Vert f(u_{\NL}) \right\Vert_{L^1((T, T + 1), L^2)} + \left\Vert f(\varphi u_{\NL}) \right\Vert_{L^1((T, T + 1), L^2)} \nonumber \\
\lesssim \ & \left\Vert u_{\NL} \right\Vert_{L^{\alpha_0}((T, T + 1), L^{2 \alpha_0})}^{\alpha_0} + \left\Vert u_{\NL} \right\Vert_{L^{\alpha_1}((T, T + 1), L^{2 \alpha_1})}^{\alpha_1} \nonumber \\
\lesssim \ & \epsilon. \label{eq_proof_thm_null_control_scattering_8}
\end{align}
Coming back to (\ref{eq_proof_thm_null_control_scattering_7}), and using (\ref{eq_proof_thm_null_control_scattering_1}), (\ref{eq_proof_thm_null_control_scattering_2}) and (\ref{eq_proof_thm_null_control_scattering_8}), one obtains (\ref{eq_proof_thm_null_control_scattering_6}).

To complete the proof, we show that (\ref{eq_proof_thm_null_control_scattering_6}) implies (\ref{eq_proof_thm_null_control_scattering_5}). For $\tau \in [0, 1]$, set $Y_\tau = Y_{[T, T + \tau]}$, where $Y_{[T, T + \tau]}$ is defined by (\ref{eq_def_Y_T}). Note that $h$ is the solution of 
\[\left \{
\begin{array}{rcccl}
\square h + \beta h & = & f^\prime(v) h + \NL_v(h) + (a \chi - 1)g & \quad & \text{in } (T, T + 1) \times \Omega, \\
(h(T), \partial_t h(T)) & = & 0 & \quad & \text{in } \Omega, \\
h & = & 0 & \quad & \text{on } (T, T + 1) \times \partial \Omega.
\end{array}
\right.\]
where $\NL_v(h)$ is defined by (\ref{eq_def_NL}). By Proposition \ref{prop_LKG_time_dependent_potential}-\emph{(ii)} and Remark \ref{rem_const_temps_petits}, there exists a constant independent of $\tau \in [0, 1]$ such that 
\[\left\Vert h \right\Vert_{Y_\tau} \lesssim \left\Vert \NL_v(h) + ( a \chi - 1 ) g \right\Vert_{L^1((T, T + \tau), L^2)}.\]
Using Lemma \ref{lem_basic_estimate_f_scatt}-\emph{(ii)} and (\ref{eq_proof_thm_null_control_scattering_6}), one obtains
\begin{equation}\label{eq_proof_thm_null_control_scattering_9}
\left\Vert h \right\Vert_{Y_\tau} \leq c \left( \epsilon + \left\Vert h \right\Vert_{Y_\tau}^2 + \left\Vert h \right\Vert_{Y_\tau}^{\alpha_1} \right),
\end{equation}
for some $c > 0$ independent of $\epsilon$ and of $\tau \in [0, 1]$. Set $\theta(s) = c \left( s^2 + s^{\alpha_1} \right) - s$, for $s \geq 0$. If $\epsilon$ is sufficiently small, then by the mean value theorem, either $\theta^\prime \left( \left\Vert h \right\Vert_{Y_\tau} \right) \leq 0$ for all $\tau \in [0, 1]$, or $\theta^\prime \left( \left\Vert h \right\Vert_{Y_\tau} \right) \geq 0$ for all $\tau \in [0, 1]$. As 
\[\left\Vert h \right\Vert_{Y_\tau} \xrightarrow{\tau \rightarrow 0^+} 0,\]
this latter case cannot occur. Hence, one has 
\[\left\Vert h \right\Vert_{Y_\tau} \geq c \left( 2 \left\Vert h \right\Vert_{Y_\tau}^2 + \alpha_1 \left\Vert h \right\Vert_{Y_\tau}^{\alpha_1} \right) \geq 2 c \left( \left\Vert h \right\Vert_{Y_\tau}^2 + \left\Vert h \right\Vert_{Y_\tau}^{\alpha_1} \right),\]
for $\tau \in [0, 1]$. In particular, (\ref{eq_proof_thm_null_control_scattering_9}) gives
\[\left\Vert h \right\Vert_{Y_1} \leq 2 c \epsilon,\]
yielding (\ref{eq_proof_thm_null_control_scattering_5}). This completes the proof.
\end{proof}

\appendix

\section{Statement of the Christ-Kiselev lemma}\label{appendix_christ_kiselev}

We recall the statement of the Christ-Kiselev lemma (see \cite{CHRIST-Kiselev}).

\begin{lem}\label{lem_christ_kiselev}
Let $X$ and $Y$ be Banach spaces. Consider $1 \leq p < q \leq \infty$, and let $T: L^p(\mathbb{R}, X) \rightarrow  L^q(\mathbb{R}, Y)$ be a continuous linear operator. Then, the operator 
\[\begin{array}{cccc} 
\tilde{T}: & L^p(\mathbb{R}, X) & \longrightarrow & L^q(\mathbb{R}, Y) \\
& F & \longmapsto & \left( t \longmapsto T \left( \mathds{1}_{(- \infty, t)} F \right) \right)
\end{array}\]
is well-defined and continuous.
\end{lem}

\section{A general lemma}\label{appendix_lemme_compact_term}

\begin{lem}\label{lem_general_compact_term_2}
Let $\mathcal{X}$ be a Banach space, and $\mathcal{Y}$ and $\mathcal{Z}$ be normed spaces. Let $A: \mathcal{X} \rightarrow \mathcal{Y}$ and $K: \mathcal{X} \rightarrow \mathcal{Z}$ be continuous linear operators. Assume that $K$ is compact, $A$ is one-to-one, and that there exists a constant $C > 0$ such that for all $x \in \mathcal{X}$, one has
\[\Vert x \Vert_\mathcal{X} \leq C \left\Vert A x \right\Vert_\mathcal{Y} + \left\Vert K x \right\Vert_\mathcal{Z}.\]
Then, there exists a constant $C^\prime > 0$ such that for all $x \in \mathcal{X}$, one has
\[\Vert x \Vert_\mathcal{X} \leq C^\prime \left\Vert A x \right\Vert_\mathcal{Y}.\]
\end{lem}

The proof of Lemma \ref{lem_general_compact_term_2} is standard and we omit it. 

\section{An extension of Rellich's theorem}\label{appendix_proof_rellich}

Here, we prove Lemma \ref{lem_rellich_extension}. Consider $s \in \mathbb{R}$, $U$ a (possibly empty) smooth bounded open subset of $\mathbb{R}^d$, $V \in \mathscr{C}^\infty(\mathbb{R}^d \backslash U)$ such that 
\[\sum_{\vert \beta \vert \leq \vert s - 1 \vert} \left\vert \partial_x^\beta V(x) \right\vert \xrightarrow{|x| \rightarrow \infty} 0,\]
and $\chi \in \mathscr{C}^\infty(\mathbb{R}^d, [0, 1])$ such that $\chi = 1$ on $B(0, 1)$ and $\chi = 0$ on $\mathbb{R}^d \backslash B(0, 2)$. Write $\Omega = \mathbb{R}^d \backslash U$. 

Let $(u_n)_n$ be a bounded sequence in $H^s(\Omega)$. For all $k \in \mathbb{N}$ sufficiently large, by the usual Rellich theorem, there exists a subsequence of $\left( \chi\left( \frac{\cdot}{k} \right) u_n \right)_n$ converging in $H^{s - 1}(\Omega \cap B(0, 2k))$. Using a diagonal argument, one proves that up to a subsequence, there exists $u_\infty \in H^{s - 1}_\mathrm{loc}(\Omega)$ such that 
\begin{equation}\label{eq_proof_lem_rellich_extension_1}
\chi\left( \frac{\cdot}{k} \right) u_n \xrightarrow{n \rightarrow \infty} u_\infty
\end{equation}
in $H^{s - 1}_\mathrm{loc}(\Omega)$, for all $k \in \mathbb{N}$ sufficiently large. 

We show that $(V u_n)_n$ is a Cauchy sequence in $H^{s - 1}(\Omega)$. Consider $\epsilon > 0$, and let $k \in \mathbb{N}$ be such that 
\begin{equation}\label{eq_proof_lem_rellich_extension_2}
\sum_{\vert \beta \vert \leq \vert s - 1 \vert} \left\vert \partial_x^\beta V(x) \right\vert \leq \epsilon
\end{equation}
for all $x \in \Omega \backslash B(0, k)$. Write
\[\left\Vert V u_n - V u_m \right\Vert_{H^{s - 1}(\Omega)} \leq \left\Vert \chi\left( \frac{\cdot}{k} \right) \left( V u_n - V u_m \right) \right\Vert_{H^{s - 1}(\Omega)} + \left\Vert \left( 1 - \chi\left( \frac{\cdot}{k} \right) \right) \left( V u_n - V u_m \right) \right\Vert_{H^{s - 1}(\Omega)}.\]
By (\ref{eq_proof_lem_rellich_extension_1}), one has
\[\left\Vert \chi\left( \frac{\cdot}{k} \right) \left( V u_n - V u_m \right) \right\Vert_{H^{s - 1}(\Omega)} \lesssim \left\Vert \chi\left( \frac{\cdot}{k} \right) \left( u_n - u_m \right) \right\Vert_{H^{s - 1}(\Omega)} \leq \epsilon\]
for $n$ and $m$ sufficiently large. For $\phi \in H^{s - 1}(\Omega)$, one can prove that 
\[\Vert V \phi \Vert_{H^{s - 1}(\Omega)} \lesssim \sum_{\vert \beta \vert \leq \vert s - 1 \vert} \left\Vert \partial_x^\beta V \right\Vert_{L^\infty(\Omega)} \Vert \phi \Vert_{H^{s - 1}(\Omega)}.\]
As $(u_n)_n$ is bounded in $H^{s - 1}(\Omega)$, this implies
\[\left\Vert \left( 1 - \chi\left( \frac{\cdot}{k} \right) \right) \left( V u_n - V u_m \right) \right\Vert_{H^{s - 1}(\Omega)} \lesssim \sum_{\vert \beta \vert \leq \vert s - 1 \vert} \left\Vert \partial_x^\beta \left( \left( 1 - \chi\left( \frac{\cdot}{k} \right) \right) V \right) \right\Vert_{L^\infty(\Omega)}.\]
As $\chi\left( \frac{x}{k} \right) = 1$ for $x \in B(0, k)$, (\ref{eq_proof_lem_rellich_extension_2}) gives
\[\left\Vert \left( 1 - \chi\left( \frac{\cdot}{k} \right) \right) \left( V u_n - V u_m \right) \right\Vert_{H^{s - 1}(\Omega)} \lesssim \epsilon,\]
implying that $(V u_n)_n$ is a Cauchy sequence in $H^{s - 1}(\Omega)$. This completes the proof of Lemma \ref{lem_rellich_extension}.

\section{Proof of Lemma \ref{lem_basic_estimate_f_local} and of Lemma \ref{lem_basic_estimate_f_scatt}}\label{appendix_proof_lem_Holder}

The following consequence of Hölder's inequality will be repeatedly used in the proofs below:
\begin{equation}\label{proof_lem_basic_f_eq_1_bis}
\left\Vert \vert u \vert \vert v \vert^{\alpha - 1} \right\Vert_{L^1((0, T), L^2)} \leq \left\Vert u \right\Vert_{L^\alpha((0, T), L^{2 \alpha})} \left\Vert v \right\Vert_{L^\alpha((0, T), L^{2 \alpha})}^{\alpha - 1}, \quad u, v \in L^\alpha((0, T), L^{2 \alpha}(\Omega)).
\end{equation}
Note that one must be careful in the case $\alpha \in (1, 2)$, as the exponent appearing in (\ref{eq_def_nonlinearity_local}) is negative in that case.

\begin{proof}[Proof of Lemma \ref{lem_basic_estimate_f_local}.]
We prove \emph{(i)}, \emph{(ii)} and \emph{(iii)} separately.

\underline{Proof of \emph{(i)}.}
Integrating (\ref{eq_def_nonlinearity_local}), one finds
\begin{equation}
    \left\vert f^\prime(s) \right\vert \lesssim \int_0^{\vert s \vert} \left( 1 + \vert \tau \vert \right)^{\alpha - 2} \mathrm{d} \tau = \frac{\left( 1 + \vert s \vert \right)^{\alpha - 1} - 1}{\alpha - 1} \lesssim 1 + \vert s \vert^{\alpha - 1}, \quad s \in \mathbb{R},\label{proof_lem_basic_f_eq_0}
\end{equation}
yielding
\[\left\vert f(u) - f(v) \right\vert \lesssim \vert u - v \vert \left( 1 + \vert u \vert^{\alpha - 1} + \vert v \vert^{\alpha - 1} \right), \quad u, v \in \mathbb{R}.\]
Together with (\ref{proof_lem_basic_f_eq_1_bis}), this gives
\begin{align*}
\left\Vert f(u) - f(v) \right\Vert_{L^1((0, T), L^2)} \lesssim \left\Vert u - v \right\Vert_{L^1((0, T), L^2)} & + \left\Vert u - v \right\Vert_{L^\alpha((0, T), L^{2 \alpha})} \left\Vert u \right\Vert_{L^\alpha((0, T), L^{2 \alpha})}^{\alpha - 1} \\
& + \left\Vert u - v \right\Vert_{L^\alpha((0, T), L^{2 \alpha})} \left\Vert v \right\Vert_{L^\alpha((0, T), L^{2 \alpha})}^{\alpha - 1}, \quad u, v \in X_T,
\end{align*}
implying \emph{(i)}. 

\underline{Proof of \emph{(ii)} in the case $\alpha \geq 2$.}
For $s, h_1, h_2 \in \mathbb{R}$ such that $h_2 \leq h_1$, one has
\begin{align}
\left\vert f(s + h_1) - f(s + h_2) \right. & \left. - f^\prime(s) (h_1 - h_2) \right\vert = \left\vert \int_{s + h_2}^{s + h_1} \int_s^t f^{\prime \prime}(\tau) \mathrm{d}\tau \mathrm{d}t \right\vert \nonumber \\
\lesssim & \left( 1 + \left\vert s + h_1 \right\vert^{\alpha - 2} + \left\vert s + h_2 \right\vert^{\alpha - 2} \right) \int_{s + h_2}^{s + h_1} \vert t - s \vert \mathrm{d}t \nonumber \\
\lesssim & \left( 1 + \left\vert s + h_1 \right\vert^{\alpha - 2} + \left\vert s + h_2 \right\vert^{\alpha - 2} \right) \vert h_1 - h_2 \vert \left( \vert h_1 \vert + \vert h_2 \vert \right), \label{proof_lem_basic_f_eq_1}
\end{align}
and (\ref{proof_lem_basic_f_eq_1}) also holds true if $h_1 \leq h_2$. It implies
\begin{align}
& \left\Vert \NL_\mathbf{u}(u) - \NL_\mathbf{u}(v) \right\Vert_{L^1((0, T), L^2)} \nonumber \\
\lesssim \ & \left\Vert (u - v) u \right\Vert_{L^1((0, T), L^2)} + \left\Vert (u - v) v \right\Vert_{L^1((0, T), L^2)} \nonumber \\
& + \left\Vert (u - v) u \vert \mathbf{u} + u \vert^{\alpha - 2} \right\Vert_{L^1((0, T), L^2)} + \left\Vert (u - v) u \vert \mathbf{u} + v \vert^{\alpha - 2} \right\Vert_{L^1((0, T), L^2)} \nonumber \\
& + \left\Vert (u - v) v \vert \mathbf{u} + u \vert^{\alpha - 2} \right\Vert_{L^1((0, T), L^2)} + \left\Vert (u - v) v \vert \mathbf{u} + v \vert^{\alpha - 2} \right\Vert_{L^1((0, T), L^2)} \label{proof_lem_basic_f_eq_2}
\end{align}
for $u, v \in X_T$.

One the one hand, Hölder's inequality gives
\[\left\Vert (u - v) u \right\Vert_{L^1((0, T), L^2)} \leq \left\Vert u - v \right\Vert_{L^\alpha((0, T), L^{2 \alpha})} \left\Vert u \right\Vert_{L^{\alpha^\prime}((0, T), L^{2 \alpha^\prime})}.\]
One has $1 \leq \alpha^\prime \leq \alpha$, implying
\[\left\Vert u \right\Vert_{L^{\alpha^\prime}((0, T), L^{2 \alpha^\prime})} \leq \left\Vert u \right\Vert_{L^1((0, T), L^2)}^{\theta_1} \left\Vert u \right\Vert_{L^\alpha((0, T), L^{2 \alpha})}^{1 - \theta_1},\]
where $\theta_1 = \frac{\alpha - 2}{\alpha - 1}$ is given by $\frac{1}{\alpha^\prime} = \theta_1 + \frac{1 - \theta_1}{\alpha}$. This gives
\begin{equation}\label{proof_lem_basic_f_eq_3}
\left\Vert (u - v) u \right\Vert_{L^1((0, T), L^2)} \lesssim \left\Vert u - v \right\Vert_{X_T} \left\Vert u \right\Vert_{X_T}.
\end{equation}

On the other hand, as above, one has
\begin{equation}\label{proof_lem_basic_f_eq_4}
\left\Vert (u - v) u \vert \mathbf{u} + v \vert^{\alpha - 2} \right\Vert_{L^1((0, T), L^2)} \leq \left\Vert u - v \right\Vert_{L^\alpha((0, T), L^{2 \alpha})} \left\Vert u \vert \mathbf{u} + v \vert^{\alpha - 2} \right\Vert_{L^{\alpha^\prime}((0, T), L^{2 \alpha^\prime})}.
\end{equation}
Note that $\theta_2 = \frac{1}{\alpha - 1}$ satisfies $\frac{2 \alpha^\prime (\alpha - 2)}{1 - \theta_2} = \frac{2 \alpha^\prime}{\theta_2} = 2 \alpha$. Hence, applying Hölder's inequality with $1 = \frac{1}{1/\theta_2} + \frac{1}{1/(1 - \theta_2)}$, one obtains
\begin{equation}\label{proof_lem_basic_f_eq_5}
\left\Vert u \vert \mathbf{u} + v \vert^{\alpha - 2} \right\Vert_{L^{\alpha^\prime}((0, T), L^{2 \alpha^\prime})} \leq \left\Vert u \right\Vert_{L^\alpha((0, T), L^{2 \alpha})} \left\Vert \mathbf{u} + v \right\Vert_{L^\alpha((0, T), L^{2 \alpha})}^{\alpha - 2}.
\end{equation}
Using $\mathbf{u} \in L^\alpha((0, T), L^{2 \alpha}(\Omega))$, (\ref{proof_lem_basic_f_eq_4}) and (\ref{proof_lem_basic_f_eq_5}), one finds
\begin{equation}\label{proof_lem_basic_f_eq_6}
\left\Vert (u - v) u \vert \mathbf{u} + v \vert^{\alpha - 2} \right\Vert_{L^1((0, T), L^2)} \lesssim \left\Vert u - v \right\Vert_{L^\alpha((0, T), L^{2 \alpha})} \left\Vert u \right\Vert_{L^\alpha((0, T), L^{2 \alpha})} \left( 1 + \left\Vert v \right\Vert_{L^\alpha((0, T), L^{2 \alpha})}^{\alpha - 2}\right).
\end{equation}

Combining (\ref{proof_lem_basic_f_eq_2}), (\ref{proof_lem_basic_f_eq_3}) and (\ref{proof_lem_basic_f_eq_6}), one obtains
\begin{align*}
& \left\Vert \NL_\mathbf{u}(u) - \NL_\mathbf{u}(v) \right\Vert_{L^1((0, T), L^2)} \\
\lesssim \ & \left\Vert u - v \right\Vert_{X_T} \left( \left\Vert u \right\Vert_{X_T} + \left\Vert v \right\Vert_{X_T} + \left\Vert u \right\Vert_{X_T} \left\Vert v \right\Vert_{X_T}^{\alpha - 2} + \left\Vert v \right\Vert_{X_T} \left\Vert u \right\Vert_{X_T}^{\alpha - 2} + \left\Vert v \right\Vert_{X_T}^{\alpha - 1} + \left\Vert v \right\Vert_{X_T}^{\alpha - 1} \right),
\end{align*}
implying \emph{(ii)}.

\underline{Proof of \emph{(ii)} in the case $\alpha \in (1, 2)$.} 
Write $x \wedge y = \min(x, y)$ and $x \vee y = \max(x, y)$, for $(x, y) \in \mathbb{R}^2$. For $s, h_1, h_2 \in \mathbb{R}$, one has
\begin{align}
\left\vert f(s + h_2) - f(s + h_1) - f^\prime(s) (h_2 - h_1) \right\vert & = \left\vert \int_{s + h_1}^{s + h_2} \int_s^t f^{\prime \prime}(\tau) \mathrm{d}\tau \mathrm{d}t \right\vert \nonumber \\
& \lesssim \int_{s + h_1 \wedge h_2}^{s + h_1 \vee h_2} \int_{s \wedge t}^{s \vee t} \left( 1 + \vert \tau \vert \right)^{\alpha - 2} \mathrm{d}\tau \mathrm{d}t. \label{eq_techniq_alpha_entre_1_et_2__1}
\end{align}
We claim that 
\begin{equation}\label{eq_techniq_alpha_entre_1_et_2__2}
\int_{s + h_1 \wedge h_2}^{s + h_1 \vee h_2} \int_{s \wedge t}^{s \vee t} \left( 1 + \vert \tau \vert \right)^{\alpha - 2} \mathrm{d}\tau \mathrm{d}t \leq \frac{2^{2 - \alpha}}{\alpha - 1} \left\vert h_1 - h_2 \right\vert \left( \left\vert h_1 \right\vert^{\alpha - 1} + \left\vert h_2 \right\vert^{\alpha - 1} \right).
\end{equation}
Using (\ref{eq_techniq_alpha_entre_1_et_2__2}), one can complete the proof of \emph{(ii)}. Indeed, Hölder's inequality, (\ref{eq_techniq_alpha_entre_1_et_2__1}) and (\ref{eq_techniq_alpha_entre_1_et_2__2}) yield
\begin{align*}
& \left\Vert \NL_\mathbf{u}(u) - \NL_\mathbf{u}(v) \right\Vert_{L^1((0, T), L^2)} \\
\lesssim & \left( \left\Vert \left( u - v \right) \left\vert u \right\vert^{\alpha - 1} \right\Vert_{L^1((0, T), L^2)} + \left\Vert \left( u - v \right) \left\vert v \right\vert^{\alpha - 1} \right\Vert_{L^1((0, T), L^2)} \right) \\
\lesssim & \left\Vert u - v \right\Vert_{L^\alpha((0, T), L^{2 \alpha})} \left( \left\Vert u \right\Vert_{L^\alpha((0, T), L^{2 \alpha})}^{\alpha - 1} + \left\Vert v \right\Vert_{L^\alpha((0, T), L^{2 \alpha})}^{\alpha - 1} \right).
\end{align*}

Now, we prove (\ref{eq_techniq_alpha_entre_1_et_2__2}). Note that we can assume $h_1 < h_2$. One has 
\begin{align}
\int_{s + h_1}^{s + h_2} \int_{s \wedge t}^{s \vee t} & \left( 1 + \vert \tau \vert \right)^{\alpha - 2} \mathrm{d}\tau \mathrm{d}t = \int_\mathbb{R} \left( 1 + \vert \tau \vert \right)^{\alpha - 2} \int_{s + h_1}^{s + h_2} \mathds{1}_{s \wedge t \leq \tau \leq s \vee t} \mathrm{d} t \mathrm{d} \tau \nonumber \\
= & \int_\mathbb{R} \left( 1 + \vert \tau \vert \right)^{\alpha - 2} \int_{s + h_1}^{s + h_2} \left( \mathds{1}_{s \leq t} \mathds{1}_{s \leq \tau \leq t} + \mathds{1}_{s > t} \mathds{1}_{t \leq \tau \leq s} \right) \mathrm{d} t \mathrm{d} \tau \nonumber \\
\leq & \left( h_2 - h_1 \right) \int_\mathbb{R} \left( 1 + \vert \tau \vert \right)^{\alpha - 2} \left( \mathds{1}_{s \leq \tau \leq s + h_2} + \mathds{1}_{s + h_1 \leq \tau \leq s} \right) \mathrm{d} \tau \nonumber \\
= & \left( h_2 - h_1 \right) \left( \mathds{1}_{h_2 > 0} \int_s^{s + h_2} \left( 1 + \vert \tau \vert \right)^{\alpha - 2} \mathrm{d} \tau + \mathds{1}_{h_1 < 0} \int_{s + h_1}^s \left( 1 + \vert \tau \vert \right)^{\alpha - 2} \mathrm{d} \tau \right) \nonumber \\
= & \left( h_2 - h_1 \right) \left( I(s, h_2) + I(-s, -h_1) \right), \label{eq_proof_lem_alpha_small_3}
\end{align}
where $I$ is defined by
\[I(s, h) = \mathds{1}_{h > 0} \int_s^{s + h} \left( 1 + \vert \tau \vert \right)^{\alpha - 2} \mathrm{d} \tau, \quad s, h \in \mathbb{R}.\]
A change of variable gives
\begin{equation}\label{eq_proof_lem_alpha_small_4}
I(s, h) = h^{\alpha - 1} \int_0^1 \left( \frac{1}{h} + \left\vert \frac{s}{h} + \tau \right\vert \right)^{\alpha - 2} \mathrm{d} \tau, \quad s \in \mathbb{R}, h > 0.
\end{equation}
We use the following elementary lemma.

\begin{lem}\label{lem_technique_elementary_small_alpha}
For $1 < \alpha \leq 2$, $a > 0$ and $b \in \mathbb{R}$, one has
\[\int_0^1 \left( a + \left\vert \tau - b \right\vert \right)^{\alpha - 2} \mathrm{d} \tau \leq \frac{2^{2 - \alpha}}{\alpha - 1}.\]
\end{lem}

A proof of Lemma \ref{lem_technique_elementary_small_alpha} is given below. Together with (\ref{eq_proof_lem_alpha_small_3}) and (\ref{eq_proof_lem_alpha_small_4}), it implies (\ref{eq_techniq_alpha_entre_1_et_2__2}).

\begin{proof}[Proof of Lemma \ref{lem_technique_elementary_small_alpha}.]
Firstly, if $b \leq 0$, then 
\[\int_0^1 \left( a + \left\vert \tau - b \right\vert \right)^{\alpha - 2} \mathrm{d} \tau \leq \int_0^1 \tau^{\alpha - 2} \mathrm{d} \tau  = \frac{1}{\alpha - 1}.\]
Secondly, if $b \geq 1$, then 
\[\int_0^1 \left( a + \left\vert \tau - b \right\vert \right)^{\alpha - 2} \mathrm{d} \tau \leq \int_0^1 \left( 1 - \tau \right)^{\alpha - 2} \mathrm{d} \tau  = \frac{1}{\alpha - 1}.\]
Thirdly, if $0 < b < 1$, then
\begin{align*}
\int_0^1 \left( a + \left\vert \tau - b \right\vert \right)^{\alpha - 2} \mathrm{d} \tau & \leq \int_0^b \left( b - \tau \right)^{\alpha - 2} \mathrm{d} \tau + \int_b^1 \left( \tau - b \right)^{\alpha - 2} \mathrm{d} \tau \\
& = \frac{b^{\alpha - 1} + \left( 1 - b \right)^{\alpha - 1}}{\alpha - 1} \leq \frac{2^{2 - \alpha}}{\alpha - 1},
\end{align*}
and this completes the proof.
\end{proof}

\underline{Proof of \emph{(iii)}.}
Using (\ref{proof_lem_basic_f_eq_0}), one finds
\[\left\Vert f^\prime(\mathbf{u}) u \right\Vert_{L^1((0, T), L^2)} \lesssim \left\Vert u \right\Vert_{L^1((0, T), L^2)} + \left\Vert \left\vert \mathbf{u} \right\vert^{\alpha - 1} u \right\Vert_{L^1((0, T), L^2)},\]
for all $u \in X_T$. As above, using Hölder's inequality, this gives \emph{(iii)}. This completes the proof of Lemma \ref{lem_basic_estimate_f_local}.
\end{proof}

\begin{proof}[Proof of Lemma \ref{lem_basic_estimate_f_scatt}]
The proof is similar to that of Lemma \ref{lem_basic_estimate_f_local}, so we only sketch it. Using
\[\left\vert f(s_1) - f(s_2) \right\vert \lesssim \vert s_1 - s_2 \vert \left( \vert s_1 \vert^{\alpha_0 - 1} + \vert s_2 \vert^{\alpha_0 - 1} + \vert s_1 \vert^{\alpha_1 - 1} + \vert s_2 \vert^{\alpha_1 - 1} \right), \quad s_1, s_2 \in \mathbb{R},\]
one obtains
\[\left\Vert f(u) - f(v) \right\Vert_{L^1((0, T), L^2)} \lesssim \sum_{i = 0,1} \left( \left\Vert \vert u - v \vert \vert u \vert^{\alpha_i - 1} \right\Vert_{L^1((0, T), L^2)} + \left\Vert \vert u - v \vert \vert v \vert^{\alpha_i - 1} \right\Vert_{L^1((0, T), L^2)} \right),\]
and this yields \emph{(i)}, by (\ref{proof_lem_basic_f_eq_1_bis}).

Consider $u, v \in L^{\alpha_0}((0, T), L^{2 \alpha_0}(\Omega)) \cap L^{\alpha_1}((0, T), L^{2 \alpha_1}(\Omega))$. For $s, h_1, h_2 \in \mathbb{R}$, one has
\begin{align*}
& \left\vert f(s + h_1) - f(s + h_2) - f^\prime(s) (h_1 - h_2) \right\vert \\
\lesssim & \left( \left\vert s + h_1 \right\vert^{\alpha_0 - 2} + \left\vert s + h_2 \right\vert^{\alpha_0 - 2} + \left\vert s + h_1 \right\vert^{\alpha_1 - 2} + \left\vert s + h_2 \right\vert^{\alpha_1 - 2} \right) \vert h_1 - h_2 \vert \left( \vert h_1 \vert + \vert h_2 \vert \right), 
\end{align*}
implying
\begin{align*}
& \left\Vert \NL_\mathbf{u}(u) - \NL_\mathbf{u}(v) \right\Vert_{L^1((0, T), L^2)} \\
\lesssim \ & \sum_{i = 0,1} \left( \left\Vert (u - v) u \vert \mathbf{u} + u \vert^{\alpha_i - 2} \right\Vert_{L^1((0, T), L^2)} + \left\Vert (u - v) u \vert \mathbf{u} + v \vert^{\alpha_i - 2} \right\Vert_{L^1((0, T), L^2)} \right. \\
& \hspace{2.55em} \left. \left\Vert (u - v) v \vert \mathbf{u} + u \vert^{\alpha_i - 2} \right\Vert_{L^1((0, T), L^2)} + \left\Vert (u - v) v \vert \mathbf{u} + v \vert^{\alpha_i - 2} \right\Vert_{L^1((0, T), L^2)} \right).
\end{align*}
As in the proof of Lemma \ref{lem_basic_estimate_f_local}, one has
\[\left\Vert (u - v) u \vert \mathbf{u} + v \vert^{\alpha_i - 2} \right\Vert_{L^1((0, T), L^2)} \lesssim \left\Vert u - v \right\Vert_{L^{\alpha_i}((0, T), L^{2 \alpha_i})} \left\Vert u \right\Vert_{L^{\alpha_i}((0, T), L^{2 \alpha_i})} \left( 1 + \left\Vert v \right\Vert_{L^{\alpha_i}((0, T), L^{2 \alpha_i})}^{\alpha_i - 2}\right),\]
and this gives the inequality of \emph{(ii)}. The constant can be expressed as a nondecreasing function of $\left\Vert \mathbf{u} \right\Vert_{L^{\alpha_i}((0, T), L^{2 \alpha_0})} + \left\Vert \mathbf{u} \right\Vert_{L^{\alpha_i}((0, T), L^{2 \alpha_1})}$, implying the remark of \emph{(ii)} about the constant. 

Finally, one has
\begin{align*}
\left\Vert f^\prime(\mathbf{u}) u \right\Vert_{L^1((0, T), L^2)} \lesssim \left\Vert u \vert \mathbf{u} \vert^{\alpha_0 - 1} \right\Vert_{L^1((0, T), L^2)} + \left\Vert u \vert \mathbf{u} \vert^{\alpha_1 - 1} \right\Vert_{L^1((0, T), L^2)},
\end{align*}
and as above, it implies \emph{(iii)}. 
\end{proof}

\section{Observability on an unbounded domain}\label{appendix_proof_obs_unbounded}



Here, we prove Theorem \ref{thm_obs_constant_coeff}, in the case of an unbounded domain $\Omega$. The proof is decomposed into two steps.

\paragraph{Step 1:} We start with the case $\Omega = \mathbb{R}^d$, with the Euclidean metric, and $a = 1$, in which Theorem \ref{thm_obs_constant_coeff} can be proved by a direct Fourier computation. Consider $\left( u^0, u^1 \right) \in L^2(\mathbb{R}^d) \times H^{-1}(\mathbb{R}^d)$ and write $u$ for the solution of
\[\left \{
\begin{array}{rcccl}
\square u + \beta u & = & 0 & \quad & \text{in } (0, T) \times \mathbb{R}^d, \\
(u(0), \partial_t u(0)) & = & \left( u^0, u^1 \right) & \quad & \text{in } \mathbb{R}^d. 
\end{array}
\right.\]
Write $\widehat{u}$ for the Fourier transform of $u$ with respect to the space variable $x$. There exist some complex functions $A$ and $B$ such that
\[\widehat{u}(t, \xi) = A(\xi) e^{i \langle \xi \rangle t} + B(\xi) e^{- i \langle \xi \rangle t}\]
for $t \in [0, T]$ and $\xi \in \mathbb{R}^d$, where $\langle \xi \rangle^2 = 1 + \vert \xi \vert^2$. On the one hand, one has
\begin{align*}
\left\Vert \left( u^0, u^1 \right) \right\Vert_{L^2(\mathbb{R}^d) \times H^{-1}(\mathbb{R}^d)}^2 \lesssim & \int_{\mathbb{R}^d} \left( \left\vert \widehat{u^0}(\xi) \right\vert^2 + \langle \xi \rangle^{-2} \left\vert \widehat{u^1}(\xi) \right\vert^2 \right) \mathrm{d}\xi \\
= & \int_{\mathbb{R}^d} \left( \left\vert A(\xi) + B(\xi) \right\vert^2 + \left\vert A(\xi) - B(\xi) \right\vert^2 \right) \mathrm{d}\xi \\
= & \ 2 \int_{\mathbb{R}^d} \left( \left\vert A(\xi) \right\vert^2 + \left\vert B(\xi) \right\vert^2 \right) \mathrm{d}\xi.
\end{align*}
On the other hand, a direct computation gives
\begin{align*}
\left\Vert u \right\Vert_{L^2((0, T) \times \mathbb{R}^d)}^2 \gtrsim & \int_{\mathbb{R}^d} \left( T \left\vert A(\xi) \right\vert^2 + T \left\vert B(\xi) \right\vert^2 + 2 \Re \left( A(\xi) \overline{B(\xi)} \frac{e^{2 i \langle \xi \rangle T} - 1}{2 i \langle \xi \rangle} \right) \right) \mathrm{d}\xi \\
\geq & \int_{\mathbb{R}^d} \left( \left\vert A(\xi) \right\vert^2 + \left\vert B(\xi) \right\vert^2 \right) \left(T - \frac{ \left\vert \sin \left( \langle \xi \rangle T \right) \right\vert}{\langle \xi \rangle} \right) \mathrm{d}\xi.
\end{align*}
For all $T > 0$, there exists a constant $C = C(T)$ such that for all $\alpha \geq 1$, one has
\[T - \frac{ \left\vert \sin \left( \alpha T \right) \right\vert}{\alpha} \geq C.\]
This completes the proof of Theorem \ref{thm_obs_constant_coeff} in the case $\Omega = \mathbb{R}^d$ and $a = 1$. Note that in that particular case, $T$ can be arbitrary small.

\paragraph{Step 2:}
Now, we prove Theorem \ref{thm_obs_constant_coeff} in the case of an unbounded domain $\Omega$. Write $U$ for the (possibly empty) smooth bounded open subset of $\mathbb{R}^d$ such that $\Omega = \mathbb{R}^d \backslash U$. Up to increasing $R_0$, we can assume that the metric of $\Omega \backslash B(0, R_0)$ is the Euclidean one. Let $\chi \in \mathscr{C}^\infty(\Omega, [0, 1])$ be such that $\chi = 1$ on $B(0, R_0 + T)$ and $\chi = 0$ on $\mathbb{R}^d \backslash B(0, R_0 + 2 T)$. Consider $\left( u^0, u^1 \right) \in L^2(\Omega) \times H^{-1}(\Omega)$ and write $u$ for the solution of 
\[\left \{
\begin{array}{rcccl}
\square u + \beta u & = & 0 & \quad & \text{in } (0, T) \times \Omega, \\
(u(0), \partial_t u(0)) & = & \left( u^0, u^1 \right) & \quad & \text{in } \Omega, \\
u & = & 0 & \quad & \text{on } (0, T) \times \partial \Omega.
\end{array}
\right.\]
One has $u = \chi u + (1 - \chi) u$, $\chi u = v + \tilde{v}$ and $(1 - \chi) u = w + \tilde{w}$, with
\[\left \{
\begin{array}{rcccl}
\square v + \beta u & = & - 2 \nabla \chi \nabla u - \Delta \chi u & \quad & \text{in } (0, T) \times \Omega, \\
(v(0), \partial_t v(0)) & = & 0 & \quad & \text{in } \Omega, \\
v & = & 0 & \quad & \text{on } (0, T) \times \partial \Omega,
\end{array}
\right.\]
\[\left \{
\begin{array}{rcccl}
\square \tilde{v} + \beta \tilde{v} & = & 0 & \quad & \text{in } (0, T) \times \Omega, \\
(\tilde{v}(0), \partial_t \tilde{v}(0)) & = & \left( \chi u^0, \chi u^1 \right) & \quad & \text{in } \Omega, \\
\tilde{v} & = & 0 & \quad & \text{on } (0, T) \times \partial \Omega,
\end{array}
\right.\]
\[\left \{
\begin{array}{rcccl}
\square w + \beta w & = & 2 \nabla \chi \nabla u + \Delta \chi u & \quad & \text{in } (0, T) \times \Omega, \\
(w(0), \partial_t w(0)) & = & 0 & \quad & \text{in } \Omega, \\
w & = & 0 & \quad & \text{on } (0, T) \times \partial \Omega,
\end{array}
\right.\]
\[\left \{
\begin{array}{rcccl}
\square \tilde{w} + \beta \tilde{w} & = & 0 & \quad & \text{in } (0, T) \times \Omega, \\
(\tilde{w}(0), \partial_t \tilde{w}(0)) & = & \left( (1 - \chi) u^0, (1 - \chi) u^1 \right) & \quad & \text{in } \Omega, \\
\tilde{w} & = & 0 & \quad & \text{on } (0, T) \times \partial \Omega.
\end{array}
\right.\]

As $(1 - \chi) u^0$ and $(1 - \chi) u^1$ are supported in $\mathbb{R}^d \backslash B(0, R_0 + 2 T)$, $\tilde{w}$ is the solution of
\[\left \{
\begin{array}{rcccl}
\square \tilde{w} + \beta \tilde{w} & = & 0 & \quad & \text{in } (0, T) \times \mathbb{R}^d, \\
(\tilde{w}(0), \partial_t \tilde{w}(0)) & = & \left( (1 - \chi) u^0, (1 - \chi) u^1 \right) & \quad & \text{in } \mathbb{R}^d,
\end{array}
\right.\]
by finite speed of propagation. Hence, the case $\Omega = \mathbb{R}^d$ (with the Euclidean metric) and $a = 1$ treated above gives
\begin{align}
\left\Vert \left( (1 - \chi) u^0, (1 - \chi) u^1 \right) \right\Vert_{L^2(\Omega) \times H^{-1}(\Omega)} & = \left\Vert \left( (1 - \chi) u^0, (1 - \chi) u^1 \right) \right\Vert_{L^2(\mathbb{R}^d) \times H^{-1}(\mathbb{R}^d)} \nonumber \\
& \lesssim \left\Vert \tilde{w} \right\Vert_{L^2((0, T) \times \mathbb{R}^d)} \nonumber \\
& = \left\Vert a \tilde{w} \right\Vert_{L^2((0, T) \times \Omega)}. \label{proof_thm_obs_constant_coeff_eq_1}
\end{align}
By finite speed of propagation again, $\tilde{v}$ is the solution of 
\[\left \{
\begin{array}{rcccl}
\square \tilde{v} + \beta \tilde{v} & = & 0 & \quad & \text{in } (0, T) \times \left( \Omega \cap B(0, R_0 + 2T) \right), \\
(\tilde{v}(0), \partial_t \tilde{v}(0)) & = & \left( \chi u^0, \chi u^1 \right) & \quad & \text{in } \Omega \cap B(0, R_0 + 2T), \\
\tilde{v} & = & 0 & \quad & \text{on } (0, T) \times \partial\left( \Omega \cap B(0, R_0 + 2T) \right).
\end{array}
\right.\]
Hence, Theorem \ref{thm_obs_constant_coeff} in the case of a compact domain gives
\begin{equation}\label{proof_thm_obs_constant_coeff_eq_2}
\left\Vert \left( \chi u^0, \chi u^1 \right) \right\Vert_{L^2(\Omega) \times H^{-1}(\Omega)} \lesssim \left\Vert a \tilde{v} \right\Vert_{L^2((0, T) \times \Omega)}. 
\end{equation}

Using (\ref{proof_thm_obs_constant_coeff_eq_1}) and (\ref{proof_thm_obs_constant_coeff_eq_2}), together with the continuity estimate of Proposition \ref{prop_LKG_time_dependent_potential}-\emph{(iii)} (with $V = 0$), one obtains
\begin{align*}
\left\Vert \left( u^0, u^1 \right) \right\Vert_{L^2(\Omega) \times H^{-1}(\Omega)} & \leq \left\Vert \left( \chi u^0, \chi u^1 \right) \right\Vert_{L^2(\Omega) \times H^{-1}(\Omega)} + \left\Vert \left( (1 - \chi) u^0, (1 - \chi) u^1 \right) \right\Vert_{L^2(\Omega) \times H^{-1}(\Omega)} \\
& \lesssim \left\Vert a u \right\Vert_{L^2((0, T) \times \Omega)} + \left\Vert a v \right\Vert_{L^2((0, T) \times \Omega)} + \left\Vert a w \right\Vert_{L^2((0, T) \times \Omega)} \\
& \lesssim \left\Vert a u \right\Vert_{L^2((0, T) \times \Omega)} + \left\Vert 2 \nabla \chi \nabla u + \Delta \chi u \right\Vert_{L^1((0, T), H^{-1})}.
\end{align*}
Consider $\phi \in H_0^1(\Omega)$, with $\Vert \phi \Vert_{H_0^1(\Omega)} \leq 1$. As $\nabla \chi$ and $\Delta \chi$ are supported in $B(0, R_0 + 2T) \backslash B(0, R_0 + T)$, one has
\[\left\vert \left\langle 2 \nabla \chi \nabla u(t) + \Delta \chi u(t), \phi \right\rangle_{H^{- 1}(\Omega) \times H_0^1(\Omega)} \right\vert = \left\vert \left\langle a u(t), \phi \Delta \chi - 2 \Div (\phi \nabla \chi) \right\rangle_{L^2(\Omega)} \right\vert \lesssim \left\Vert a u(t) \right\Vert_{L^2(\Omega)}\]
for all $t \in (0, T)$. Hence, the Cauchy-Schwarz inequality gives
\[\left\Vert \left( u^0, u^1 \right) \right\Vert_{L^2(\Omega) \times H^{-1}(\Omega)} \lesssim \left\Vert a u \right\Vert_{L^2((0, T) \times \Omega)}\]
and this completes the proof.

\printbibliography

\noindent
\textsc{Perrin Thomas:} \texttt{thomas.perrin@ens-rennes.fr}

\end{document}